\documentclass[11pt,reqno]{amsart}
\usepackage{a4wide}
\usepackage{amsfonts,amsmath,amssymb,amscd}
\usepackage{amsthm}
\usepackage{amsmath,todonotes}
\usepackage{amscd}
\usepackage{cancel}
\usepackage{verbatim}
\usepackage{color}
\usepackage[all]{xy}
\usepackage[mathscr]{eucal}
\usepackage{mathrsfs}
\usepackage{fullpage}
\usepackage{enumitem}
\usepackage{hyperref}
\usepackage{bbm}
\usepackage{qtree}
\usepackage{mathrsfs}
\usepackage{bm}
\usepackage{bigints}
\usepackage{url}
\allowdisplaybreaks

\renewcommand{\d}{{\delta}}

\renewcommand{\r}{{\rho}}

\usepackage[headheight=110pt,top=0.6in, bottom=0.7in, left=0.6in, right=0.6in]{geometry}
\numberwithin{equation}{section}

\newtheorem{remark}{Remark}

\newtheorem{theorem}{Theorem}
\newtheorem{entry}{ Entry}

\newtheorem{lemma}[theorem]{Lemma}

\newtheorem{corollary}[theorem]{Corollary}

\newtheorem{proposition}[theorem]{Proposition}

\theoremstyle{definition}

\numberwithin{theorem}{section} 
\numberwithin{equation}{section}
\numberwithin{table}{section}

\makeatletter
\def\proof{\@ifnextchar[{\@oproof}{\@nproof}}
\def\@oproof[#1][#2]{\trivlist\item[\hskip\labelsep\textit{#2 Proof of\
		#1.}~]\ignorespaces}
\def\@nproof{\trivlist\item[\hskip\labelsep\textit{Proof.}~]\ignorespaces}

\begin{document}
 
\title[]{Trigonometric analogues of the identities associated with twisted sums of divisor functions
 }
  \author{Debika Banerjee}
\address{Debika Banerjee\\ Department of Mathematics\\
Indraprastha Institute of Information Technology IIIT, Delhi\\
Okhla, Phase III, New Delhi-110020, India.} 
\email{debika@iiitd.ac.in}

\author{Khyati Khurana}
\address{Khyati Khurana\\ Department of Mathematics\\
Indraprastha Institute of Information Technology IIIT, Delhi\\
Okhla, Phase III, New Delhi-110020, India} 
\email{khyatii@iiitd.ac.in}

  
 \maketitle
 \begin{abstract}
Inspired by two entries published in Ramanujan's lost notebook on page 355, B. C. Berndt, S. Kim, and A. Zaharescu presented Riesz sum identities for the twisted divisor sums. Subsequently, S. Kim derived analogous results by replacing twisted divisor sums with twisted sums of divisor functions. In a recent work, the authors of the present paper deduced Cohen-type identities and Vorono\"i summation formulas associated with these twisted sums of divisor functions. 
 The present paper aims to derive an equivalent version of the results offered in the previous article in terms of identities involving finite sums of trigonometric functions and the doubly infinite series.
As an application, the authors provide an identity for $r_6(n)$, which is analogous to Hardy's famous result where $r_6(n)$ denotes the number of representations of natural number $n$ as a sum of six squares. 
 
 \end{abstract}
 

		\section{Introduction}
  The lost notebook \cite{MR947735} of Ramanujan contains several beautiful identities. Some  are intimately connected with the famous {\it circle} and {\it divisor} problems. Among these results, on page 355 in his lost notebook, we encounter the following two important identities involving a finite trigonometric sum and a doubly infinite series of Bessel functions.
\begin{entry}\label{entry1}
If $0<\theta<1$ and $x>0$, then
\begin{align}\label{Entry1}
\sum_{n=1}^\infty F\left(\frac{x}{n} \right)\sin(2\pi n \theta)=\pi x \left(\frac{1}{2}-\theta\right)-\frac{\cot(\pi \theta)}{4}  +\frac{\sqrt{x}}{2} \sum_{m=1}^\infty \sum_{n=0}^\infty \left\{  \frac{J_1(4\pi \sqrt{m(n+\theta) x} )  }{ \sqrt{m(n+\theta)} }  - \frac{J_1(4\pi \sqrt{m(n+1-\theta) x} )  }{ \sqrt{m(n+1-\theta)} }   \right\} .
\end{align}
 
\end{entry}

\begin{entry}\label{entry2}
If $0<\theta<1$ and $x>0$, then
\begin{align}\label{Entry2}
\sum_{n=1}^\infty F\left(\frac{x}{n} \right)\cos(2\pi n \theta)= \frac{1}{4}-x \log(2\sin(\pi \theta))+\frac{\sqrt{x} }{2} \sum_{m=1}^\infty \sum_{n=0}^\infty \left\{  \frac{M_1(4\pi \sqrt{m(n+\theta) x} )  }{ \sqrt{m(n+\theta)} }  + \frac{M_1(4\pi \sqrt{m(n+1-\theta) x} )  }{ \sqrt{m(n+1-\theta)} }   \right\},
\end{align}
\end{entry}
where 
\begin{align*}
F(x)=\begin{cases}
\lfloor x \rfloor, \ \ \quad \text{if }  x \text{ is not an integer},\\\
x-\frac{1}{2},\ \ \quad \text{if } x \text{ is an integer;  }
\end{cases}  
\end{align*}
and $J_\nu$ is the Bessel function  of the first kind of order $\nu$ \cite[p.~ 40]{MR1349110}
\begin{align}\label{J_bessel}
  J_{\nu}(z):=\sum_{n=0}^\infty \frac{(-1)^n(\frac{1}{2}z)^{\nu+2n}}{n! \Gamma{(\nu+n+1)}}, \ z\in \mathbb{C},
  \end{align}
  and
\begin{align}\label{bessel}
    M_\nu(z)=-Y_\nu(z)-\frac{2}{\pi } K_\nu(z),
\end{align}
where $Y_\nu$ is the Bessel function of the second kind  of order $\nu$
 and $ K_\nu(z)$ denotes the modified Bessel function of order $\nu$ \cite[p.~64, 78]{MR1349110}, which are defined as the following 

  \begin{align}\label{Y_bessel}
  & Y_{\nu}(z):= \frac{J_{\nu}(z)\cos{\pi \nu}-J_{-\nu}(z)}{\sin{\pi \nu}}, \ z\in \mathbb{C}, \nu \notin \mathbb{Z}, 
  \end{align} 
   \begin{align}\label{Kbessel}
  & K_{\nu}(z):=\frac{\pi}{2}\frac{I_{-\nu}(z)-I_{\nu}(z)}{\sin{\pi \nu}}, \ z\in \mathbb{C}, \nu \notin \mathbb{Z},
  \end{align}
  with $I_{\nu}$ being the Bessel function of the imaginary argument of order $\nu$ \cite[p.~77]{MR1349110} given by  
  \begin{align}\label{Inu}
  & I_{\nu}(z):=\sum_{n=0}^\infty \frac{(\frac{1}{2}z)^{\nu+2n}}{n! \Gamma{(\nu+n+1)}}, \ z\in \mathbb{C}.
  \end{align}

 These identities have three different interpretations. The double series in Entry \ref{entry1} and Entry \ref{entry2} can be interpreted as an iterated series in two possible ways. 
 It can also be interpreted symmetrically where the products of the indices tend to infinity. It is important to note that the identity \eqref{Entry1} in Entry \ref{entry1} with the order of the double sum interchanged was established by B. C. Berndt and A. Zaharescu in \cite{MR2221114}. A few years later, B. C. Berndt, S. Kim and A. Zaharescu in \cite{MR3019715} provided the proof of \eqref{Entry1} as recorded by Ramanujan.  On the other hand, B. C. Berndt, S. Kim and A. Zaherescu in \cite{MR2871168} provided  proof of \eqref{Entry2} of Entry \ref{entry2} with the order of the summation reversed and with additional conditions. Moreover, they gave the proof of \eqref{Entry1} of Entry \ref{entry1} and \eqref{Entry2} of Entry \ref{entry2} under the symmetric interpretation in \cite{MR2871168}. Recently, in 2019, B. C. Berndt, J. Li and A. Zaharescu in \cite{MR4017155} offered the proof of the identity \eqref{Entry2} in Entry \ref{entry2} as given by Ramanujan.   
As an application of Entry \ref{entry1} in \cite{MR2221114} with the order of the double sum reversed, B. C. Berndt and Zaharescu  derived the following beautiful identity associated with $r_2(n)$
where $r_k(n)$ denotes the number of representations of $n$ as a sum of $k$ squares for any $k\geq 2$ 

\begin{align}\label{newrep}
\sideset{}{'}\sum_{0\leq n\leq x} r_2(n)  =\pi x +2\sqrt{x}\sum_{n=0}^\infty\sum_{m=1}^\infty\left\{  \frac{J_1\left(4\pi \sqrt{m(n+\frac{1}{4})x } \ \right)}{\sqrt{m (n+\frac{1}{4})\ } }- \frac{J_1\left(4\pi \sqrt{m(n+\frac{3}{4})x } \ \right)}{\sqrt{m (n+\frac{3}{4})\ } }\right\}.
\end{align}
  The prime $ \prime $ on the summation sign on the left-hand side of \eqref{newrep} implies that weight $1/2$ is considered if $x$ is an integer. It is interesting to note that \eqref{newrep} implies 
  Hardy's formula \cite{hardy1915expression} 
\begin{align}\label{hardy}
\sideset{}{'} \sum_{0\leq n\leq x}  r_2(n)  =\pi x +P(x),
\end{align}
where $P(x)$ is the error term given by the following infinite series
\begin{align}\label{Px}
 P(x)=   \sum_{n=1}^\infty r_2(n)\left(\frac{x}{n}\right)^{1/2}J_1(2\pi \sqrt{nx}).
\end{align}
Dixon and Ferrar \cite{dixon1934some} in 1934 derived another such identity. They proved that, for $\Re(\nu)>0$ and $ x>0$,
 \begin{align}\label{ferrar}
    \sum_{n=0}^\infty r_2(n) n^{\nu/2} K_{\nu}(2\pi \sqrt{nx})=\frac{\Gamma(\nu+1)}{2\pi^{\nu+1}}x^{\frac{\nu}{2}}\sum_{n=0}^\infty  \frac{r_2(n)}{(n+x)^{\nu+1}}.
 \end{align}
 Plugging $\nu=1/2$ in \eqref{ferrar} and using the property of the Bessel function  \cite[p. ~ 80, eq. 13]{MR1349110}  
 \begin{align}\label{property}
K_{\frac{1}{2}}(z)=\sqrt{\frac{\pi }{2z}}e^{-z},
  \end{align}
  we obtain another famous Hardy's result \cite[eq. (2.12)] {hardy1915expression}  
  \begin{align}\label{hardysecond}
 \sum_{n=1}^\infty r_2(n)e^{-s\sqrt{n}}=\frac{2\pi }{s^2}-1+2\pi s\sum_{n=1}^\infty\frac{r_2(n)}{(s^2+4\pi^2n)^\frac{3}{2}},
  \end{align} 
where $\Re (s)>0$. In his paper, Hardy \cite{hardy1915expression} used \eqref{hardysecond} to deduce a lower bound for the error term $P(x)$ given in \eqref{Px}. However, K. Chandrasekharan and R. Narasimhan \cite[eq. (56)]{MR171761}, obtained an interesting identity analogous to \eqref{hardysecond} for Ramanujan tau function. For $\Re( s)>0$,
 \begin{align}
     \sum_{n=1}^\infty \tau(n)e^{-s\sqrt{n}}=2^{36}\pi^{23/2}\Gamma\left(\frac{25}{2}\right)\sum_{n=1}^\infty \frac{\tau(n)}{(s^2+16\pi^2n)^{25/2}}.
 \end{align}

Similar to Entry \ref{entry1}, Entry \ref{entry2} is associated with the Dirichlet divisor problem. The divisor problem is the estimation of the error term $\Delta(x)$ that arises in the summatory function of $d(n)=\sum_{d|n}1$. 
Vorono\"i \cite{voronoi1904fonction} in 1904, expressed  $\Delta(x)$ in terms of  Bessel functions
\begin{align}\label{vor bessel}
 \sideset{}{'}\sum_{n\leq x} d(n)  = x \log x+(2\gamma-1)x + \Delta(x),
\end{align}
where $\gamma$ is the Euler's constant  and
\begin{align}\label{deltax}
   \Delta(x)=  \frac{1}{4}+\sum_{n=1}^\infty d(n ) \left(\frac{x}{n}\right)^{1/2}M_1(4\pi \sqrt{nx}),
\end{align}
 and $M_1(z)$ is defined in \eqref{bessel}. Vorono\"i employed \eqref{deltax} to prove the result $\Delta(x)=O(x^{1/3+\epsilon})$   for each fixed $\varepsilon>0.$ Since then, many number theorists dealt with the above expression to improve the bound of $\Delta(x)$. The occurrence of $M_1(z)$ in \eqref{deltax} indicates that there exists some correlation between Entry \ref{entry2} and \eqref{vor bessel}. Considering this fact, B. C. Berndt and A. Zaharescu derived an identity equivalent to Entry \ref{entry2} in \cite{MR2221114} by introducing twisted divisor sum $d_\chi (n)$ defined by
\begin{align*}
 d_\chi (n) = \sum_{d|n}\chi (d),  
\end{align*}
where $\chi$ is a primitive Dirichlet character modulo $q$ and $\tau(\chi)=\sum_{h=1}^q\chi(q)e^{2\pi i h/q}$. Their identity reads as the following
\begin{align}\label{chid}
  \sideset{}{'}\sum_{n\leq x}  d_\chi (n) =- \frac{x}{\tau({\Bar{\chi}})}\sum_{h=1}^{q-1}\Bar{\chi}(h)\log (2\sin(\pi h/q))+\frac{\sqrt{q}}{\tau({\Bar{\chi}})}\sum_{n=1}^\infty d_{\bar{\chi}} (n)\left( \frac{x}{n}\right)^{1/2}M_1(4\pi \sqrt{nx/q}),
\end{align}
where $\chi$ is a non-principal, even primitive character modulo $q$ and $\tau({ {\chi}})$ is the Gauss sum  of a Dirichlet character modulo $q$ 
is defined by
\begin{align}\label{Gauss}
\tau(\chi)=\sum_{h=1}^q \chi(h)e^{2\pi i h/q} .
\end{align} Hence, \eqref{chid} can be considered a character analogue of Entry \ref{entry2}. The authors in their follow-up papers \cite{MR2994091,MR3351542} generalized Ramanujan's identities by studying Riesz sums for twisted divisor sums. Subsequently, S. Kim in 2017, \cite{MR3541702} extended the definition of twisted sum into twisted sums of divisor functions
\begin{align}\label{twisted sum divisor}
   \sigma_{k,\chi}(n):=\sum_{d|n}d^k  \chi(d), \ \ \ \ \bar{\sigma}_{k,\chi}(n):=\sum_{d|n}d^k  \chi(n/d), \ \ \ \ \sigma_{\chi_1,\chi_2}(n):=\sum_{d|n}d^k  \chi_1(d) \chi_2(n/d),
\end{align}
where $\chi, \chi_1$ and $\chi_2$ are the Dirichlet characters and derived a trigonometric analogue of Riesz sum identities for the twisted sums of divisor functions. As a corollary of the main results, the author obtained a Riesz sum identity for $r_6(n)$ where $r_6(n)$ signifies the number of representations of $n$ as a sum of six squares denoted by $r_6(n)$. The function $r_6(n)$ can be expressed  as (see \cite[p.~ 63]{MR2246314} )
\begin{align}\label{berndt}
    r_6(n)=16\sum_{\substack{d|n \\ \frac{n}{d}\ odd   } }(-1)^{  (n/d-1)/2 }d^2-4\sum_{\substack{d|n \\ d \ odd   } }(-1)^{ (d-1)/2 }d^2.
\end{align}
 Popov \cite[equation (6)]{popov1935uber} extended \eqref{ferrar} for $r_k(n)$. The formula he deduced is as follows  
\begin{align*}
     \sum_{n=0}^\infty r_k(n) n^{\nu/2} K_{\nu}(2\pi \sqrt{n\beta})=\frac{\Gamma(\nu+\frac{k}{2})}{2\pi^{\nu+\frac{k}{2}}}\beta^{\frac{\nu}{2}}\sum_{n=0}^\infty  \frac{r_k(n)}{(n+\beta)^{\nu+\frac{k}{2}}},
\end{align*}
 where $\Re{(\sqrt{\beta})}$ and $\Re{(\nu)}>0$. Identities associated with the $K$-Bessel function have been studied by many number theorists due to its importance. 
  The famous Ramanujan-Guinand \cite[p. 253]{MR947735},  \cite{MR70661} identity involves   the general divisor function 
\begin{align}\label{sigmas}\sigma_s(n)=\sum_{d|n}d^s, \end{align} 
and the $K$-Bessel function 
\begin{align}\label{guinand}
  4x^{\frac{1}{2}}\sum_{n=1}^{\infty}  \frac{\sigma_{\nu}(n)}{n^{\nu/2}}   K_{\nu/2}(2\pi n x)\ +\   \Lambda(s)(x^{(1-\nu)/2} - x^{(\nu-1)/2} )& = 4x^{-\frac{1}{2}}\sum_{n=1}^{\infty}  \frac{\sigma_{\nu}(n)}{n^{\nu/2}}  K_{\nu/2}\left( \frac{2 \pi n}{x}\right)\notag\\
 &+\   \Lambda(-s)( x^{-(1+\nu)/2} - x^{(1+\nu)/2}), 
 \end{align}
where $\Lambda(s)= \pi^{-\frac{s}{2}}\Gamma\left(  \frac{s}{2}\right)\zeta(s)$. Cohen in 2010,   reproved \eqref{guinand} in \cite{MR2744771}. He obtained several beautiful identities as an application of \eqref{guinand}. One of them is mentioned below.
  \begin{proposition}\label{Cohen-type} \cite[p.~62, Theorem 3.4]{MR2744771}  For $\nu \notin \mathbb{Z}$ such that $\Re(\nu) \geq 0$ and any integer N such that $N \geq \lfloor \frac{\Re(\nu)+1}{2}\rfloor$  then 
  \small         
        \begin{align} \label{Cohen Identity}
     &8 \pi x^{\frac{\nu}{2}}\sum_{n=1}^{\infty} \sigma_{-\nu}(n) n^{\nu/2} K_{\nu}(4\pi \sqrt{nx})
    =-\frac{ \Gamma(\nu) \zeta(\nu)}{(2\pi)^{\nu-1} }  +  \frac{\Gamma(1+\nu) \zeta(1+\nu)}{\pi^{\nu+1} 2^\nu x} + \left\{ \frac{\zeta(\nu)x^{\nu-1}}{\sin\left(\frac{\pi \nu}{2}\right)}  + \frac{2}{  \sin \left(\frac{\pi \nu}{2}\right)}\sum_{j=1}^{N} \zeta(2j)\ \zeta(2j-\nu)x^{2j-1} 
    \right.\notag\\&\left.\ \  
     \hspace{5cm}-\pi\frac{\zeta(\nu+1) x^{\nu}} {\cos(\frac{\pi\nu}{2}) }  +\frac{2   }{  \sin \left(\frac{\pi \nu}{2}\right)}x^{2N+1}\sum_{n=1}^{\infty}{\sigma}_{-\nu}(n)  \left(\frac{n^{\nu-2N}-x^{\nu-2N}}  {    n^2-x^2 }\right)
    \right\}.  
    	 \end{align} 
      
      \end{proposition}

 Cohen's identity \eqref{Cohen Identity} proved to be very useful in deriving a Vorono\"i summation formula for the divisor function $\sigma_s(n)$ in \eqref{sigmas}. B. C. Berndt et al. \cite{MR3558223}, as a matter of fact utilised \eqref{Cohen Identity} to obtain the following Vorono\"i summation formula. 
  

\begin{proposition}\label{vorlemma1} \cite[p.~841, Theorem 6.1]{MR3558223} 
Let $0< \alpha< \beta$ and $\alpha, \beta \notin \mathbb{Z}$. Let $f$ denote a function analytic inside a closed contour strictly containing $[\alpha, \beta ]$. Assume $-\frac{1}{2}< \Re{(\nu)}<\frac{1}{2}$,    then 
           \begin{align*}     
          &  \sum_{\alpha<j <\beta}  {\sigma}_{-\nu}(j)f(j)   =   \int_{\alpha} ^{  \beta   }f(t) \left\{   \zeta(1-\nu) \ t^{-\nu}   +    \zeta(\nu+1) \         \right\}  dt  \notag\\
&+ 2\pi \sum_{n=1}^{\infty} \sigma_{-\nu}(n) n^{\nu/2}    \int_{\alpha} ^{  \beta   }f(t)  (t)^{-\frac{\nu}{2}}   \left\{  \left(   \frac{2}{\pi} K_{\nu}(4\pi   \sqrt{nt}) - Y_{\nu}(4\pi   \sqrt{nt})\right) \cos \left(\frac{\pi \nu}{2}\right) 
  -  J_{\nu}(4\pi   \sqrt{nt})  \sin \left(\frac{\pi \nu}{2}\right)  \right\}  dt.
           \end{align*}
\end{proposition}
In a recent study, B. C. Berndt et al. \cite{MR4489312} investigated the general version of \eqref{ferrar} for a certain class of arithmetical function studied by K. Chandrasekharan and R. Narasimhan \cite{MR171761}. More precisely, they considered those arithmetical functions whose functional equation consists of only a single gamma factor. 
   However, the first author and Maji \cite{MR4570432}  investigated identities analogous to \eqref{ferrar} for general divisor function defined by
   \begin{align*}
   \sigma_{z}^{(k)}(n)=\sum_{d^k|n} d^z.    
   \end{align*}
Interestingly, the Dirichlet series associated with $\sigma_{z}^{(k)}(n)$ fails to fit within the setting proposed by K. Chandrasekharan and R. Narasimhan in \cite{MR171761}. Significantly, Maji and the first author \cite{MR4570432} deduced multiple identities, including the majority of Cohen's findings in \cite{MR2744771}. They provided the following identity as a direct consequence of their primary results:
 \begin{proposition}\label{first paper}
    Let $a$ and $x$ be two positive real numbers and $k\geq 1$ be an odd integer. For $\Re(\nu)>0$, we have 
   \begin{align*}
				\left(	\frac{a^2x}{4}\right)^{\frac{\nu}{2}+k+1} \sum_{n=1}^\infty\sigma_k(n)n^{\frac{\nu}{2}}K_{\nu}(a\sqrt{nx})=\frac{(-1)^{\frac{k+1}{2}}}{2}  \Gamma(\nu+k+1)\left(  2\pi \right)^{k+1}\sum_{n=1}^\infty \frac{ \sigma_{k }(n)}{\left(\frac{16\pi^2}{a^2}\frac{n}{x}+1\right)^{\nu+k+1} }+Q_{\nu}(x),
				\end{align*}
    where \begin{align*}
					& Q_{\nu}(x)=-\frac{a^{2k+2}\Gamma(\nu)\zeta(-k)}{2^{2k+4} }x^{k+1}+\frac{a^{2k}\Gamma(1+\nu)\zeta(1-k)}{2^{2k+1}}x^{k}+\frac{1}{2}\Gamma(1+k+\nu)\Gamma(1+k)\zeta(1+k) .
				\end{align*}
    \end{proposition}
    The above identity was also obtained by  B. C. Berndt et al. in \cite[equation (6.11)]{MR4489312} as a particular case of their main result. Proposition \ref{Cohen-type}, \ref{vorlemma1}, and \ref{first paper} can be regarded as identities corresponding to character modulo $1$. We have extended these results to character modulo $q$ in our paper \cite{devika2023}. In  \cite{devika2023}, we have established the Cohen-type identities given in \eqref{Cohen Identity} for these twisted sums of the divisor functions defined in \eqref{twisted sum divisor}. In that same paper, we have derived the Vorono\"i summation formulas by appealing to their Cohen-type identities. 
   
   In this paper, we focus on deriving trigonometric analogues of the identities derived in \cite{devika2023}. More precisely, we offer the identities associated with the $K$-Bessel function and the following weighted sums of divisor functions 
\begin{align}\label{weighted sums}
\sum_{d|n}d^{z} \sin \left( 2\pi d \theta \right), \ \ \sum_{d|n}d^z \sin \left(\frac{2\pi n \theta}{d}\right), \ \ \sum_{d|n}d^z \cos \left( 2\pi d \theta \right), \  \
\sum_{d|n}d^z \cos \left(\frac{2\pi n \theta}{d}\right), 
\end{align}
etc. We prove these identities are equivalent to our previous results in \cite{devika2023}.  Moreover, we present formulas for the following two infinite series,
    \begin{align}\label{r66}
        \sum_{n=1}^\infty r_6(n)n^{\nu/2} K_{\nu}(a \sqrt{nx}),\ \ \ \ \ \ \  \sum_{n=1}^{\infty} r_6(n) e^{-4\pi \sqrt{nx}}.   
    \end{align}
 We also derive an identity from our two main results that give rise to \eqref{r66}.  
   
  The paper is organized as follows:		 
  In the next section, we state the main results.   Section \ref{cohen identities...} and Section \ref{voronoi identities...} provide some special kind of Cohen-type identities and Vorono\"i summation formulas, respectively. In Section \ref{preliminary}, we state some significant results which we need to prove our results. Sections \ref{proof of integer nu} and \ref{proof of cohen identities...} and \ref{proof of voronoi...} are devoted to the proof of identities mentioned in Sections \ref{integer results}, \ref{cohen identities...} and \ref{voronoi identities...} respectively. 

 \vspace{.5cm}

		\section{Main Results  }\label{integer results}	
   
  As mentioned earlier, our results are trigonometric analogue or equivalent versions of the identities associated with twisted sums of divisor functions, which we derived in \cite{devika2023}. Therefore, it is essential to state those identities immediately after their corresponding equivalent version involving finite trigonometric sum and $K$-Bessel function. 
\par Throughout this section, we will set $z=k$ in \eqref{weighted sums} where $k \in \mathbb{Z}_+$.  Before stating our main results, let us recall some important facts.
\par The Hurwitz zeta function is defined by
\begin{align} \label{Hurwitz1}
\zeta(s,\alpha)= \sum_{n=0}^\infty\frac{ 1}{(n+\alpha)^s}, \ \Re(s)>1 \ \mbox{and} \ 0<\alpha<1.
\end{align}
  It is well known that the Dirichlet $L$-function $\sum_{n=1}^{\infty}\frac{\chi(n)}{n^s}$ for $\Re(s)>1$  can be expressed in terms of the Hurwitz zeta function   at rational points \cite[p. 71, equation (16)]{{MR1790423}}
   \begin{align}\label{Hurwitz} 
L(s, \chi) =\frac{1}{q^s}\sum_{r=1}^{q-1} \zeta\left(s,\frac{r}{q}\right)\chi(r),
  \end{align} 
  where $\chi$ is the Dirichlet character modulo $q$ with  $q > 2$. Conversely,  we have  
\begin{align}\label{Lbound}
  \zeta\left(s,\frac{h}{q}\right) =\frac{q^s}{\phi(q)} \sum_{\chi} \bar{\chi}(h) L(s, \chi),
\end{align}
 for $(h,q) = 1$ and $0 < h < q$.
\par Since the identities offered in this paper involve the modified $K$-Bessel function, it is important to mention some related results. The asymptotic estimate for the $K$-Bessel function defined in \eqref{Kbessel} is \cite[p.~202]{MR1349110}
\begin{align*}
 K_{\nu}(x) =\left(\frac{\pi}{2x}  \right)^{\frac12} e^{-x} +O\left(\frac{e^{-x}}{x^{\frac32}}  \right) \ \mbox{as}\ x \rightarrow \infty.
\end{align*}
 Throughout this paper, we will consider $\Re (\nu)\geq 0$ as $K_{-\nu}(x)=K_\nu(x)$ \cite[p.~79]{MR1349110}. We will assume $a$ and $x$ to be two positive real numbers and $\theta,\psi $ are positive real numbers such that $ 0<\theta,\psi<1$. 
   
   \begin{theorem}\label{odd1_based}
    
				Let $k\geq 0 $ be an   even integer. Then for any $\Re{(\nu)}>0$ we have
				\begin{align}\label{p7}
					&\left(	\frac{a^2x}{4}\right)^{\frac{\nu}{2}+k+1}  \sum_{n=1}^{\infty} n^{\nu/2} K_{\nu}(a\sqrt{nx}) \sum_{d|n}d^k \sin \left( 2\pi d \theta \right)\nonumber \\
     &= -\frac{(-1)^{\frac{k}{2}}a^{2k+2}k!}{2^{2k+4}(2\pi)^{k+1}}  \Gamma(\nu)     \left(\zeta(1+k,\theta) - \zeta(1+k, 1-\theta)\right) x^{k+1} 
        +\delta_k \frac{\pi \Gamma(1+\nu)}{4}  \left(\zeta(0, \theta) - \zeta(0, 1-\theta)\right) \notag\\
					&+\frac{(-1)^{\frac{k}{2}}}{4}  (2 \pi )^{k+1} \sum_{d=1}^{\infty}d^k \sum_{ m=0}^{\infty}\left\{  \frac{ \Gamma(\nu+k+1)  }{(1+\frac{16\pi^2d}{a^2x}(m+\theta))^{1+\nu+k} } - \frac{ \Gamma(\nu+k+1)  }{(1+\frac{16\pi^2d}{a^2x}(m+1-\theta))^{1+\nu+k} } \right\}  ,
				\end{align} where $\delta_k$ is given by
    \begin{align}\label{del_k0}
\delta_k=\begin{cases}
1\quad \ \ \text{if } k=0 ,\\
0\quad \ \ \text{else } .
\end{cases}
\end{align}
\end{theorem}
 
We remark that Theorem \ref{odd1_based} is equivalent to \cite[Theorem 2.1]{devika2023}. But for the sake of completeness, we would like to state it here. 
\begin{theorem}\label{M1}
\cite[Theorem 2.1]{devika2023}
				Let  $k\geq 0$ be an even integer, and $\chi$ be an odd primitive Dirichlet character modulo $q$. Then for any $\Re{(\nu)}>0$ we have\begin{align}\label{M11} 
    \left(	\frac{a^2x}{4}\right)^{\frac{\nu}{2}+k+1}	\sum_{n=1}^\infty\sigma_{k,\chi}(n)n^{\frac{\nu}{2}}K_\nu(a\sqrt{nx})=&\frac{(-1)^{\frac{k}{2}}i k!q^ka^{2k+2}}{2^{2k+3}(2\pi)^{k+1}}   \Gamma(\nu)\tau(\chi) L(1+k,\Bar{\chi})\ x^{k+1}+\delta_{k}\frac{\Gamma(1+\nu)L(1,\chi)}{2} \notag\\
					&-\frac{(-1)^{\frac{k}{2}}i}{2q} \tau(\chi) (2 \pi )^{k+1}\sum_{n=1}^\infty \Bar{\sigma}_{k,\Bar{\chi}}(n)  \frac{\Gamma(\nu+k+1)}{\left(\frac{16\pi^2}{a^2q}\frac{n}{x}+1\right)^{\nu+k+1} } ,
      \end{align} where $\delta_k$ is defined in \eqref{del_k0}.
\end{theorem}

\begin{theorem}\label{odd2_based}
    Let $k\geq 2 $ be an   even integer. Then for any $\Re{(\nu)}>0$ we have
				\begin{align}\label{r7}
					& \left(	\frac{a^2x}{4}\right)^{\frac{\nu}{2}+k+1} \sum_{n=1}^{\infty} n^{\nu/2} K_{\nu}(a\sqrt{nx}) \sum_{d|n}d^k \sin \left(\frac{2\pi n \theta}{d}\right)=\frac{(-1)^\frac{k}{2}2^k\pi^{k+1}}{4}\Gamma(\nu+k+1)     \left(\zeta(-k,\theta) - \zeta(-k, 1-\theta)\right)\notag \\&
					+\frac{(-1)^{\frac{k}{2}}}{4}   (2 \pi )^{k+1}  \Gamma(\nu+k+1)\sum_{r=1}^{\infty}  \sum_{ m=0}^{\infty}\left\{  \frac{  (m+\theta)^k  }{(1+\frac{16\pi^2r}{a^2x}(m+\theta))^{1+\nu+k} } - \frac{ (m+1-\theta)^k    }{(1+\frac{16\pi^2r}{a^2x}(m+1-\theta))^{1+\nu+k} } \right\} .
				\end{align}\end{theorem}
  Analogous to Theorem \ref{odd1_based}, one can demonstrate that \cite[Theorem 2.3]{devika2023} is an equivalent version of Theorem \ref{odd2_based}. We will include the result here.

\begin{theorem}\label{M2}
\cite[Theorem 2.3]{devika2023}
				Let  $k\geq 2$ be an even integer, and $\chi$ be an odd primitive Dirichlet character modulo $q$. Then  for any $\Re{(\nu)}>0$  we have\begin{align}\label{ll1} 
					\left(	\frac{a^2x}{4}\right)^{\frac{\nu}{2}+k+1} \sum_{n=1}^\infty\Bar{\sigma}_{k,\chi}(n)n^{\frac{\nu}{2}}K_\nu(a\sqrt{nx})=&\frac{k!}{2}\Gamma(\nu+k+1)L(1+k,\chi) \notag\\
     &- \frac{(-1)^{\frac{k}{2}}i}{2}\tau(\chi) \left(\frac{ 2\pi }{q}\right)^{k+1}\sum_{n=1}^\infty\sigma_{k,\Bar{\chi}}(n)\frac{\Gamma(\nu+k+1)}{\left(\frac{16\pi^2}{a^2q}\frac{n}{x}+1\right)^{\nu+k+1} }. 
				\end{align} 
			\end{theorem}

  Theorem \ref{odd1_based}, together with Theorem \ref{odd2_based}, gives rise to the following beautiful identity associated with $r_6(n)$.
\begin{corollary}\label{cor1r}
     For any $\Re{(\nu)}>0$ we have
    \begin{align}
        &\left(	\frac{a^2x}{4}\right)^{\frac{\nu}{2}+3} \sum_{n=1}^{\infty}  n^{\nu/2} K_{\nu}(a\sqrt{nx}) \sum_{d|n}d^2  \left\{ 16\sin \left(\frac{2\pi n \theta}{d}\right) -4\sin \left( 2\pi d \theta \right)\right\} \notag\\
        &=\frac{16}{3} \pi^3\Gamma(\nu+3)  (\theta-3\theta^2+2\theta^3   ) -\frac{a^6}{256}  \Gamma(\nu)(\cot(\pi \theta)+\cot^3(\pi \theta))x^3\notag\\
        &+ (2 \pi )^3\Gamma(\nu+3)\sum_{n=1}^\infty\sum_{m=0}^\infty   \left\{   \frac{  n^2-4(m+\theta)^2  }{(1+\frac{16\pi^2n}{a^2x}(m+\theta))^{\nu+3} }    -  \frac{  n^2-4(m+1-\theta)^2  }{(1+\frac{16\pi^2n}{a^2x}(m+1-\theta))^{\nu+3} }   \right\}.
    \end{align}
\end{corollary} 
One can also obtain an interesting identity analogous to Hardy's result in \eqref{ferrar}
by substituting $\theta=1/4$.  
\begin{corollary}\label{cor2r}
  For any $\Re{(\nu)}>0$ we have
    \begin{align}
        &\left(	\frac{a^2x}{4}\right)^{\frac{\nu}{2}+3} \sum_{n=1}^{\infty} r_6(n) n^{\nu/2} K_{\nu}(a\sqrt{nx})   = \frac{\pi^3}{2}  \Gamma(\nu+3)   -\frac{a^6}{{128} }  \Gamma(\nu) x^{3}\notag\\
        &+ (2 \pi )^3\Gamma(\nu+3)\sum_{n=1}^\infty\sum_{m=0}^\infty   \left\{   \frac{  n^2-4(m+1/4)^2  }{(1+\frac{16\pi^2n}{a^2x}(m+1/4))^{\nu+3} }    -  \frac{  n^2-4(m+3/4)^2  }{(1+\frac{16\pi^2n}{a^2x}(m+3/4))^{\nu+3} }   \right\}.
    \end{align}
     
\end{corollary}
In particular, $\nu=1/2,a=4\pi $  yields an identity analogous to \eqref{hardysecond}.
 \begin{corollary}\label{cor3r}  We have
 \begin{align}
      &    \sum_{n=1}^{\infty} r_6(n) e^{-4\pi \sqrt{nx}}    = \frac{15}{512 \pi^3} x^{-3} -1  + \frac{15}{32\pi^3} x^{-3}\sum_{n=1}^\infty\sum_{m=0}^\infty   \left\{   \frac{  n^2-4(m+1/4)^2  }{(1+\frac{ n}{ x}(m+1/4))^{ \frac{7}{2}} }    -  \frac{  n^2-4(m+3/4)^2  }{(1+\frac{ n}{ x}(m+3/4))^{\frac{7}{2}} }   \right\}.
    \end{align}
\end{corollary}
 \vspace{.5cm}
			 The corresponding cosine versions of the above identities are listed below.

     \begin{theorem}\label{even1_based}
      Let $k\geq 1 $ be an odd integer. Then for any $\Re{(\nu)}>0$ we have
     \begin{align}\label{1290}
					&	\left(	\frac{a^2x}{4}\right)^{\frac{\nu}{2}+k+1} \sum_{n=1}^{\infty} n^{\nu/2} K_{\nu}(a\sqrt{nx}) \sum_{d|n}d^k \cos \left( 2\pi d \theta \right) = \frac{(-1)^{\frac{k-1}{2}}a^{2k+2}k! }{2^{2k+4}(2\pi)^{k+1} }  \Gamma(\nu)\left\{  
					 \zeta(k+1,\theta)+\zeta(k+1,1-\theta)   \right\}x^{k+1}  \notag\\
					&  
     +\frac{(-1)^{\frac{k+1}{2}}(2 \pi )^{k+1} \Gamma(\nu+k+1)}{4}  \sum_{d=1}^\infty d^k \sum_{m=0 }^\infty \left\{   \frac{1}{\left(\frac{16\pi^2}{a^2}\frac{d(m+\theta)}{x}+1\right)^{\nu+k+1} }  
					 +\frac{1}{\left(\frac{16\pi^2}{a^2}\frac{d(m+1-\theta)}{x}+1\right)^{\nu+k+1} }      \right\}\notag\\
					& \ \  
      -\delta_{k,1}\ \frac{ a^2 }{16}\Gamma(1+\nu)x     ,
				\end{align}
 where $\delta_{k,1}$ is given by
 \begin{align}\label{del_k}
\delta_{k,1}=\begin{cases}
1\quad \ \ \text{if } k=1 ,\\
0\quad \ \ \text{else } .
\end{cases}
\end{align}
\end{theorem} 
 
 Theorem \ref{even1_based} is equivalent to the following theorem, mentioned below.  
\begin{theorem}\label{thmeven1} 
\cite[Theorem 2.5]{devika2023} 
  	Let $k$ be a positive odd integer and let $\chi$ be a non-principal even primitive Dirichlet character modulo q, then for any $\Re{(\nu)}>0$  we have \begin{align}\label{thm3}
					\left(	\frac{a^2x}{4}\right)^{\frac{\nu}{2}+k+1} \sum_{n=1}^\infty\sigma_{k,\chi}(n)n^{\frac{\nu}{2}}K_\nu(a\sqrt{nx})=&\frac{(-1)^{\frac{k-1}{2}}\Gamma(k+1)a^{2k+2}q^k}{2^{2k+3}(2\pi)^{k+1}}  \tau(\chi)\Gamma(\nu)L(1+k,\Bar{\chi}) x^{k+1}\notag\\
					&+\frac{(-1)^{\frac{k+1}{2}}}{2q}\tau(\chi) (2 \pi )^{k+1}\sum_{n=1}^\infty \Bar{\sigma}_{k,\Bar{\chi}}(n) \frac{\Gamma(\nu+k+1)}{\left(\frac{16\pi^2}{a^2q}\frac{n}{x}+1\right)^{\nu+k+1} }.
				\end{align} 
    \end{theorem}

		\vspace{.3cm}	 
 
\begin{theorem}\label{even2_based}
Let $k\geq 1 $ be an odd integer. Then for any $\Re{(\nu)}>0$ we have 
			\begin{align*}
					&\left(	\frac{a^2x}{4}\right)^{\frac{\nu}{2}+k+1} \sum_{n=1}^{\infty} n^{\nu/2} K_{\nu}(a\sqrt{nx}) \sum_{d|n}d^k \cos \left(\frac{ 2\pi n \theta }{d}\right) =    
					  \frac{(-1)^{\frac{k+1}{2}}\left(  2\pi  \right)^{k+1} }{8 }\Gamma(\nu+k+1) \left\{ \zeta(-k,\theta)+\zeta(-k,1-\theta)    \right\} 
       \\
					&+\frac{(-1)^{\frac{k+1}{2}}(2 \pi )^{k+1}}{4 }   \Gamma(\nu+k+1)  \sum_{r=1}^\infty  \sum_{m=0 }^\infty \left\{   \frac{(m+ \theta)^k}{\left(\frac{16\pi^2}{a^2}\frac{r(m+\theta)}{x}+1\right)^{\nu+k+1} }   +\frac{(m+1-\theta)^k}{\left(\frac{16\pi^2}{a^2}\frac{r(m+1-\theta)}{x}+1\right)^{\nu+k+1} }        \right\}
     \notag\\&-\frac{a^{2k+2} }{2^{2k+4} }  \zeta(-k)\Gamma(\nu) x^{k+1}.
				\end{align*}	 
\end{theorem} 
 
Theorem \ref{even2_based} is equivalent to the following result. 
      \begin{theorem}\label{even2}
      \cite[Theorem 2.7]{devika2023}
				Let k be a positive odd integer and $\chi$ be a non-principal even primitive Dirichlet  character modulo q. Then  for any $\Re{(\nu)}>0$ we have
				 \begin{align} \label{thm4} 
					\left(	\frac{a^2x}{4}\right)^{\frac{\nu}{2}+k+1}\sum_{n=1}^\infty\Bar{\sigma}_{k,\chi}(n)n^{\frac{\nu}{2}}K_\nu(a\sqrt{nx})=&\frac{\Gamma(k+1)}{2}\Gamma(\nu+k+1)L(1+k,\chi) \notag\\
     &+\frac{(-1)^{\frac{k+1}{2}}}{2}\tau(\chi) \left(\frac{ 2\pi }{q}\right)^{k+1}\sum_{n=1}^\infty\sigma_{k,\Bar{\chi}}(n)\frac{\Gamma(\nu+k+1)}{\left(\frac{16\pi^2}{a^2q}\frac{n}{x}+1\right)^{\nu+k+1} }.
				\end{align} 
			\end{theorem}

   Next, we state the identities involving two trigonometric functions,  which are the following:
 \begin{theorem}\label{botheven_odd1_based}
Let $k\geq 1 $ be an odd integer. Then for any $\Re{(\nu)}>0$ we have 
		\begin{align}\label{rr3}
					&\left(	\frac{a^2x}{4}\right)^{\frac{\nu}{2}+k+1}\sum_{n=1}^{\infty} n^{\nu/2} K_{\nu}(a\sqrt{nx}) \sum_{d|n}d^k  \sin \left( 2\pi d \theta\right)\sin \left( \frac{2\pi n \psi}{d}\right)\notag\\
     =&-\frac{(-1)^{\frac{k+1}{2}}}{8 \left( 2\pi \right)^{-k-1}} \Gamma(\nu+k+1)           \sum_{m,n\geq 0} \left\{ \frac{ (n+\psi)^k}{\left(\frac{16\pi^2}{a^2}\frac{(n+\psi)(m+\theta)}{x}+ 1 \right)^{\nu+k+1}}  -\frac{ (n+1- \psi)^k}{\left(\frac{16\pi^2}{a^2}\frac{(n+1-\psi)(m+\theta)}{x}+ 1 \right)^{\nu+k+1}}   \right.\notag\\&\left.\ \ 
        - \frac{ (n+\psi)^k}{\left(\frac{16\pi^2}{a^2}\frac{(n+\psi)(m+1-\theta)}{x}+ 1 \right)^{\nu+k+1}}                                +  \frac{ (n+1-\psi)^k}{\left(\frac{16\pi^2}{a^2}\frac{(n+1-\psi)(m+1-\theta)}{x}+ 1 \right)^{\nu+k+1}}         \right\}     .
				\end{align}
 \end{theorem}

 \begin{theorem}\label{botheven_odd2_based}
Let $k\geq 1 $ be an odd integer. Then for any $\Re{(\nu)}>0$ we have 
			\begin{align*}
					&\left(	\frac{a^2x}{4}\right)^{\frac{\nu}{2}+k+1} \sum_{n=1}^{\infty} n^{\nu/2} K_{\nu}(a\sqrt{nx}) \sum_{d|n}d^k  \cos \left( 2\pi d \theta\right)\cos \left( \frac{2\pi n \psi}{d}\right) 
					= \frac{(-1)^{\frac{k-1}{2}}a^{2k+2}k! }{2^{2k+4}(2\pi)^{k+1} }   \Gamma(\nu)\left\{  
					 \zeta(k+1, \theta)
       \right.\notag\\&\left.\ 
       +\zeta(k+1,1- \theta)    \right\} x^{k+1}+\frac{(-1)^{\frac{k+1}{2}}}{8 }  \left( {2\pi } \right)^{k+1}  \Gamma(\nu+k+1)\sum_{n,m\geq 0}\left\{   
				\frac{(n+\psi)^k }{\left(\frac{16\pi^2}{a^2}\frac{(m+\theta)(n+\psi)}{x}+1\right)^{\nu+k+1} }
    \right.\notag\\&\left.\ 
    +\frac{(n+1-\psi)^k }{\left(\frac{16\pi^2}{a^2}\frac{(m+\theta)(n+1-\psi)}{x}+1\right)^{\nu+k+1} }
    +\frac{(n+\psi)^k }{\left(\frac{16\pi^2}{a^2}\frac{(m+1-\theta)(n+\psi)}{x}+1\right)^{\nu+k+1} }
    +\frac{(n+1-\psi)^k }{\left(\frac{16\pi^2}{a^2}\frac{(m+1-\theta)(n+1-\psi)}{x}+1\right)^{\nu+k+1} }
    \right\}.
    \end{align*}
 \end{theorem}
 Our next result  \cite[Theorem 2.10]{devika2023} is the equivalent version of Theorem \ref{botheven_odd1_based}. But to show the equivalence of Theorem \ref{botheven_odd2_based} with our next result \cite[Theorem 2.10]{devika2023}, we need the help of Theorem \ref{even1_based} and Theorem \ref{even2_based}.

			\begin{theorem}\label{botheven_odd}
   \cite[Theorem 2.10]{devika2023}
				Let k be a positive odd integer. Let $\chi_1$ and $\chi_2$ be primitive characters modulo p and q, respectively, such that either both are non-principal even characters or both are odd characters. Then for any $\Re{(\nu)}>0$ we have
				\begin{align}\label{THM5}
					\left(	\frac{a^2x}{4}\right)^{\frac{\nu}{2}+k+1} \sum_{n=1}^\infty \sigma_{k,\chi_1,\chi_2}(n)n^{\frac{\nu}{2}}K_\nu(a\sqrt{nx})=&\frac{(-1)^{\frac{k+1}{2}}}{2p}  \left(\frac{2\pi }{q}\right)^{k+1} \tau(\chi_1)\tau(\chi_2)  \sum_{n=1}^\infty \sigma_{k,\Bar{\chi_2},\Bar{\chi_1}}(n)\frac{\Gamma(\nu+k+1)}{\left(\frac{16\pi^2}{a^2pq}\frac{n}{x}+1\right)^{\nu+k+1} }.
				\end{align}
     \end{theorem}
Our next results concern both sine and cosine functions.
 \begin{theorem}\label{even-odd1_based}
Let $k\geq 2 $ be an   even integer. Then for any $\Re{(\nu)}>0$ we have		\begin{align*}
					&\left(	\frac{a^2x}{4}\right)^{\frac{\nu}{2}+k+1}\sum_{n=1}^{\infty} n^{\nu/2} K_{\nu}(a\sqrt{nx}) \sum_{d|n}d^k  \cos \left( 2\pi d \theta\right)\sin \left( \frac{2\pi n \psi}{d}\right) \\&=\frac{(-1)^{\frac{k}{2}}\left( 2\pi \right)^{k+1}}{8 }    \Gamma(\nu+k+1) \sum_{m,n\geq 0}^\infty  \left\{\frac{(n+\psi)^k}{\left(\frac{16\pi^2}{a^2}\frac{(n+\psi)(m+\theta)}{x}+1\right)^{\nu+k+1} }-\frac{(n+1-\psi)^k}{\left(\frac{16\pi^2}{a^2}\frac{(n+1-\psi)(m+\theta)}{x}+1\right)^{\nu+k+1} }
					\right.\notag\\&\left.\ \  \ \ \ \ \ \ \ \ \ \ \ \ \ \ \ \ \  \ \ \ \ \ \ \ \ \ \ \ \ \ \ \ \ \ \ \ \ \ \ -\frac{(n+1-\psi)^k}{\left(\frac{16\pi^2}{a^2}\frac{(n+1-\psi)(m+1-\theta)}{x}+1\right)^{\nu+k+1} }
					  +\frac{(n+\psi)^k}{\left(\frac{16\pi^2}{a^2}\frac{(n+\psi)(m+1-\theta)}{x}+1\right)^{\nu+k+1} }           \right\} .   \notag \\ 
				\end{align*}
 \end{theorem}
  \begin{theorem}\label{even-odd2_based}
Let $k\geq 0 $ be an   even integer. Then for any $\Re{(\nu)}>0$ we have	\begin{align*}
					&\left(	\frac{a^2x}{4}\right)^{\frac{\nu}{2}+k+1} \sum_{n=1}^{\infty} n^{\nu/2} K_{\nu}(a\sqrt{nx}) \sum_{d|n}d^k  \sin \left( 2\pi d \theta\right)\cos \left( \frac{2\pi n \psi}{d}\right)  
					=-\frac{(-1)^{\frac{k}{2}}a^{2k+2}k!}{ 2^{2k+4}(2\pi)^{k+1}}   \Gamma(\nu)     \left(\zeta(1+k,\theta)
      \right.\notag\\&\left.\ - \zeta(1+k, 1-\theta)\right)x^{k+1}  +\frac{(-1)^{\frac{k}{2}}\left( 2\pi \right)^{k+1}}{8 }     \Gamma(\nu+k+1) \sum_{m,n\geq 0}^\infty  \left\{\frac{(n+\psi)^k}{\left(\frac{16\pi^2}{a^2}\frac{(n+\psi)(m+\theta)}{x}+1\right)^{\nu+k+1} }
					   \right.\notag\\&\left.\  
      +\frac{(n+1-\psi)^k}{\left(\frac{16\pi^2}{a^2}\frac{(n+1-\psi)(m+\theta)}{x}+1\right)^{\nu+k+1} } -\frac{(n+1-\psi)^k}{\left(\frac{16\pi^2}{a^2}\frac{(n+1-\psi)(m+1-\theta)}{x}+1\right)^{\nu+k+1} }
	 -\ \frac{(n+\psi)^k}{\left(\frac{16\pi^2}{a^2}\frac{(n+\psi)(m+1-\theta)}{x}+1\right)^{\nu+k+1} }  \right\}. 
	\end{align*}
 \end{theorem}
  With the aid of Theorem \ref{odd2_based}, one can show the equivalence of  
 Theorem \ref{even-odd1_based} and our next result \cite[Theorem 2.16]{devika2023}. Similarly, with the help of Theorem \ref{odd1_based}, the equivalence of Theorem \ref{even-odd2_based} and the following result \cite[Theorem 2.16]{devika2023} can be proved.

  \begin{theorem}\label{even-odd}
  \cite[Theorem 2.16]{devika2023}
				Let $k\geq 0$ be an even integer. Let $\chi_1$ and $\chi_2$ be primitive characters modulo $p$ and $q $, respectively, such that one of them is a non-principal even character and the other is an odd character. Then for $\Re{(\nu)}>0$ we have 
				\begin{align}\label{THM6}
					\left(	\frac{a^2x}{4}\right)^{\frac{\nu}{2}+k+1}\sum_{n=1}^\infty \sigma_{k,\chi_1,\chi_2}(n)n^{\frac{\nu}{2}}K_\nu(a\sqrt{nx})=\frac{(-1)^{\frac{k}{2}}}{2ip}  \left(\frac{2\pi }{q}\right)^{k+1} \tau(\chi_1)\tau(\chi_2)  \sum_{n=1}^\infty \sigma_{k,\Bar{\chi_2},\Bar{\chi_1}}(n)\frac{\Gamma(\nu+k+1)}{\left(\frac{16\pi^2}{a^2pq}\frac{n}{x}+1\right)^{\nu+k+1} }. 
				\end{align} 
			\end{theorem}
\begin{remark}
 It is useful to mention that we have obtained Theorem \ref{M1}, Theorem \ref{M2}, Theorem \ref{thmeven1}, Theorem \ref{even2}, Theorem \ref{botheven_odd} and Theorem \ref{even-odd} in our previous paper \cite{devika2023}, without using any known identities. The proofs were entirely based on analytic techniques.
\end{remark}


\section{Cohen-Type Identity}\label{cohen identities...}

In this section, we focus on identities which deal with  $z=-\nu$ such that $\Re{(\nu)}\geq0$. We will assume that $x>0$ and $0<\theta, \psi<1$. Let us define
\begin{align}\label{Halpha}
H_{\alpha, \beta}= \{ (n-\alpha)(r-\beta),\ n, r \in \mathbb{N}\},    
\end{align}
where $0\leq \alpha, \beta <1$.
\begin{theorem}\label{oddcohen based}
				Let $\nu \notin \mathbb{Z}$ such that $\Re{(\nu)}\geq0$. Then, for any integer $ N$ such that 
 $  N\geq \lfloor\frac{\Re{(\nu)}+1}{2}\rfloor$, we have
				\begin{align*}
					&8 \pi x^{\nu/2} \sum_{n=1}^{\infty} n^{\nu/2} K_{\nu}(4 \pi\sqrt{nx}) \sum_{d|n}d^{-\nu} \sin \left( 2\pi d \theta \right)= \frac{1 }{\cos\left(\frac{\pi \nu}{2}\right)}     \zeta(\nu+1) \left(\zeta(1,\theta) - \zeta(1, 1-\theta)\right)x^\nu\\
     & -\frac{\pi }{ 2 \sin\left(\frac{\pi \nu}{2}\right)}      \left(\zeta(1-\nu,\theta) - \zeta(1-\nu, 1-\theta)\right)
     +\frac{1 }{ 2 x\cos\left(\frac{\pi \nu}{2}\right)}      \left(\zeta(-\nu,\theta) - \zeta(-\nu, 1-\theta)\right)\\
     & -\frac{1 }{  \cos\left(\frac{\pi \nu}{2}\right)}     \sum_{j=1}^N \zeta(2j) \left(\zeta(2j-\nu,\theta) - \zeta(2j-\nu, 1-\theta)\right)x^{2j-1}\\
      & 
 -\frac{1 }{  \cos\left(\frac{\pi \nu}{2}\right)} x^{2N+1}\sum_{d=1}^{\infty}d^{-\nu-1
 } \sum_{ m=0}^{\infty}\left\{  
 \frac{ \left(  d(m+\theta)  \right)^{\nu+1-2N}-    x    ^{\nu+1-2N}   }{ (m+\theta)\left(d^2(m+\theta)^2-x^2 \right)} 
 -  \frac{ \left(  d(m+1-\theta)  \right)^{\nu+1-2N}-   x   ^{\nu+1-2N}   }{ (m+1-\theta)\left(d^2(m+1-\theta)^2-x^2 \right)} 
 \right\},
				\end{align*}
     provided   $x\notin  H_{\theta,0}$.
			\end{theorem} 
     The equivalent version of Theorem \ref{oddcohen based} is the following \cite[Theorem 3.5]{devika2023}.
\begin{theorem}\label{oddcohen}
\cite[Theorem 3.5]{devika2023}
Let  $\nu \notin \mathbb{Z}$ such that $\Re{(\nu)}\geq 0 .$ Let $\chi$ be an odd primitive character modulo $ q.$  If $N$ is any integer such that  $ N\geq \lfloor\frac{\Re{(\nu)}+1}{2}\rfloor$, then 
   \begin{align}\label{THM7}
  & 8\pi x^{\nu/2} \sum_{n=1}^{\infty} \sigma_{-\nu,  {\chi}}(n) n^{\nu/2} K_{\nu}(4\pi \sqrt{nx})
    =-\frac{ \Gamma(\nu) L(\nu,  {\chi})}{(2\pi)^{\nu-1} }  +  \frac{2\Gamma(1+\nu) L(1+\nu,  {\chi})}{(2\pi)^{\nu+1} }x^{-1} 
     + \frac{iq^{1-\nu}  }{\tau(\chi)} \left\{ \frac{2\zeta(\nu+1)L(1,\bar{\chi})(qx)^{\nu}} {\cos(\frac{\pi\nu}{2}) } 
    \right.\notag\\&\left.\ \ 
     \ \ \ \ \ \ \ 
      - \frac{2}{  \cos \left(\frac{\pi \nu}{2}\right)}\sum_{j=1}^{N} \zeta(2j)\ L(2j-\nu, \bar{\chi})(qx)^{2j-1} 
     -\frac{2  }{ \cos\left(\frac{\pi \nu}{2}\right)   }(qx)^{2N+1}\sum_{n=1}^{\infty}\bar{\sigma}_{-\nu, \bar{\chi}}(n)  \frac{ \left( n^{\nu+1-2N}-(qx)^{\nu+1-2N}\right)}{ n\ (n^2-(qx)^2)}
    \right\}, 
    	 \end{align}
 provided   $qx\notin \mathbb{Z}_+$.
 \end{theorem}
  
\begin{theorem}\label{oddcohen2 based}
 Let $\nu \notin \mathbb{Z}$ such that $\Re{(\nu)}\geq0$. Then, for any integer $ N$ such that 
 $  N\geq \lfloor\frac{\Re{(\nu)}+1}{2}\rfloor$, we have
				\begin{align*}
					&8 \pi x^{\nu/2} \sum_{n=1}^{\infty} n^{\nu/2} K_{\nu}(4 \pi\sqrt{nx}) \sum_{d|n}d^{-\nu} \sin \left( \frac{2\pi n \theta}{d}  \right)= 
\frac{2}{(2\pi)^\nu  } \Gamma(\nu) \zeta(\nu)  \left(\zeta(1, \theta) - \zeta(1, 1-\theta)\right)\notag\\
    + & \frac{\pi }{ 2 \sin\left(\frac{\pi \nu}{2}\right)}   x^\nu   \left(\zeta(1+\nu, \theta) - \zeta(1+\nu, 1-\theta)\right)
     +\frac{x^{\nu-1} }{ 2 \cos\left(\frac{\pi \nu}{2}\right)}    \left(\zeta(\nu,\theta) - \zeta(\nu, 1-\theta)\right)\notag\\
    + & \frac{1 }{  \cos\left(\frac{\pi \nu}{2}\right)}     \sum_{j=1}^{N-1} \zeta(2j+1-\nu) \left(\zeta(2j+1,\theta) - \zeta(2j+1 , 1-\theta)\right)x^{2j} 
       \notag\\
+&  \frac{x^{2N}  }{  \cos\left(\frac{\pi \nu}{2}\right)} \sum_{r=1}^{\infty}     \sum_{  m=0}^\infty
      \left\{  (m+\theta)^{-\nu}  \frac{   (r(m+\theta))^{\nu+1-2N}-x^{\nu+1-2N} }{ r^2(m+\theta)^2-x^2} 
      -  (m+1-\theta)^{-\nu}  \frac{   (r(m+1-\theta))^{\nu+1-2N}-x^{\nu+1-2N} }{   r^2(m+1-\theta)^2-x^2} 
      \right\},    
     \end{align*}
      provided   $x\notin  H_{\theta,0}$.
     \end{theorem}

Analogous to Theorem \ref{oddcohen based}, one can  prove that the equivalent version of Theorem \ref{oddcohen2 based} is \cite[Theorem 3.7]{devika2023}. 
 \par
Next, we state the corresponding cosine version of the above identities.

\begin{theorem} \label{evencohen based}
Let $\nu \notin \mathbb{Z}$ such that $\Re{(\nu)}\geq0$. Then, for any integer $ N$ such that 
 $  N\geq \lfloor\frac{\Re{(\nu)}+1}{2}\rfloor$, we have
				\begin{align*}
					&8 \pi x^{\nu/2} \sum_{n=1}^{\infty} n^{\nu/2} K_{\nu}(4 \pi\sqrt{nx}) \sum_{d|n}d^{-\nu} \cos \left( 2\pi d \theta \right)= -\frac{\pi }{\cos\left(\frac{\pi \nu}{2}\right)}     \zeta(\nu+1)  x^\nu\\
     & -\frac{\pi }{ 2 \cos\left(\frac{\pi \nu}{2}\right)}      \left(\zeta(1-\nu,\theta) + \zeta(1-\nu, 1-\theta)\right)
     -\frac{1 }{ 2 x\sin\left(\frac{\pi \nu}{2}\right)}      \left(\zeta(-\nu,\theta) + \zeta(-\nu, 1-\theta)\right)\\
     & +\frac{1 }{  \sin\left(\frac{\pi \nu}{2}\right)}     \sum_{j=1}^N \zeta(2j) \left(\zeta(2j-\nu,\theta) + \zeta(2j-\nu, 1-\theta)\right)x^{2j-1}\\
      & 
 +\frac{1 }{  \sin\left(\frac{\pi \nu}{2}\right)} x^{2N+1}\sum_{d=1}^{\infty}d^{-\nu
 } \sum_{ m=0}^{\infty}\left\{  
 \frac{ \left(  d(m+\theta)  \right)^{\nu-2N}-    x^{\nu-2N}   }{  \left(d^2(m+\theta)^2-x^2 \right)} 
 -  \frac{ \left(  d(m+1-\theta)  \right)^{\nu-2N}-   x^{\nu-2N}   }{  \left(d^2(m+1-\theta)^2-x^2 \right)} 
 \right\},
				\end{align*}
     provided   $x\notin  H_{\theta,0}$.
			\end{theorem}
For deriving Theorem \ref{evencohen based}, one needs to use \cite[Theorem 3.1]{devika2023} and Proposition \ref{Cohen-type}. Conversely, Theorem \ref{evencohen based}  will imply \cite[Theorem 3.1]{devika2023}.

\begin{theorem} \label{evencohen2 based}
 Let $\nu \notin \mathbb{Z}$ such that $\Re{(\nu)}\geq0$. Then, for any integer $ N$ such that 
 $  N\geq \lfloor\frac{\Re{(\nu)}+1}{2}\rfloor$, we have
				\begin{align*}
					&8 \pi x^{\nu/2} \sum_{n=1}^{\infty} n^{\nu/2} K_{\nu}(4 \pi\sqrt{nx}) \sum_{d|n}d^{-\nu} \cos \left( \frac{2\pi n \theta}{d}  \right) 
=
     -\frac{ \Gamma(\nu) \zeta(\nu )}{(2\pi)^{\nu-1} }   +\frac{x^{\nu-1}}{2\sin \left(\frac{\pi \nu}{2} \right)}   \left\{   \zeta(\nu, \theta )+ \zeta(\nu, 1- \theta)  \right\}    \notag\\
              & -\frac{\pi \ x^\nu}{2\cos \left(\frac{\pi \nu}{2} \right)} \left\{    \zeta(1+\nu, \theta )+ \zeta(1+\nu, 1- \theta )        \right\}  +\frac{1}{\sin\left(\frac{\pi \nu}{2} \right)}\sum_{j=1}^N \zeta(2j-\nu) x^{2j-1}
      \left\{   \zeta(2j, \theta )+ \zeta(2j , 1- \theta )    \right\}\notag\\
      &+\frac{ x^{2N+1} }{  \sin \left(\frac{\pi \nu}{2}\right)} \sum_{r=1}^{\infty}    \sum_{  m=0}^\infty   
      \left\{ (m+\theta)^{-\nu}  \frac{   (r(m+\theta))^{\nu-2N}-x^{\nu-2N} }{ r^2(m+\theta)^2-x^2} 
      + (m+1-\theta)^{-\nu}  \frac{   (r(m+1-\theta))^{\nu-2N}-x^{\nu-2N} }{  r^2(m+1-\theta)^2-x^2} 
      \right\},  
\end{align*}
   provided   $x\notin  H_{\theta,0}$.
\end{theorem}
Analogous to Theorem \ref{evencohen based},  the equivalent version of Theorem \ref{evencohen2 based} is \cite[Theorem 3.3]{devika2023}.

  \par
Now, we state the identities involving two trigonometric functions: 
 \begin{theorem} \label{cohen2 even-odd1based}
 Let $\nu \notin \mathbb{Z}$ such that $\Re{(\nu)}\geq0$. Then, for any integer $ N$ such that 
 $  N\geq \lfloor\frac{\Re{(\nu)}+1}{2}\rfloor$, we have
				\begin{align*}
					8\pi x^{\frac{\nu}{2}}&\sum_{n=1}^{\infty} n^{\nu/2} K_{\nu}(a\sqrt{nx}) \sum_{d|n}d^{-\nu}  \sin \left( 2\pi d \theta\right)\sin \left( \frac{2\pi n \psi}{d}\right)\notag\\
      &=\frac{  1}{2 \sin \left(\frac{\pi \nu}{2}\right)} \left(\zeta(1-\nu, \theta) - \zeta(1-\nu, 1-\theta)\right) \left(\zeta(1, \psi) - \zeta(1, 1-\psi)\right)\notag\\
      &-   \frac{  1}{2 \sin \left(\frac{\pi \nu}{2}\right)}x^{\nu} \left(\zeta(1, \theta) - \zeta(1 , 1-\theta)\right) \left(\zeta(\nu+1, \psi) - \zeta(\nu+1, 1-\psi)\right)\notag\\
      &+\frac{  1}{2 \sin \left(\frac{\pi \nu}{2}\right)} \sum_{j=1}^{N-1} x^{2j}\left(\zeta(2j+1-\nu, \theta) - \zeta(2j+1-\nu, 1-\theta)\right)    \left(\zeta(2j+1, \psi) - \zeta(2j+1, 1-\psi)\right)\notag\\
      &+\frac{  1}{2 \sin \left(\frac{\pi \nu}{2}\right)}x^{2N}\sum_{m,n\geq 0}^{\infty} \left\{ \frac{(n+\psi)^{-\nu-1} }{(m+\theta)}\left( \frac{((n+\psi)(m+\theta))^{\nu-2N+2}-x^{\nu-2N+2}}{ ((n+\psi)(m+\theta))^2-x^2} \right) 
   \right.\notag\\&\left.\ \ \ \ \ \ \ \  \ \ \  \ \ 
     -\ \frac{(n+1-\psi)^{-\nu-1} }{(m+\theta)}\left( \frac{((n+1-\psi)(m+\theta))^{\nu-2N+2}-x^{\nu-2N+2}}{ ((n+1-\psi)(m+\theta))^2-x^2} \right) 
    \right.\notag\\&\left.\ \ \ \ \ \ \ \  \ \ \  \ \ 
    -\ \frac{(n+\psi)^{-\nu-1} }{(m+1-\theta)}\left( \frac{((n+\psi)(m+1-\theta))^{\nu-2N+2}-x^{\nu-2N+2}}{ ((n+\psi)(m+1-\theta))^2-x^2} \right) 
       \right.\notag\\&\left.\ \  \ \ \ \ \ \  \ \ \  \ \  
      +\ \frac{(n+1-\psi)^{-\nu-1} }{(m+1-\theta)}\left( \frac{((n+1-\psi)(m+1-\theta))^{\nu-2N+2}-x^{\nu-2N+2}}{ ((n+1-\psi)(m+1-\theta))^2-x^2 } \right) 
   \right\},  \end{align*}
    provided   $x\notin  H_{\theta,\psi}$.
\end{theorem}
 \begin{theorem} \label{cohen2 even-odd2based}
 Let $\nu \notin \mathbb{Z}$ such that $\Re{(\nu)}\geq0$. Then, for any integer $ N$ such that 
 $  N\geq \lfloor\frac{\Re{(\nu)}+1}{2}\rfloor$, we have
				\begin{align*}
					&8\pi x^{\frac{\nu}{2}}\sum_{n=1}^{\infty} n^{\nu/2} K_{\nu}(a\sqrt{nx}) \sum_{d|n}d^{-\nu}  \cos \left( 2\pi d \theta\right)\cos \left( \frac{2\pi n \psi}{d}\right)
 \notag\\
 &=-\frac{\pi \ x^\nu}{2\cos \left(\frac{\pi \nu}{2} \right)} \left\{    \zeta(1+\nu, \psi)+ \zeta(1+\nu, 1- \psi )   \right\} 
 - \frac{\pi }{ 2 \cos\left(\frac{\pi \nu}{2}\right)}      \left(\zeta(1-\nu,\theta) + \zeta(1-\nu, 1-\theta)\right)   \notag\\   
 &+\frac{1}{2  \sin\left(\frac{\pi \nu}{2}\right) }\sum_{j=1}^N  x^{2j-1} \left\{\zeta(2j,\psi)+\zeta(2j,1-\psi)       \right\}  \left\{\zeta(2j-\nu, \theta) + \zeta(2j-\nu, 1-\theta)  \right\} \notag\\  
  &+\frac{  1}{2 \sin \left(\frac{\pi \nu}{2}\right)}x^{2N+1}\sum_{m,n\geq 0}^{\infty} \left\{  {(n+\psi)^{-\nu} } \left( \frac{((n+\psi)(m+\theta))^{\nu-2N}-x^{\nu-2N}}{ ((n+\psi)(m+\theta))^2-x^2} \right) 
   \right.\notag\\&\left.\ \ \ \ \ \ \ \  \ \ \  \ \ 
     +\  {(n+1-\psi)^{-\nu} } \left( \frac{((n+1-\psi)(m+\theta))^{\nu-2N}-x^{\nu-2N}}{ ((n+1-\psi)(m+\theta))^2-x^2} \right) 
    \right.\notag\\&\left.\ \ \ \ \ \ \ \  \ \ \  \ \ 
    +\  {(n+\psi)^{-\nu} } \left( \frac{((n+\psi)(m+1-\theta))^{\nu-2N}-x^{\nu-2N}}{ ((n+\psi)(m+1-\theta))^2-x^2} \right) 
       \right.\notag\\&\left.\ \  \ \ \ \ \ \  \ \ \  \ \  
      +\  {(n+1-\psi)^{-\nu} } \left( \frac{((n+1-\psi)(m+1-\theta))^{\nu-2N}-x^{\nu-2N}}{ ((n+1-\psi)(m+1-\theta))^2-x^2 } \right) 
   \right\},   
 \end{align*}
  provided   $x\notin  H_{\theta,\psi}$.
 \end{theorem}
The equivalent version of the Theorem \ref{cohen2 even-odd1based} is \cite[Theorem 3.11]{devika2023}. To prove Theorem \ref{cohen2 even-odd2based}, one requires   \cite[Theorem 3.9]{devika2023}, Theorem \ref{evencohen based}, Theorem \ref{evencohen2 based} and Proposition \ref{Cohen-type}. But  Theorem   \ref{cohen2 even-odd2based} implies  \cite[Theorem 3.9]{devika2023} independently.


 \begin{theorem} \label{cohen2 even-odd3based}
  Let $\nu \notin \mathbb{Z}$ such that $\Re{(\nu)}\geq0$. Then, for any integer $ N$ such that 
 $  N\geq \lfloor\frac{\Re{(\nu)}+1}{2}\rfloor$, we have
				\begin{align*}
 &8\pi x^{\frac{\nu}{2}}\sum_{n=1}^{\infty}  n^{\nu/2} K_{\nu}(4\pi \sqrt{nx}) \sum_{d|n}d^{-\nu}  \cos \left(  {2\pi d \theta } \right)\sin \left( \frac{2\pi n\psi}{d}\right)=  \frac{\pi }{ 2 \sin\left(\frac{\pi \nu}{2}\right)}   x^\nu   \left(\zeta(1+\nu, \psi) - \zeta(1+\nu, 1-\psi)\right)\notag\\
  & \ \ \ \ \ \  +\frac{  1}{  2 \cos \left(\frac{\pi \nu}{2}\right)}   \left(\zeta(1, \psi) - \zeta(1, 1-\psi)\right) \left\{ \zeta(1-\nu, \theta)+\zeta(1-\nu,1- \theta)  \right\} \notag\\ 
    &\ \ \ \ \ \ +\frac{   1}{  2 \cos \left(\frac{\pi \nu}{2}\right)}  \sum_{j=1}^{N-1} x^{2j} \left( \zeta(2j+1,\psi)-\zeta(2j+1,1-\psi)         \right) \left\{\zeta(2j+1-\nu, \theta) + \zeta(2j+1-\nu, 1- \theta)  
     \right\} \notag\\
     &+\frac{  1}{2 \cos \left(\frac{\pi \nu}{2}\right)}x^{2N}\sum_{m,n\geq 0}^{\infty} \left\{  {(n+\psi)^{-\nu} } \left( \frac{((n+\psi)(m+\theta))^{\nu-2N+1}-x^{\nu-2N+1}}{ ((n+\psi)(m+\theta))^2-x^2} \right) 
   \right.\notag\\&\left.\ \ \ \ \ \ \ \  \ \ \  \ \ 
     -\  {(n+1-\psi)^{-\nu} } \left( \frac{((n+1-\psi)(m+\theta))^{\nu-2N+1}-x^{\nu-2N+1}}{ ((n+1-\psi)(m+\theta))^2-x^2} \right) 
    \right.\notag\\&\left.\ \ \ \ \ \ \ \  \ \ \  \ \ 
    +\  {(n+\psi)^{-\nu} } \left( \frac{((n+\psi)(m+1-\theta))^{\nu-2N+1}-x^{\nu-2N+1}}{ ((n+\psi)(m+1-\theta))^2-x^2} \right) 
       \right.\notag\\&\left.\ \  \ \ \ \ \ \  \ \ \  \ \  
      -\  {(n+1-\psi)^{-\nu} } \left( \frac{((n+1-\psi)(m+1-\theta))^{\nu-2N+1}-x^{\nu-2N+1}}{ ((n+1-\psi)(m+1-\theta))^2-x^2 } \right) 
   \right\},
\end{align*}
 provided   $x\notin  H_{\theta,\psi}$.

\end{theorem}
 \begin{theorem} \label{cohen2 even-odd4based}
Let $\nu \notin \mathbb{Z}$ such that $\Re{(\nu)}\geq0$. Then, for any integer $ N$ such that 
 $  N\geq \lfloor\frac{\Re{(\nu)}+1}{2}\rfloor$, we have
			 \begin{align*}
  & 8 \pi x^{\nu/2}\sum_{n=1}^{\infty}  n^{\nu/2} K_{\nu}(4\pi \sqrt{nx}) \sum_{d|n}d^{-\nu}  \sin \left(  {2\pi d\theta }  \right)\cos \left( \frac{2\pi n\psi}{d}\right)   =-   \frac{\pi }{ 2 \sin\left(\frac{\pi \nu}{2}\right)}      \left(\zeta(1-\nu,\theta) - \zeta(1-\nu, 1-\theta)\right)\notag\\
     & \ \ \ \ \ \ +\frac{   x^{\nu}}{  2 \cos \left(\frac{\pi \nu}{2}\right)}   \left(\zeta(1, \theta) - \zeta(1, 1-\theta)\right) \left\{ \zeta(1+\nu,\psi)+\zeta(1+\nu,1-\psi)  \right\} \notag\\ 
   & \ \ \ \ \ \ -\frac{   1}{  2 \cos \left(\frac{\pi \nu}{2}\right)}  \sum_{j=1}^N x^{2j-1} \left(\zeta(2j-\nu, \theta) - \zeta(2j-\nu, 1-\theta)\right) \left\{ \zeta(2j,\psi)+\zeta(2j,1-\psi) 
     \right\} \notag\\
  &\ \ \ \ \ \ -\frac{  1}{2 \cos \left(\frac{\pi \nu}{2}\right)}x^{2N+1}\sum_{m,n\geq 0}^{\infty} \left\{ \frac{(n+\psi)^{-\nu-1} }{(m+\theta)}\left( \frac{((n+\psi)(m+\theta))^{\nu-2N+1}-x^{\nu-2N+1}}{ ((n+\psi)(m+\theta))^2-x^2} \right) 
   \right.\notag\\&\left.\ \ \ \ \ \ \ \  \ \ \  \ \ \ \ \ \ \ \  
     +\ \frac{(n+1-\psi)^{-\nu-1} }{(m+\theta)}\left( \frac{((n+1-\psi)(m+\theta))^{\nu-2N+1}-x^{\nu-2N+1}}{ ((n+1-\psi)(m+\theta))^2-x^2} \right) 
    \right.\notag\\&\left.\ \ \ \ \ \ \ \  \ \ \  \ \ \ \ \ \ \ \ 
    -\ \frac{(n+\psi)^{-\nu-1} }{(m+1-\theta)}\left( \frac{((n+\psi)(m+1-\theta))^{\nu-2N+1}-x^{\nu-2N+1}}{ ((n+\psi)(m+1-\theta))^2-x^2} \right) 
       \right.\notag\\&\left.\ \  \ \ \ \ \ \  \ \ \  \ \ \ \ \ \ \ \   
      -\ \frac{(n+1-\psi)^{-\nu-1} }{(m+1-\theta)}\left( \frac{((n+1-\psi)(m+1-\theta))^{\nu-2N+1}-x^{\nu-2N+1}}{ ((n+1-\psi)(m+1-\theta))^2-x^2 } \right) 
   \right\},  
\end{align*}
 provided   $x\notin  H_{\theta,\psi}$.
\end{theorem}

To prove Theorem \ref{cohen2 even-odd3based}, one requires  \cite[Theorem 3.13]{devika2023} and Theorem \ref{oddcohen2 based}. On the other hand, Theorem \ref{cohen2 even-odd4based} is based on   \cite[Theorem 3.14]{devika2023} and Theorem \ref{oddcohen based}. Conversely, Theorem \ref{cohen2 even-odd3based} and Theorem \ref{cohen2 even-odd4based} imply   \cite[Theorem 3.13, Theorem 3.14]{devika2023}, respectively.


 
\section{Voronoi summation formula}\label{voronoi identities...}

 Here, we state the identities analogous to Entry \ref{entry1} which is the following:
\begin{theorem}\label{vor1.1}
Let $0< \alpha< \beta$ and $\alpha, \beta \notin \mathbb{Z}$. Let $f$ denote a function analytic inside a closed contour strictly containing $[\alpha, \beta ]$. Assume $0< \Re{(\nu)}<\frac{1}{2}$. Then, we have
\begin{align}\label{vv}
&\sum_{\alpha<j<\beta} \sum_{d|j}d^{-\nu}  \sin \left( {2\pi d \theta}\right) f(j)
    =- {(2\pi)^{\nu}} \Gamma(-\nu) \sin\left(\frac{\pi \nu}{2}\right)\left\{   \zeta(-\nu,\theta)- \zeta(-\nu,1-\theta)\right\} \int_\alpha^\beta    {f(t) }   \mathrm{d}t 
 \notag\\
 -&\pi \int_\alpha^\beta    \frac{f(t)}{t^{\frac{\nu}{2}}  } \sum_{d=1}^{\infty}   d^{-\frac{\nu}{2} }   \sum_{  m=0}^\infty  \left[ \left(m+\theta\right)^\frac{\nu}{2}\left\{  \left(   \frac{2}{\pi} K_{\nu}\left(4\pi  \sqrt {d\left(m+\theta\right)t}\ \right) + Y_{\nu}\left(4\pi \sqrt {d\left(m+\theta\right)t}\ \right)   \right) \sin \left(\frac{\pi \nu}{2}\right) 
        \right.\right.\notag\\&\left.\left.\ \ \ \ \ \ \ \ \ \ \ \ \ \ \ \ \ \ \ \ \ \ \ \ \ \ \ \ \ \ \  \ \ \ \ \ \ \ \ \ \ \ \ \ \ \  \ \ \ \ \  \ \ \  \ \ \ \  \ \  \ 
    -  J_{\nu}\left(4\pi  \sqrt {d\left(m+\theta\right)t}\ \right)  \cos \left(\frac{\pi \nu}{2}\right)  \right\}
    \right.\notag\\&\left.\  \ \ \ \ \  \ \ 
   - \left(m+1-\theta\right)^\frac{\nu}{2}\left\{  \left(   \frac{2}{\pi} K_{\nu}\left(4\pi  \sqrt {d\left(m+1-\theta\right)t}\ \right) + Y_{\nu}\left(4\pi \sqrt {d\left(m+1-\theta\right)t}\ \right)   \right) \sin \left(\frac{\pi \nu}{2}\right) 
        \right.\right.\notag\\&\left.\left.\ \ \ \ \ \ \ \ \ \ \ \ \ \ \ \ \ \ \ \ \ \ \ \ \ \ \ \ \ \ \  \ \ \ \ \ \ \ \ \ \ \ \ \ \ \  \ \ \ \ \  \ \ \  \ \ \ \  \ \  \ 
    -  J_{\nu}\left(4\pi  \sqrt {d\left(m+1-\theta\right)t}\ \right)  \cos \left(\frac{\pi \nu}{2}\right)  \right\}
     \right] \mathrm{d}t.
   \end{align}
           
 \end{theorem}
We demonstrate that Theorem \ref{vor1.1} is equivalent to the following theorem.  
  
\begin{theorem}\label{voro2}
 \cite[Theorem 4.4]{devika2023}
Let $0< \alpha< \beta$ and $\alpha, \beta \notin \mathbb{Z}$. Let $f$ denote a function analytic inside a closed contour strictly containing $[\alpha, \beta ]$. Assume that   $\chi$  is an odd primitive character modulo  $q$. For $0< \Re{(\nu)}<\frac{1}{2}$, we have  
\begin{align}\label{vv1}
&\frac{   q^{1+\frac{\nu}{2}}  } {\tau(\chi )}    \sum_{\alpha<j <\beta}    { {\sigma}_{-\nu, \chi }(j)} f(j)   =   \frac{ q^{1+\frac{\nu}{2}}}{\tau(\chi)} L(1+\nu,\chi) \int_\alpha^\beta     {f(t) } \mathrm{d}t
 + 2 \pi i \sum_{n=1}^{\infty}\bar{\sigma}_{-\nu, \bar{\chi }}(n) \ n^{\nu/2} \int_{\alpha} ^{  \beta   }f(t) (t)^{-\frac{\nu}{2}}  \notag\\ 
  &   \times     \left\{  \left(   \frac{2}{\pi} K_{\nu}\left(4\pi  \sqrt{\frac{nt}{q}}\ \right) + Y_{\nu}\left(4\pi  \sqrt{\frac{nt}{q}}\ \right)   \right)  \sin \left(\frac{\pi \nu}{2}\right)  
    -  J_{\nu}\left(4\pi  \sqrt{\frac{nt}{q}}\ \right)  \cos \left(\frac{\pi \nu}{2}\right)   \right\}  dt.
\end{align}
\end{theorem}

 \begin{theorem}\label{vor1.2}
 Let $0< \alpha< \beta$ and $\alpha, \beta \notin \mathbb{Z}$. Let $f$ denote a function analytic inside a closed contour strictly containing $[\alpha, \beta ]$. Assume $0< \Re{(\nu)}<\frac{1}{2}$. Then, we have
\begin{align*}
    &\sum_{\alpha<j<\beta} \sum_{d|j}d^{-\nu}  \sin \left(\frac{2\pi j \theta}{d}\right)\frac{f(j)}{j}= \frac{ \Gamma(\nu)\sin{\left( \frac{\pi \nu}{2} \right)}}{(2\pi )^\nu }\{ \zeta(\nu,\theta)-\zeta(\nu,1-\theta)\}\int_\alpha^\beta    \frac{f(t) }{t^{\nu+1}}  \mathrm{d}t   \notag\\
     + &\pi \int_{\alpha} ^{  \beta   }\frac{f(t)}{ t^{\frac{\nu}{2}+1} }     \sum_{r=1}^{\infty}   r^{\frac{\nu}{2}} \sum_{  m=0}^\infty  \left[ \left\{  \left(   \frac{2}{\pi}   \frac{K_{\nu}\left(4\pi  \sqrt {r\left(m+\theta \right)t} \ \right)}{\left(m+\theta \right)^{\frac{\nu}{2}}}  -   \frac{Y_{\nu}\left(4\pi  \sqrt {r\left(m+\theta \right)t} \ \right)}{\left(m+\theta \right)^{\frac{\nu}{2}}} \right) \sin \left(\frac{\pi \nu}{2}\right) 
    \right.\right.\notag\\&\left.\left.\ \ \ \ \ \ \ \ \ \ \ \ \ \ \ \ \ \ \ \ \ \ \ \ \ \ \ \ \ \ \  \ \ \ \ \ \ \ \ \ \ \ \ \ \ \  \ \ \ \ \  \ \ \  \ \ \ \  \ \  \ \   
  +    \frac{J_{\nu}\left(4\pi  \sqrt {r\left(m+\theta \right)t} \ \right)}{\left(m+\theta \right)^{\frac{\nu}{2}}}   \cos \left(\frac{\pi \nu}{2}\right)  \right\}
     \right.\notag\\&\left.  \ \ \ \ \  \ \  
     -\left\{  \left(   \frac{2}{\pi}   \frac{K_{\nu}\left(4\pi  \sqrt {r\left(m+1-\theta \right)t} \ \right)}{\left(m+1-\theta \right)^{\frac{\nu}{2}}}  -   \frac{Y_{\nu}\left(4\pi  \sqrt {r\left(m+1-\theta \right)t} \ \right)}{\left(m+1-\theta \right)^{\frac{\nu}{2}}} \right) \sin \left(\frac{\pi \nu}{2}\right) 
     \right.\right.\notag\\&\left.\left.\ \ \ \ \ \ \ \ \ \ \ \ \ \ \ \ \ \ \ \ \ \ \ \ \ \ \ \ \ \  \ \ \ \ \ \ \ \ \ \ \ \ \ \ \  \ \ \ \ \  \ \ \  \ \ \ \  
  +    \frac{J_{\nu}\left(4\pi  \sqrt {r\left(m+1-\theta \right)t} \ \right)}{\left(m+1-\theta \right)^{\frac{\nu}{2}}}   \cos \left(\frac{\pi \nu}{2}\right)  \right\} \right]\mathrm{d}t.
\end{align*}
 \end{theorem}

  Similar to Theorem \ref{vor1.1}, one can show that the equivalent version of Theorem \ref{vor1.2} is \cite[Theorem 4.3]{devika2023}. The identities in the next two theorems are analogous to Entry \ref{entry2}.
\begin{theorem}\label{vor1.3}
Let $0< \alpha< \beta$ and $\alpha, \beta \notin \mathbb{Z}$. Let $f$ denote a function analytic inside a closed contour strictly containing $[\alpha, \beta ]$. Assume $0< \Re{(\nu)}<\frac{1}{2}$, then
\begin{align}\label{vor1234}
&\sum_{\alpha<j<\beta} f(j) \sum_{d|j}d^{-\nu}  \cos \left( {2\pi d \theta}\right)
    = {(2\pi)^{\nu}} \Gamma(-\nu) \cos\left(\frac{\pi \nu}{2}\right)\left\{   \zeta(-\nu,\theta)+ \zeta(-\nu,1-\theta)\right\} \int_\alpha^\beta    {f(t) }   \mathrm{d}t 
 \notag\\
 +&\pi \int_\alpha^\beta    \frac{f(t)}{t^{\frac{\nu}{2}}  } \sum_{d=1}^{\infty}   d^{-\frac{\nu}{2} }   \sum_{  m=0}^\infty  \left[ \left(m+\theta\right)^\frac{\nu}{2}\left\{  \left(   \frac{2}{\pi} K_{\nu}\left(4\pi  \sqrt {d\left(m+\theta\right)t}\ \right) - Y_{\nu}\left(4\pi \sqrt {d\left(m+\theta\right)t}\ \right)   \right) \cos \left(\frac{\pi \nu}{2}\right) 
        \right.\right.\notag\\&\left.\left.\ \ \ \ \ \ \ \ \ \ \ \ \ \ \ \ \ \ \ \ \ \ \ \ \ \ \ \ \ \ \  \ \ \ \ \ \ \ \ \ \ \ \ \ \ \  \ \ \ \ \  \ \ \  \ \ \ \  \ \  \ 
    -  J_{\nu}\left(4\pi  \sqrt {d\left(m+\theta\right)t}\ \right)  \sin \left(\frac{\pi \nu}{2}\right)  \right\}
    \right.\notag\\&\left.\  \ \ \ \ \  \ \ 
   + \left(m+1-\theta\right)^\frac{\nu}{2}\left\{  \left(   \frac{2}{\pi} K_{\nu}\left(4\pi  \sqrt {d\left(m+1-\theta\right)t}\ \right) - Y_{\nu}\left(4\pi \sqrt {d\left(m+1-\theta\right)t}\ \right)   \right) \cos \left(\frac{\pi \nu}{2}\right) 
        \right.\right.\notag\\&\left.\left.\ \ \ \ \ \ \ \ \ \ \ \ \ \ \ \ \ \ \ \ \ \ \ \ \ \ \ \ \ \ \  \ \ \ \ \ \ \ \ \ \ \ \ \ \ \  \ \ \ \ \  \ \ \  \ \ \ \  \ \  \ 
    -  J_{\nu}\left(4\pi  \sqrt {d\left(m+1-\theta\right)t}\ \right)  \sin \left(\frac{\pi \nu}{2}\right)  \right\}
     \right] \mathrm{d}t.
   \end{align}
           
 \end{theorem}

 Theorem \ref{vor1.3} is proved using  the following  \cite[Theorem 4.2]{devika2023} and Proposition \ref{vorlemma1}. But the following identity can be deduced directly using Theorem \ref{vor1.3}.
\begin{theorem}\label{vore2}
\cite[Theorem 4.2]{devika2023}
Let $0< \alpha< \beta$ and $\alpha, \beta \notin \mathbb{Z}$. Let $f$ denote a function analytic inside a closed contour strictly containing $[\alpha, \beta ]$. Assume that   $\chi$  is a  non-principal, even primitive character modulo $ q$.  For  $0< \Re{(\nu)}<\frac{1}{2}$, we have  
\begin{align}\label{fann}
&\frac{   q^{1+\frac{\nu}{2}}  } {\tau(\chi )}    \sum_{\alpha<j <\beta}    { {\sigma}_{-\nu, \chi }(j)} f(j)   =   \frac{ q^{1+\frac{\nu}{2}}}{\tau(\chi)} L(1+\nu,\chi) \int_\alpha^\beta     {f(t) } \mathrm{d}t
 + 2 \pi \sum_{n=1}^{\infty}\bar{\sigma}_{-\nu, \bar{\chi }}(n) \ n^{\nu/2} \int \frac{f(t) }{t^{\frac{\nu}{2}}}  \notag\\ 
  &   \times     \left\{  \left(   \frac{2}{\pi} K_{\nu}\left(4\pi  \sqrt{\frac{nt}{q}}\ \right) - Y_{\nu}\left(4\pi  \sqrt{\frac{nt}{q}}\ \right)   \right) \cos \left(\frac{\pi \nu}{2}\right) 
    -  J_{\nu}\left(4\pi  \sqrt{\frac{nt}{q}}\ \right)  \sin \left(\frac{\pi \nu}{2}\right)  \right\} \mathrm{d}t.
\end{align}

\end{theorem}

\begin{theorem}\label{vor1.4}
Let $0< \alpha< \beta$ and $\alpha, \beta \notin \mathbb{Z}$. Let $f$ denote a function analytic inside a closed contour strictly containing $[\alpha, \beta ]$. Assume $0< \Re{(\nu)}<\frac{1}{2}$, then
\begin{align*}
 &  \sum_{\alpha<j <\beta}  \sum_{d/j}  {d} ^{-\nu} \cos\left(\frac{2 \pi j \theta }{d}\right)f(j)   =\frac{ \Gamma(\nu)\cos{\left( \frac{\pi \nu}{2} \right)}}{(2\pi )^\nu }\{ \zeta(\nu,\theta)+\zeta(\nu,1-\theta)\} \int_{\alpha} ^{  \beta   }\frac{f(t)}{ t^{{\nu} } } \mathrm{d}t
 \notag\\
 + &\pi \int_{\alpha} ^{  \beta   }\frac{f(t)}{ t^{\frac{\nu}{2}} }     \sum_{r=1}^{\infty}   r^{\frac{\nu}{2}} \sum_{  m=0}^\infty  \left[ \left\{  \left(   \frac{2}{\pi}   \frac{K_{\nu}\left(4\pi  \sqrt {r\left(m+\theta \right)t} \ \right)}{\left(m+\theta \right)^{\frac{\nu}{2}}}  -   \frac{Y_{\nu}\left(4\pi  \sqrt {r\left(m+\theta \right)t} \ \right)}{\left(m+\theta \right)^{\frac{\nu}{2}}} \right) \cos \left(\frac{\pi \nu}{2}\right) 
    \right.\right.\notag\\&\left.\left.\ \ \ \ \ \ \ \ \ \ \ \ \ \ \ \ \ \ \ \ \ \ \ \ \ \ \ \ \ \ \  \ \ \ \ \ \ \ \ \ \ \ \ \ \ \  \ \ \ \ \  \ \ \  \ \ \ \  \ \  \ \   
  -    \frac{J_{\nu}\left(4\pi  \sqrt {r\left(m+\theta \right)t} \ \right)}{\left(m+\theta \right)^{\frac{\nu}{2}}}   \sin \left(\frac{\pi \nu}{2}\right)  \right\}
     \right.\notag\\&\left.  \ \ \ \ \  \ \  
     +\left\{  \left(   \frac{2}{\pi}   \frac{K_{\nu}\left(4\pi  \sqrt {r\left(m+1-\theta \right)t} \ \right)}{\left(m+1-\theta \right)^{\frac{\nu}{2}}}  -   \frac{Y_{\nu}\left(4\pi  \sqrt {r\left(m+1-\theta \right)t} \ \right)}{\left(m+1-\theta \right)^{\frac{\nu}{2}}} \right) \cos \left(\frac{\pi \nu}{2}\right) 
     \right.\right.\notag\\&\left.\left.\ \ \ \ \ \ \ \ \ \ \ \ \ \ \ \ \ \ \ \ \ \ \ \ \ \ \ \ \ \  \ \ \ \ \ \ \ \ \ \ \ \ \ \ \  \ \ \ \ \  \ \ \  \ \ \ \  
  -    \frac{J_{\nu}\left(4\pi  \sqrt {r\left(m+1-\theta \right)t} \ \right)}{\left(m+1-\theta \right)^{\frac{\nu}{2}}}   \sin \left(\frac{\pi \nu}{2}\right)  \right\} \right]\mathrm{d}t.
           \end{align*}
 \end{theorem}

Analogous to Theorem \ref{vor1.3}, the equivalent version of Theorem \ref{vor1.4} is \cite[Theorem 4.1]{devika2023}.
     Next, we state the identities involving two trigonometric functions,  which are the following:

\begin{theorem}\label{vor1.6}
Let $0< \alpha< \beta$ and $\alpha, \beta \notin \mathbb{Z}$. Let $f$ denote a function analytic inside a closed contour strictly containing $[\alpha, \beta ]$. Assume $0< \Re{(\nu)}<\frac{1}{2}$, then 
\begin{align}\label{1234}
   &\sum_{\alpha<j <\beta}  \sum_{d/j}  {d} ^{-\nu} \cos \left( {2\pi d \theta}\right)\cos\left(\frac{2 \pi j \psi }{d}\right)f(j) \notag\\
 =\frac{\pi}{2}    &   \int_{\alpha} ^{  \beta   }\frac{f(t) }{t^{\frac{\nu}{2}}}\sum_{m,n\geq 0}^{\infty}  \left[ \left( \frac{m+\theta}{n+\psi}\right)^{\nu/2} \left\{  \left(   \frac{2}{\pi} K_{\nu}(4\pi   \sqrt{(n+\psi)(m+\theta)t}) - Y_{\nu}(4\pi    \sqrt{(n+\psi)(m+\theta)t})\right) \cos \left(\frac{\pi \nu}{2}\right)
 \right.\right.\notag\\&\left.\left.\ \ \ \ \ \ \ \ \ \ \ \ \ \ \ \ \ \ \ \ \ \ \ \ \ \ \ \ \  \ \ \ \ \ \ \ \ \ \ \ \ \ \ \  \ \ \ \ \  \ \ \  \ \ \ \  
  -  J_{\nu}(4\pi    \sqrt{(n+\psi)(m+\theta)t})  \sin \left(\frac{\pi \nu}{2}\right)   
    \right\}   
   \right.\notag\\&\left. 
   +\left( \frac{m+\theta}{n+1-\psi}\right)^{\nu/2} \left\{  \left(   \frac{2}{\pi} K_{\nu}(4\pi   \sqrt{(n+1-\psi)(m+\theta)t}) - Y_{\nu}(4\pi    \sqrt{(n+1-\psi)(m+\theta)t})\right) \cos \left(\frac{\pi \nu}{2}\right)
 \right.\right.\notag\\&\left.\left.\ \ \ \ \ \ \ \ \ \ \ \ \ \ \ \ \ \ \ \ \ \ \ \ \ \ \ \ \  \ \ \ \ \ \ \ \ \ \ \ \ \ \ \  \ \ \ \ \ \ \ \  \  
  -  J_{\nu}(4\pi    \sqrt{(n+1-\psi)(m+\theta)t})  \sin \left(\frac{\pi \nu}{2}\right)   
    \right\}   
   \right.\notag\\&\left.  
   +\left( \frac{m+1-\theta}{n+\psi}\right)^{\nu/2} \left\{  \left(   \frac{2}{\pi} K_{\nu}(4\pi   \sqrt{(n+\psi)(m+1-\theta)t}) - Y_{\nu}(4\pi    \sqrt{(n+\psi)(m+1-\theta)t})\right) \cos \left(\frac{\pi \nu}{2}\right)
 \right.\right.\notag\\&\left.\left.\ \ \ \ \ \ \ \ \ \ \ \ \ \ \ \ \ \ \ \ \ \ \ \ \ \ \ \ \  \ \ \ \ \ \ \ \ \ \ \ \ \ \ \  \ \ \ \ \ \ \  \ \ 
  -  J_{\nu}(4\pi    \sqrt{(n+\psi)(m+1-\theta)t})  \sin \left(\frac{\pi \nu}{2}\right)   
    \right\}   
   \right.\notag\\&\left.  
   +\left( \frac{m+1-\theta}{n+1-\psi}\right)^{\nu/2} \left\{  \left(   \frac{2}{\pi} K_{\nu}(4\pi   \sqrt{(n+1-\psi)(m+1-\theta)t}) - Y_{\nu}(4\pi    \sqrt{(n+1-\psi)(m+1-\theta)t})\right) \cos \left(\frac{\pi \nu}{2}\right)
 \right.\right.\notag\\&\left.\left.\ \ \ \ \ \ \ \ \ \ \ \ \ \ \ \ \ \ \ \ \ \ \ \ \ \ \ \ \  \ \ \ \ \ \ \ \ \ \ \ \ \ \ \  \ \ \  
  -  J_{\nu}(4\pi    \sqrt{(n+1-\psi)(m+1-\theta)t})  \sin \left(\frac{\pi \nu}{2}\right)   
    \right\}  
 \right]\mathrm{d}t. 
 \end{align}
\end{theorem}
 For deriving Theorem \ref{vor1.6}, one requires   \cite[Theorem 4.5]{devika2023}, Theorem \ref{vor1.3}, Theorem \ref{vor1.4} and Proposition \ref{vorlemma1}. But  Theorem   \ref{vor1.6} implies  \cite[Theorem 4.5]{devika2023} independently.

\begin{theorem}\label{vor1.5}
Let $0< \alpha< \beta$ and $\alpha, \beta \notin \mathbb{Z}$. Let $f$ denote a function analytic inside a closed contour strictly containing $[\alpha, \beta ]$. Assume $0< \Re{(\nu)}<\frac{1}{2}$, then 
\begin{align*}
   &\sum_{\alpha<j <\beta}  \sum_{d/j}  {d} ^{-\nu} \sin \left( {2\pi d \theta}\right)\sin \left(\frac{2 \pi j \psi }{d}\right)\frac{f(j)}{j}  \notag\\
   =\frac{\pi}{2}    &   \int_{\alpha} ^{  \beta   }\frac{f(t) }{t^{\frac{\nu}{2}+1}} \sum_{m,n\geq 0}^{\infty}  \left[ \left( \frac{m+\theta}{n+\psi}\right)^{\nu/2} \left\{  \left(   \frac{2}{\pi} K_{\nu}(4\pi   \sqrt{(n+\psi)(m+\theta)t}) + Y_{\nu}(4\pi    \sqrt{(n+\psi)(m+\theta)t})\right) \cos \left(\frac{\pi \nu}{2}\right)
 \right.\right.\notag\\&\left.\left.\ \ \ \ \ \ \ \ \ \ \ \ \ \ \ \ \ \ \ \ \ \ \ \ \ \ \ \ \  \ \ \ \ \ \ \ \ \ \ \ \ \ \ \  \ \ \ \ \  \ \ \  \ \ \ \  
  + J_{\nu}(4\pi    \sqrt{(n+\psi)(m+\theta)t})  \sin \left(\frac{\pi \nu}{2}\right)   
    \right\}    
   \right.\notag\\&\left. 
   -\left( \frac{m+\theta}{n+1-\psi}\right)^{\nu/2} \left\{  \left(   \frac{2}{\pi} K_{\nu}(4\pi   \sqrt{(n+1-\psi)(m+\theta)t}) + Y_{\nu}(4\pi    \sqrt{(n+1-\psi)(m+\theta)t})\right) \cos \left(\frac{\pi \nu}{2}\right)
 \right.\right.\notag\\&\left.\left.\ \ \ \ \ \ \ \ \ \ \ \ \ \ \ \ \ \ \ \ \ \ \ \ \ \ \ \ \  \ \ \ \ \ \ \ \ \ \ \ \ \ \ \  \ \ \ \ \ \ \ \  \  
  +  J_{\nu}(4\pi    \sqrt{(n+1-\psi)(m+\theta)t})  \sin \left(\frac{\pi \nu}{2}\right)   
    \right\}   
   \right.\notag\\&\left.  
   -\left( \frac{m+1-\theta}{n+\psi}\right)^{\nu/2} \left\{  \left(   \frac{2}{\pi} K_{\nu}(4\pi   \sqrt{(n+\psi)(m+1-\theta)t}) + Y_{\nu}(4\pi    \sqrt{(n+\psi)(m+1-\theta)t})\right) \cos \left(\frac{\pi \nu}{2}\right)
 \right.\right.\notag\\&\left.\left.\ \ \ \ \ \ \ \ \ \ \ \ \ \ \ \ \ \ \ \ \ \ \ \ \ \ \ \ \  \ \ \ \ \ \ \ \ \ \ \ \ \ \ \  \ \ \ \ \ \ \  \ \ 
  +  J_{\nu}(4\pi    \sqrt{(n+\psi)(m+1-\theta)t})  \sin \left(\frac{\pi \nu}{2}\right)   
    \right\}  
   \right.\notag\\&\left.  
   +\left( \frac{m+1-\theta}{n+1-\psi}\right)^{\nu/2} \left\{  \left(   \frac{2}{\pi} K_{\nu}(4\pi   \sqrt{(n+1-\psi)(m+1-\theta)t}) + Y_{\nu}(4\pi    \sqrt{(n+1-\psi)(m+1-\theta)t})\right) \cos \left(\frac{\pi \nu}{2}\right)
 \right.\right.\notag\\&\left.\left.\ \ \ \ \ \ \ \ \ \ \ \ \ \ \ \ \ \ \ \ \ \ \ \ \ \ \ \ \  \ \ \ \ \ \ \ \ \ \ \ \ \ \ \  \ \ \  
  +  J_{\nu}(4\pi    \sqrt{(n+1-\psi)(m+1-\theta)t})  \sin \left(\frac{\pi \nu}{2}\right)   
    \right\}   
 \right]\mathrm{d}t.
   \end{align*}
 \end{theorem}
 We remark that Theorem \ref{vor1.5} is equivalent to  
  \cite[Theorem 4.7]{devika2023}.
\begin{theorem}\label{vor1.7}
Let $0< \alpha< \beta$ and $\alpha, \beta \notin \mathbb{Z}$. Let $f$ denote a function analytic inside a closed contour strictly containing $[\alpha, \beta ]$. Assume $0< \Re{(\nu)}<\frac{1}{2}$, then 
\begin{align*}
    &\sum_{\alpha<j <\beta}  \sum_{d/j}  {d} ^{-\nu}  \cos \left( 2\pi d \theta\right)\sin \left( \frac{2\pi j \psi}{d}\right) \frac{f(j)}{j}   \notag\\
    =\frac{\pi}{2}    &   \int_{\alpha} ^{  \beta   }\frac{f(t) }{t^{\frac{\nu}{2}+1}}    \sum_{m,n\geq 0}^{\infty}  \left[ \left( \frac{m+\theta}{n+\psi}\right)^{\nu/2} \left\{  \left(   \frac{2}{\pi} K_{\nu}(4\pi   \sqrt{(n+\psi)(m+\theta)t}) - Y_{\nu}(4\pi    \sqrt{(n+\psi)(m+\theta)t})\right) \sin \left(\frac{\pi \nu}{2}\right)
 \right.\right.\notag\\&\left.\left.\ \ \ \ \ \ \ \ \ \ \ \ \ \ \ \ \ \ \ \ \ \ \ \ \ \ \ \ \  \ \ \ \ \ \ \ \ \ \ \ \ \ \ \  \ \ \ \ \  \ \ \  \ \ \ \  
  + J_{\nu}(4\pi    \sqrt{(n+\psi)(m+\theta)t})  \cos \left(\frac{\pi \nu}{2}\right)   
    \right\}   
   \right.\notag\\&\left. 
   -\left( \frac{m+\theta}{n+1-\psi}\right)^{\nu/2} \left\{  \left(   \frac{2}{\pi} K_{\nu}(4\pi   \sqrt{(n+1-\psi)(m+\theta)t}) - Y_{\nu}(4\pi    \sqrt{(n+1-\psi)(m+\theta)t})\right) \sin \left(\frac{\pi \nu}{2}\right)
 \right.\right.\notag\\&\left.\left.\ \ \ \ \ \ \ \ \ \ \ \ \ \ \ \ \ \ \ \ \ \ \ \ \ \ \ \ \  \ \ \ \ \ \ \ \ \ \ \ \ \ \ \  \ \ \ \ \ \ \ \  \  
  +  J_{\nu}(4\pi    \sqrt{(n+1-\psi)(m+\theta)t})  \cos \left(\frac{\pi \nu}{2}\right)   
    \right\}    
   \right.\notag\\&\left.  
   +\left( \frac{m+1-\theta}{n+\psi}\right)^{\nu/2} \left\{  \left(   \frac{2}{\pi} K_{\nu}(4\pi   \sqrt{(n+\psi)(m+1-\theta)t}) - Y_{\nu}(4\pi    \sqrt{(n+\psi)(m+1-\theta)t})\right) \sin \left(\frac{\pi \nu}{2}\right)
 \right.\right.\notag\\&\left.\left.\ \ \ \ \ \ \ \ \ \ \ \ \ \ \ \ \ \ \ \ \ \ \ \ \ \ \ \ \  \ \ \ \ \ \ \ \ \ \ \ \ \ \ \  \ \ \ \ \ \ \  \ \ 
  +  J_{\nu}(4\pi    \sqrt{(n+\psi)(m+1-\theta)t})  \cos \left(\frac{\pi \nu}{2}\right)   
    \right\}  
   \right.\notag\\&\left.  
   -\left( \frac{m+1-\theta}{n+1-\psi}\right)^{\nu/2} \left\{  \left(   \frac{2}{\pi} K_{\nu}(4\pi   \sqrt{(n+1-\psi)(m+1-\theta)t}) - Y_{\nu}(4\pi    \sqrt{(n+1-\psi)(m+1-\theta)t})\right) \sin \left(\frac{\pi \nu}{2}\right)
 \right.\right.\notag\\&\left.\left.\ \ \ \ \ \ \ \ \ \ \ \ \ \ \ \ \ \ \ \ \ \ \ \ \ \ \ \ \  \ \ \ \ \ \ \ \ \ \ \ \ \ \ \  \ \ \  
  -  J_{\nu}(4\pi    \sqrt{(n+1-\psi)(m+1-\theta)t})  \cos \left(\frac{\pi \nu}{2}\right)   
    \right\}   
 \right]\mathrm{d}t.
   \end{align*}
\end{theorem}
\begin{theorem}\label{vor1.8}
Let $0< \alpha< \beta$ and $\alpha, \beta \notin \mathbb{Z}$. Let $f$ denote a function analytic inside a closed contour strictly containing $[\alpha, \beta ]$. Assume $0< \Re{(\nu)}<\frac{1}{2}$, then 
\begin{align*}
    &\sum_{\alpha<j <\beta}  \sum_{d/j}  {d} ^{-\nu}    \sin \left( 2\pi d \theta\right)\cos \left( \frac{2\pi j \psi}{d}\right)  f(j)  \notag\\
    =-\frac{\pi}{2}    &   \int_{\alpha} ^{  \beta   }\frac{f(t) }{t^{\frac{\nu}{2}}} \sum_{m,n\geq 0}^{\infty}  \left[ \left( \frac{m+\theta}{n+\psi}\right)^{\nu/2} \left\{  \left(   \frac{2}{\pi} K_{\nu}(4\pi   \sqrt{(n+\psi)(m+\theta)t}) + Y_{\nu}(4\pi    \sqrt{(n+\psi)(m+\theta)t})\right) \sin \left(\frac{\pi \nu}{2}\right)
 \right.\right.\notag\\&\left.\left.\ \ \ \ \ \ \ \ \ \ \ \ \ \ \ \ \ \ \ \ \ \ \ \ \ \ \ \ \  \ \ \ \ \ \ \ \ \ \ \ \ \ \ \  \ \ \ \ \  \ \ \  \ \ \ \  
  -  J_{\nu}(4\pi    \sqrt{(n+\psi)(m+\theta)t})  \cos \left(\frac{\pi \nu}{2}\right)   
    \right\}   
   \right.\notag\\&\left. 
   +\left( \frac{m+\theta}{n+1-\psi}\right)^{\nu/2} \left\{  \left(   \frac{2}{\pi} K_{\nu}(4\pi   \sqrt{(n+1-\psi)(m+\theta)t}) + Y_{\nu}(4\pi    \sqrt{(n+1-\psi)(m+\theta)t})\right) \sin \left(\frac{\pi \nu}{2}\right)
 \right.\right.\notag\\&\left.\left.\ \ \ \ \ \ \ \ \ \ \ \ \ \ \ \ \ \ \ \ \ \ \ \ \ \ \ \ \  \ \ \ \ \ \ \ \ \ \ \ \ \ \ \  \ \ \ \ \ \ \ \  \  
  -  J_{\nu}(4\pi    \sqrt{(n+1-\psi)(m+\theta)t})  \cos \left(\frac{\pi \nu}{2}\right)   
    \right\}  
   \right.\notag\\&\left.  
   -\left( \frac{m+1-\theta}{n+\psi}\right)^{\nu/2} \left\{  \left(   \frac{2}{\pi} K_{\nu}(4\pi   \sqrt{(n+\psi)(m+1-\theta)t}) + Y_{\nu}(4\pi    \sqrt{(n+\psi)(m+1-\theta)t})\right) \sin \left(\frac{\pi \nu}{2}\right)
 \right.\right.\notag\\&\left.\left.\ \ \ \ \ \ \ \ \ \ \ \ \ \ \ \ \ \ \ \ \ \ \ \ \ \ \ \ \  \ \ \ \ \ \ \ \ \ \ \ \ \ \ \  \ \ \ \ \ \ \  \ \ 
  -  J_{\nu}(4\pi    \sqrt{(n+\psi)(m+1-\theta)t})  \cos \left(\frac{\pi \nu}{2}\right)   
    \right\}  
   \right.\notag\\&\left.  
   -\left( \frac{m+1-\theta}{n+1-\psi}\right)^{\nu/2} \left\{  \left(   \frac{2}{\pi} K_{\nu}(4\pi   \sqrt{(n+1-\psi)(m+1-\theta)t}) + Y_{\nu}(4\pi    \sqrt{(n+1-\psi)(m+1-\theta)t})\right) \sin \left(\frac{\pi \nu}{2}\right)
 \right.\right.\notag\\&\left.\left.\ \ \ \ \ \ \ \ \ \ \ \ \ \ \ \ \ \ \ \ \ \ \ \ \ \ \ \ \  \ \ \ \ \ \ \ \ \ \ \ \ \ \ \  \ \ \  
  -  J_{\nu}(4\pi    \sqrt{(n+1-\psi)(m+1-\theta)t})  \cos \left(\frac{\pi \nu}{2}\right)   
    \right\}   
 \right]\mathrm{d}t.
   \end{align*}
\end{theorem}
 
 To prove Theorem \ref{vor1.7}, one requires  \cite[Theorem 4.10]{devika2023} and Theorem \ref{vor1.2}. On the other hand, Theorem \ref{vor1.8} is based on   \cite[Theorem 4.9]{devika2023} and Theorem \ref{vor1.1}. Conversely, Theorem \ref{vor1.7} and Theorem \ref{vor1.8} directly imply    \cite[Theorem 4.10, Theorem 4.9]{devika2023}, respectively.

 \section{Preliminaries}\label{preliminary}
			
			We begin this section by recalling   some important results which we will use throughout the paper. The functional equation of $\zeta(s)$ is given by  \cite[p. 259]{MR0434929}
			\begin{align}\label{1st_use}
				\Gamma(s)  \zeta(s) & = \frac{\pi^{s} \zeta(1-s)}{2^{1-s} \cos\left(\frac{\pi s}{2}\right) }. 
			\end{align}
	Next, we state the functional equation for $L(s, \chi)$ \cite[Corrolary 10.9, p.~333]{MR2378655}
			\begin{align}\label{ll(s)}
				L(s, \chi) & =i^{-\kappa} \frac{\tau(\chi)}{ \pi}\left(\frac{(2\pi)}{q}\right)^{s} \Gamma(1-s) \sin \frac{\pi (s+\kappa)}{2} L(1-s, \bar{\chi}),  
			\end{align}
			where 
			\begin{align*}
				\kappa=\kappa(\chi)=\begin{cases}
					&  0, \ \ \mbox{if} \ \chi(-1)=1,\\
					&  1,\  \ \mbox{if} \ \chi(-1)=-1.\\
				\end{cases}
			\end{align*}
Now replacing $s$ by $s-z$ in \eqref{ll(s)}, we obtain
			\begin{align*}
				L(s-z, \chi) & =i^{-\kappa} \frac{\tau(\chi)}{\pi}\left(\frac{(2\pi)}{q}\right)^{s-z} \Gamma(1+z-s) \sin \frac{\pi (s+\kappa-z) }{2} L(1+z-s, \bar{\chi}). 
			\end{align*}
			The above equation can be rewritten as
   \begin{align}\label{exact l}
       \Gamma(1+z-s)L(1+z-s, \bar{\chi})= i^{\kappa}\frac{\pi}{\tau(\chi)}\left( \frac{q}{2\pi}\right)^{s-z}\frac{L(s-z, \chi)}{\sin{\pi(\frac{s+\kappa-z}{2})}}.
   \end{align}



The Hurwitz zeta function satisfies the following functional equations \cite[p. 587]{MR3283193}
\begin{align}\label{hur1}
  &\sum_{r=1}^q \zeta\left(s,\frac{r}{q}\right) \cos\left(\frac{2\pi rh}{q}\right)= \frac{q\Gamma(1-s)}{\left(2\pi q \right)^{1-s}} \sin\left(\frac{\pi s}{2}\right) \left\{ \zeta\left(1-s,\frac{h}{q}\right)+\zeta\left(1-s, 1-{\frac{h}{q}}\right)\right\},\\
  &\sum_{r=1}^q \zeta\left(s,\frac{r}{q}\right) \sin\left(\frac{2\pi rh}{q}\right)= \frac{q\Gamma(1-s)}{\left(2\pi q \right)^{1-s}}  \cos\left(\frac{\pi s}{2}\right)\left\{ \zeta\left(1-s,\frac{h}{q}\right)-\zeta\left(1-s, 1-{\frac{h}{q}}\right)\right\}.\label{hur2}
\end{align}
{  Hurwitz zeta function also satisfies
  \cite[p. 264]{MR0434929} \begin{align}\label{SV}
            \zeta(-n,\theta)=&-\frac{B_{n+1}(\theta)}{n+1},
             \end{align}
  for each $n\ge 0$; where $B_{n}(\theta)$ is a  Bernoulli polynomial defined as follows  \cite[p. 264]{MR0434929}
  \begin{align}
    \frac{z e^{\theta z}}{e^z-1}=\sum_{n=0}^\infty \frac{B_{n}(\theta)}{n!}z^n, \text{ for }|z|<2\pi,
\end{align}
for any $\theta \in \mathbb{C}$ 
             and we have the relation \cite[p. 274]{MR0434929}
             \begin{align}
             B_{n}(1-\theta)=&(-1)^n B_{n}(\theta) \ \ \text{ for every }n \geq 0. \label{bernoulli2}
         \end{align}

         }


			We will also note that  \cite[p.69, p.71]{MR1790423}
			\begin{align}\label{both}
				\tau(\chi)\tau( \bar{\chi})= 
				\begin{cases}
					&-q, \ \text{for odd primitive }\chi \ mod \ q , \\
					& q, \ \text{for even non principal primitive }\chi \ mod \ q, 
				\end{cases} 
			\end{align}
			 	and the fact
			\begin{align}\label{prop}
				\sum_{\substack{\chi \ mod \ q\\ \chi \ odd}}\chi(a) \bar{\chi}(h)=\begin{cases}
					& \pm \frac{\phi(q)}{2} \ \ \mbox{if} \ h \equiv \pm a \ (mod \ q)\\
					&0 \ \ \mbox{otherwise ;}
				\end{cases}\\
   \sum_{\substack{\chi \ mod \ q\\ \chi \ even}}\chi(a) \bar{\chi}(h)=\begin{cases}
					&   \frac{\phi(q)}{2} \ \ \mbox{if} \ h \equiv \pm a \ (mod \ q)\\
					&0 \ \ \mbox{otherwise.}
				\end{cases} \label{prop1}
			\end{align}	
			Here we would like to mention another identity \cite[Lemma 2.5]{MR3181548} namely 
			\begin{align}\label{cos}
				 & \sin \left(\frac{2\pi h d}{q}\right)=\frac{1}{i \phi(q)} \sum_{\substack{\chi \ mod \ q\\    \chi \ odd}} \chi(d) \tau(\bar{\chi}) \chi(h), \\
					& \cos \left(\frac{2\pi h d}{q}\right)=\frac{1}{ \phi(q)} \sum_{\substack{\chi \ mod \ q\\    \chi \ even}} \chi(d) \tau(\bar{\chi}) \chi(h),    
				 \label{sin}\end{align}
			whenever $(d,q)=(h,q)=1$. The factorization theorem for Gauss sum $\tau(\chi)$ in \eqref{Gauss} is as follows \cite[p. 65]{MR1790423}
\begin{align}\label{gauss}
    \chi(n)\tau(\Bar{\chi})=\sum_{h=1}^{q-1}\Bar{\chi}(h)e^{2\pi i nh/q},
\end{align}
 for any character modulo $\chi$ modulo $q$.

{ 
 
\par 
Next, we recall the Mellin transform of a locally integrable function $f(x)$ on $(0, \infty)$ is defined by
			\begin{align}\label{Mellin}
				\mathcal{M}[f; s] = F(s) =\int_{0}^{\infty} f(t) \ t^{s-1} dt,
			\end{align}
			provided the integral converges. 
			 
For example, $\mathcal{M}[e^{-t}; s]=\Gamma(s)$ for $\Re(s)>0$, and we have the corresponding Mellin's inversion formula for any $\Re{(y)}>0$
\begin{align*}
   e^{-y}=  \frac{1}{2\pi i}\int_{(c)}\Gamma(s)y^{-s}\mathrm{d}s, 
\end{align*}
where $c>0$. Here the notation $(c)$ denotes the vertical line $[c-i\infty, c+i\infty]$. The functional relations for $\Gamma(s)$ are given by \cite[p.~73]{MR1790423}
				\begin{align} \label{relation of gamma}
					 \Gamma(s+1)=&s \Gamma(s), \ \ \ \ \ \ \ \ \Gamma(s)\Gamma\left(s+\frac{1}{2}\right)=2^{1-2s}\sqrt{\pi}\Gamma(2s),
     \\
        & \Gamma(s)\Gamma(1-s)=\frac{\pi}{\sin (\pi s
)}.\label{reflection formula for gamma add}
   \end{align}
 
 \begin{lemma}\label{lem4} 
 \cite[p.~24, Lemma 5.5]{devika2023} 
				Let $n\geq0$ be any integer and $t>0$ be any real number. Then
				\begin{align*}
				 	\frac{1} {2\pi i}\int_{(c)}  { \Gamma(s+n) \Gamma(a-s)} t^{-s}\mathrm{d}s= \frac{ \Gamma(a+n)     }{(1+t)^{a+n} }t^{n},
				\end{align*}
				for $0<c<\Re(a)$.
			\end{lemma}
The following lemma states the asymptotic behaviour of $\Gamma(s)$.
			\begin{lemma}\label{Gamma}
				\cite[p.~38]{MR0364103} In a vertical strip, for s=$\sigma+it$ with $a\leq \sigma \leq b$ and $|t|\geq 1$,
				\begin{align*}
					|\Gamma(s)|=(2\pi)^\frac{1}{2}|t|^{\sigma-\frac{1}{2}}\exp^{-\frac{1}{2}\pi |t|}\left( 1+O\left(\frac{1}{|t|}\right)\right).
				\end{align*}
			\end{lemma}
	 Now, we  record a few important results related to the modified $K$-Bessel function $K_\nu(x)$ defined by \eqref{Kbessel}.
			\begin{lemma}\label{eq:bessel}
   \cite[p.~10, Lemma 3.3]{MR4570432}
				Let $\nu \in \mathbf{C}$. For any $c> \max\{0,-\Re(\nu)\}$, we have
				\begin{align*}
					t^{\frac{\nu}{2}}K_{\nu}(a\sqrt{tx})=\frac{1}{2}\left( \frac{2}{a\sqrt{x}}\right)^\nu \frac{1}{2\pi i}
					\int_{(c)}\Gamma(s) \Gamma(s+\nu)\left( \frac{4}{a^2x}\right)^st^{-s} ds.   \end{align*}
			\end{lemma}
			 
}

 \section{Proof of  Main Results   when $z \in  \mathbb{Z}$ }\label{proof of integer nu} Throughout this section, let  $X=4/(a^2x)$ where $a,x>0$ be positive  real number.
 \par

    \begin{proof}[Theorem \rm{\ref{odd1_based}} and its equivalence with Theorem \rm{\ref{M1}}][]
First, we recall that $\Re(\nu)>0$. The double series on the right-hand side of the identity \eqref{p7} converges absolutely and uniformly on any compact interval for $\theta\in (0,1)$. Since the summands in the right-hand side of \eqref{p7} are continuous functions, the series converges to  a continuous function of $\theta$. Therefore, it is sufficient to prove the identity \eqref{p7} for $\theta=h/q$, where $q$ is prime and $0<h<q.$ If the identity holds on the dense subset of fractions, then by continuity, the identity holds for all values of $\theta$.  Employing Lemma \ref{eq:bessel} with $t=n$ and subsequently interchanging the summation and integration, we get for an even integer $k\geq 0$, 
\begin{align}\label{g1}
				\sum_{n=1}^{\infty}  \sum_{d|n}d^k \sin \left(   \frac{ 2\pi dh}{q} \right) n^{\nu/2} &K_{\nu}(a\sqrt{nx})
     = \frac{1}{2} \left(\frac{2}{a \sqrt{x}}\right)^{ \nu}  \frac{1}{2\pi i} \int_{(c)}  \Gamma(s) \Gamma(s + \nu) \left( \frac{4}{a^2 x } \right)^s  \sum_{n=1}^{\infty} \sum_{d|n}d^k \sin \left(   \frac{ 2\pi dh}{q} \right) n^{-s} \mathrm{d}s \nonumber \\ 
     & = \frac{1}{2} \left(\frac{2}{a \sqrt{x}}\right)^{ \nu}  \frac{1}{2\pi i} \int_{(c)}  \Gamma(s) \Gamma(s + \nu) \left( \frac{4}{a^2 x } \right)^s  \sum_{m=1}^{\infty}  m^{-s} \sum_{d=1}^{\infty}   {d^{k-s}}{\sin \left(   \frac{ 2\pi dh}{q} \right)}  \mathrm{d}s \nonumber \\ 
				& = \frac{q^k }{2}   X^{\nu/2} \frac{1}{2\pi i} \int_{(c)}  \Gamma(s) \Gamma(s + \nu)  \zeta(s) \sum_{r=1}^q \zeta\left(s-k,\frac{r}{q}\right) \sin\left(\frac{2\pi rh}{q}\right) (q^{-1}X)^s \mathrm{d}s,
			\end{align}
			where $c>k+1$ and $X = \frac{4}{a^2 x}$.  Next, we investigate the following integral
			\begin{align}\label{g2}
				&\mathcal{G}_{k}^{(\nu)}(X):=\frac{1}{2\pi i} \int_{(c)}  \Gamma(s + \nu)  \Gamma(s) \zeta(s) \sum_{r=1}^q \zeta\left(s-k,\frac{r}{q}\right) \sin\left(\frac{2\pi rh}{q}\right) (q^{-1}X)^s \mathrm{d}s \nonumber\\
   & = \frac{(-1)^{k/2} q^{-k}}{2(2\pi)^{k+1}}\frac{1}{2\pi i} \int_{(c)}  \Gamma(s + \nu) \zeta(1-s)  \Gamma(k+1-s)   \left\{ \zeta\left(k+1-s,\frac{h}{q}\right)-\zeta\left(k+1-s, 1-{\frac{h}{q}}\right)\right\}  ( 4\pi^2 X)^s \mathrm{d}s, 
\end{align}
   where in the last step, we used \eqref{1st_use}, \eqref{hur2}. Next, we consider the following integral
   \begin{align}\label{J1}
    \mathcal{H}_T:=   \frac{1}{2\pi i} \int_{c-iT}^{c+iT}  \Gamma(s + \nu) \zeta(1-s)  \Gamma(k+1-s)  \left\{ \zeta\left(k+1-s,\frac{h}{q}\right)-\zeta\left(k+1-s, 1-{\frac{h}{q}}\right)\right\}  ( 4\pi^2 X)^s \mathrm{d}s,
   \end{align}
   for some large positive number $T$. Let us consider the positively oriented rectangular contour $\mathcal{C}:$ consisting of the line segments $[c-iT, c+iT], [c+iT, -d+iT], [-d+iT, -d-iT]$ and $[-d-iT, c-iT]$ where the choice for $d$ is as follows: $ 0< d< \min \{1, \Re(\nu)\}$ whenever $\Re(\nu)>0$  and $0<d<1$ otherwise. One can note that the simple pole of $ \Gamma(k+1-s)$ at $s=k+1, k+3, k+5,\cdots $ will get cancelled by the trivial zeroes of $\zeta(1-s)$ except at $s=k+1$ when $ k=0$. By employing \eqref{SV} and  \eqref{bernoulli2}, one can easily derive that $\left( \zeta (s,\frac{h}{q} )-\zeta (s, 1-{\frac{h}{q}} )\right)$ encounters zeros at negative odd integers.  Hence the remaining  poles of $ \Gamma(k+1-s)$ at $s=k+2, k+4, k+6,\cdots $ will get cancelled by the simple zeroes of $\left( \zeta (k+1-s,\frac{h}{q} )-\zeta (k+1-s, 1-{\frac{h}{q}} )\right)$. Inside the contour $\mathcal{C}$, the integrand has a pole at $s=0$ and  possibly at $s=1$. By employing Cauchy's residue theorem, the integral in \eqref{J1} can be rewritten as
\begin{align}\label{I_T1}
 \mathcal{H}_T &=R_0+\delta_k~  R_{1} +\frac{1}{2\pi i} \left(  \int_{-d-i T}^{-d+iT}\pm \int_{-d\pm iT}^{c\pm iT}\right)   \Gamma(s + \nu) \zeta(1-s)  \Gamma(k+1-s)  
 \notag\\
    & \ \ \ \  \ \ \ \  \ \ \ \  \ \ \ \  \ \ \ \  \ \ \ \  \ \ \ \  \ \ \ \  \ \ \ \  \ \ \ \  \ \ \ \  \times \left\{ \zeta\left(k+1-s,\frac{h}{q}\right)-\zeta\left(k+1-s, 1-{\frac{h}{q}}\right)\right\}  ( 4\pi^2 X)^s \mathrm{d}s, 
\end{align}
 where  $\delta_k$ is defined in \eqref{del_k0}; and $R_0$ and $R_{1}$ are the residues at $s=0$ and $s=1$, respectively which are given as  
\begin{align}\label{g4}
   R_0= &- \Gamma(k+1) \Gamma( \nu)   \left\{ \zeta\left(k+1,\frac{h}{q}\right)-\zeta\left(k+1, 1-{\frac{h}{q}}\right)\right\} , \\
   R_1= & \frac{1}{2}\Gamma(\nu+1)   \left\{ \zeta\left(0,\frac{h}{q}\right)-\zeta\left(0, 1-{\frac{h}{q}}\right)\right\} (4\pi^2X).\label{g5}
\end{align}
From \cite[p.~82, eq 14]{MR1790423}, we recall that $|L(s,\chi)|\leq 2q|s|$ for $\sigma\geq 1/2$. Then by Lemma \ref{Gamma} together with functional equation \eqref{exact l}, we find that in a bounded vertical strip  $L(s, \chi) \ll_{q, \sigma} |t|^{O_{\sigma, q}(1)}$
 with $|t|>1$. Now by \eqref{Lbound} we can easily deduce that   $\zeta\left(s, \frac{h}{q}\right) \ll_{q, \sigma}|t|^{O_{\sigma, q}(1)}$ in a bounded vertical strip. Since $-1<-d<\sigma<c<k+2$, employing Lemma \ref{Gamma} and the aforementioned bound for Hurwitz zeta function, we can conclude that the integrals along the horizontal segments, i.e., the last two integrals in \eqref{I_T1}, vanish as $T \rightarrow \infty$. Hence, letting $T \rightarrow \infty$ and then
substituting back the expression \eqref{I_T1} in \eqref{g2}, we get
\begin{align}\label{g3}
				\mathcal{G}_{k}^{(\nu)}(X)  &= \frac{(-1)^{k/2} q^{-k}}{2(2\pi)^{k+1}} ~ \left( R_0+\delta_k~  R_{1}+     \mathcal{A}_{k}^{(\nu)}(X) \right), 
			\end{align}
and $\mathcal{A}_{k}^{(\nu)}(X)$ is defined by
\begin{align}\label{g6}
 \mathcal{A}_{k}^{(\nu)}(X):= & \int_{(-d)}  \Gamma(s + \nu) \zeta(1-s)  \Gamma(k+1-s)   \left\{ \zeta\left(k+1-s,\frac{h}{q}\right)-\zeta\left(k+1-s, 1-{\frac{h}{q}}\right)\right\}  ( 4\pi^2 X)^s \mathrm{d}s \nonumber\\
 =&  4\pi^2 X \frac{1}{2\pi i} \int_{(1+d)}  \Gamma(\nu+1-s) \zeta(s)  \Gamma(k+s)   \left\{ \zeta\left(k+s,\frac{h}{q}\right)-\zeta\left(k+s, 1-{\frac{h}{q}}\right)\right\}  ( 4\pi^2 X)^{-s} \mathrm{d}s \nonumber\\
=& 4\pi^2 X \sum_{r=1}^\infty \sum_{m=0}^\infty \left\{(m+h/q)^{-k} \frac{1}{2\pi i}\int_{(1+d)}  \Gamma(\nu+1-s)    \Gamma(k+s)   \left( 4\pi^2 Xr(m+\frac{h}{q})\right)^{-s}\mathrm{d}s  
\right.\notag\\&\left.\ \ 
- (m+1-h/q)^{-k} \frac{1}{2\pi i}\int_{(1+d)}  \Gamma(\nu+1-s)    \Gamma(k+s)   \left( 4\pi^2 Xr(m+1-\frac{h}{q})\right)^{-s} \mathrm{d}s     \right\}.  
    \end{align}
Employing Lemma  \ref{lem4} for $n =k$ and $a=1+\nu$ in \ref{g6},  we get
\begin{align}\label{g7}
 \mathcal{A}_{k}^{(\nu)}(X)=&(4\pi^2 X)^{k+1} \Gamma(\nu+k+1)\sum_{r=1}^\infty \sum_{m=0}^\infty \left\{ \frac{ r^k }{\left(1+4\pi^2Xr(m+h/q) \right)^{\nu+k+1}}
-  \frac{ r^k }{\left(1+4\pi^2Xr(m+1-h/q) \right)^{\nu+k+1}}  \right\}.  
 \end{align}
Inserting \eqref{g4}, \eqref{g5} and  \eqref{g7} in \eqref{g3}, and then  combining with \eqref{g1} and \eqref{g2}, we obtain the result.
\par
 Next, we demonstrate that Theorem \ref{odd1_based} is equivalent to Theorem \ref{M1}. 
    
  \textbf{Theorem \ref{M1} $\Rightarrow$ Theorem \ref{odd1_based}} 
 We will prove the theorem for $\theta=h/q$, where $q$ is prime and $0<h<q.$ Now we multiply the identity \eqref{M11} in Theorem \ref{M1} with $ \frac{1}{i\phi(q)} \chi(h)  \tau(\bar{\chi})$ and then take the sum on odd primitive character $\chi$ modulo $q$. Hence, the left-hand side of the identity in \eqref{M11} becomes
\begin{align}
					 \frac{1}{i\phi(q)} \sum_{\chi \ odd } \chi(h)  \tau(\bar{\chi})   \sum_{n=1}^\infty\sigma_{k,\chi}(n)n^{\frac{\nu}{2}}K_\nu(a\sqrt{nx}) &= \frac{1}{i\phi(q)}\sum_{n=1}^\infty n^{\frac{\nu}{2}}K_\nu(a\sqrt{nx}) \sum_{d|n}d^k \sum_{\chi \ odd } \chi(d)  \chi(h)  \tau(\bar{\chi}) \notag\\
					&=\sum_{n=1}^{\infty} n^{\nu/2} K_{\nu}(a\sqrt{nx}) \sum_{d|n}d^k \sin \left(\frac{2\pi d h}{q}\right),\label{1.1}
				\end{align}
				where in the last step, we have used \eqref{cos}. The sum in the first term of the right-hand side of \eqref{M11} becomes
     
				\begin{align}
					\frac{1}{i\phi(q)} \sum_{\chi \ odd } \chi(h)\tau(\chi) \tau(\bar{\chi}) L(1+k, \bar{\chi})  =- \frac{q^{-k}}{2i}\left(\zeta(1+k, h/q) - \zeta(1+k, 1-h/q)\right),\label{1.2}
				\end{align} 
   where we have used \eqref{Hurwitz}, \eqref{both}  and \eqref{prop}. Similarly, the sum in the second term of \eqref{M11} becomes\begin{align}
       \frac{1}{i\phi(q)} \sum_{\chi \ odd } \chi(h)  \tau(\bar{\chi}) L(1, \chi)  = -\frac{\pi }{q \phi(q)} \sum_{\chi \ odd } \chi(h)\tau(\chi) \tau(\bar{\chi}) L(0, \bar{\chi}) = \frac{\pi }{2}\left(\zeta(0, h/q) - \zeta(0, 1-h/q)\right). \label{1.22}
       \end{align}

  The infinite sum in the last term on the right-hand side of \eqref{M11} transforms into the following
				\begin{align}
					& \frac{1}{i\phi(q)} \sum_{\chi \ odd } \chi(h)\tau(\chi) \tau(\bar{\chi}) \sum_{n=1}^\infty \Bar{\sigma}_{k,\Bar{\chi}}(n)  \frac{\Gamma(\nu+k+1)}{\left(\frac{16\pi^2}{a^2q}\frac{n}{x}+1\right)^{\nu+k+1} } = - \frac{q}{i\phi(q)}\sum_{n=1}^\infty \frac{\Gamma(\nu+k+1)}{\left(\frac{16\pi^2}{a^2q}\frac{n}{x}+1\right)^{\nu+k+1} }\sum_{d|n}d^k \sum_{\chi \ odd } \chi(h) \Bar{ \chi}(n/d) \notag\\
					&=   - \frac{q}{i\phi(q)} \sum_{d=1}^\infty d^k  \sum_{r=1}^\infty  \frac{\Gamma(\nu+k+1)}{\left(\frac{16\pi^2}{a^2q}\frac{dr}{x}+1\right)^{\nu+k+1} }   \sum_{\chi \ odd } \chi(h) \Bar{ \chi}(r) \notag\\
					&= - \frac{q}{2i }  \sum_{d=1}^\infty d^k  \left( \sum_{\substack{r=1\\r\equiv h( q)}}^\infty \frac{\Gamma(\nu+k+1)}{\left(\frac{16\pi^2}{a^2q}\frac{dr}{x}+1\right)^{\nu+k+1} }
					-\sum_{\substack{r=1\\r\equiv -h( q)}}^\infty \frac{\Gamma(\nu+k+1)}{\left(\frac{16\pi^2}{a^2q}\frac{dr}{x}+1\right)^{\nu+k+1} }  \right)      \notag\\
					&= - \frac{q}{2i }   \sum_{d=1}^{\infty}d^k \sum_{ m=0}^{\infty}\left(  \frac{ \Gamma(\nu+k+1)  }{(1+\frac{16\pi^2(mq+h)d}{qa^2x})^{1+\nu+k} } - \frac{ \Gamma(\nu+k+1)  }{(1+\frac{16\pi^2(mq+q-h)d}{qa^2x})^{1+\nu+k} } \right)     \notag\\ 
					&= - \frac{q}{2i }   \sum_{d=1}^{\infty}d^k \sum_{ m=0}^{\infty}\left(  \frac{ \Gamma(\nu+k+1)  }{(1+\frac{16\pi^2d}{a^2x}(m+h/q))^{1+\nu+k} } - \frac{ \Gamma(\nu+k+1)  }{(1+\frac{16\pi^2d}{a^2x}(m+1-h/q))^{1+\nu+k} } \right).           \label{1.3}
				\end{align} 
				 Employing \eqref{1.1}, \eqref{1.2}, \eqref{1.22}, and \eqref{1.3} in  \eqref{M11}, we get the desired result.\\
  \textbf{ Theorem \ref{odd1_based} $\Rightarrow$  Theorem \ref{M1} }  Let $\theta=h/q$, and  $\chi$ be an odd primitive character modulo $q$. We first multiply the identity \eqref{p7} in Theorem \ref{odd1_based}
				 by $\bar{\chi}(h)/\tau(\bar{\chi})$, and then take summation on $h$, $ 0<h<q$. The left-hand side of the identity \eqref{p7} becomes
				\begin{align}\label{p1}
				    & \frac{1}{\tau(\bar{\chi})}\sum_{h=1}^{q-1}\bar{\chi}(h)\sum_{n=1}^{\infty} n^{\nu/2} K_{\nu}(a\sqrt{nx}) \sum_{d|n}d^k \sin \left( 2\pi d   h/q \right)\notag\\
       &= \frac{1}{2 i \tau(\bar{\chi})} \sum_{n=1}^{\infty} n^{\nu/2} K_{\nu}(a\sqrt{nx}) \sum_{d|n}d^k  \sum_{h=1}^{q-1}\bar{\chi}(h) \left( e^{2\pi i  d   h/q} - e^{-2\pi i  d   h/q}\right)\notag\\
     &= \frac{1}{2 i  } \sum_{n=1}^{\infty} n^{\nu/2} K_{\nu}(a\sqrt{nx}) \sum_{d|n}d^k    \left(  \chi(d)-\chi(-d)\right)\notag\\
     &=i^{-1}\sum_{n=1}^\infty\sigma_{k,\chi}(n)n^{\frac{\nu} {2}}K_\nu(a\sqrt{nx}),
				\end{align} 
where in the penultimate step, we have used \eqref{gauss}. With the help of \eqref{Hurwitz} and \eqref{both}, the right-hand side of \eqref{p7} transforms into  
 \begin{align}\label{p2}
 \frac{1}{\tau(\bar{\chi})}\sum_{h=1}^{q-1}\bar{\chi}(h)\left(\zeta(1+k,h/q) - \zeta(1+k, 1- h/q)\right) =-2 q^{k}\tau({\chi})L(1+k,\Bar{\chi}),
\end{align} and employing the functional equation for $L$-function \eqref{ll(s)}, the second term becomes
\begin{align}\label{p3}
   \frac{1}{\tau(\bar{\chi})}\sum_{h=1}^{q-1}\bar{\chi}(h)\left(\zeta(0,h/q) - \zeta(0, 1- h/q)\right) &=-\frac{2\tau({\chi})}{q} L(0,\Bar{\chi})
   =i^{-1} L(1,{\chi}).
\end{align}
The infinite series on the right-hand side of \eqref{p7} takes the form
\begin{align}\label{p14}
&\frac{1}{\tau(\bar{\chi})}\sum_{h=1}^{q-1}\bar{\chi}(h)\sum_{d=1}^{\infty}d^k \sum_{ m=0}^{\infty}\left(  \frac{ 1 }{(1+\frac{16\pi^2d}{a^2x}(m+h/q))^{1+\nu+k} } - \frac{ 1  }{(1+\frac{16\pi^2d}{a^2x}(m+1-h/q))^{1+\nu+k} } \right)  \notag\\
&=\frac{1}{\tau(\bar{\chi})}\sum_{h=1}^{q-1}\bar{\chi}(h)\sum_{d=1}^{\infty}d^k \sum_{\substack{r=1 \\ r \equiv h(q)}}^\infty  \frac{ 1 }{(1+\frac{16\pi^2dr}{a^2xq} )^{1+\nu+k} } 
-\frac{1}{\tau(\bar{\chi})}\sum_{h=1}^{q-1}\bar{\chi}(h)\sum_{d=1}^{\infty}d^k \sum_{\substack{r=1 \\ r \equiv -h(q)}}^\infty  \frac{ 1 }{(1+\frac{16\pi^2dr}{a^2xq} )^{1+\nu+k} } \notag\\
&=\frac{2}{\tau(\bar{\chi})}  \sum_{d=1}^{\infty}  \sum_{ r=1  }^\infty  \frac{ d^k \bar{\chi}(r) }{(1+\frac{16\pi^2dr}{a^2xq} )^{1+\nu+k} }   =-\frac{2\tau({\chi}) }{q}   \sum_{ n=1  }^\infty  \frac{  \Bar{\sigma}_{k,\Bar{\chi} }(n)}{(1+\frac{16\pi^2dr}{a^2xq} )^{1+\nu+k} }.   
\end{align}
Inserting \eqref{p1}, \eqref{p2}, \eqref{p3}  and \eqref{p14} into  \eqref{p7} we get the result.
    \end{proof}
    \begin{proof}[Theorem \rm{\ref{odd2_based}} and its equivalence with Theorem \rm{\ref{M2}} ][]
    The proof of Theorem \ref{odd2_based} is similar to the proof of Theorem \ref{odd1_based}. The equivalence of Theorem \ref{odd2_based} and Theorem \ref{M2} can be derived in a similar way as the previous one. To avoid repetitions, we skip the details of the proof.
    \end{proof}
\begin{proof}[Corollary \rm{\ref{cor1r}}][] 
Multiplying \eqref{p7} of Theorem  \ref{odd1_based} by $-4$ and  \eqref{r7} of Theorem \ref{odd2_based} by $16$, and then adding both the expressions, putting $k=2$ yield 
\begin{align}\label{m1}
    &\sum_{n=1}^{\infty}  n^{\nu/2} K_{\nu}(a\sqrt{nx}) \sum_{d|n}d^2  \left\{ 16\sin \left(\frac{2\pi n \theta}{d}\right) -4\sin \left( 2\pi d \theta \right)\right\} \notag\\
        &=- {16}  \pi^3\Gamma(\nu+3) X^{\frac{\nu}{2}+3 }( \zeta(-2,\theta) - \zeta(-2,1-\theta) ) -\frac{1}{4\pi^3}X^{\frac{\nu}{2} } \Gamma(\nu) ( \zeta(3,\theta) - \zeta(3,1-\theta) )\notag\\
        &+X^{\frac{\nu}{2} }(2 \pi X)^3\Gamma(\nu+3)\sum_{n=1}^\infty\sum_{m=0}^\infty   \left\{   \frac{  n^2-4(m+\theta)^2  }{(1+\frac{16\pi^2r}{a^2x}(m+\theta))^{\nu+3} }    -  \frac{  n^2-4(m+1-\theta)^2  }{(1+\frac{16\pi^2r}{a^2x}(m+1-\theta))^{\nu+3} }   \right\}.
        	\end{align}

  Using \eqref{SV}, we have 
         \begin{align}\label{m2}
             \zeta(-2,\theta)- \zeta(-2,1-\theta)=-\frac{1}{3}({B_{3}(\theta)}-B_{3}(1-\theta))=-\frac{1}{3}(\theta-3\theta^2+2\theta^3   ). 
         \end{align}
        By using partial fraction expansion for $\cot(\pi \theta)$, we can easily get
         \begin{align*}
             \zeta(n,1-\theta)+(-1)^n \zeta(n,\theta)=-\frac{\pi }{(n-1)!}\frac{d^{n-1}}{dx^{n-1}}\cot(\pi\theta).
         \end{align*}
         Hence by above, we get
          \begin{align}\label{m3}
             ( \zeta(3,\theta) - \zeta(3,1-\theta) )=\pi^3(\cot(\pi \theta)+\cot^3(\pi \theta)).
         \end{align}
         Substituting \eqref{m2}, \eqref{m3} in \eqref{m1}, we get the result.
      \end{proof}

    \begin{proof}[Corollary \rm{\ref{cor2r}}][] 
Proof directly follows from \eqref{berndt} and Corollary \rm{\ref{cor1r}}.
\end{proof}
  \begin{proof}[Corollary \rm{\ref{cor3r}}][] 
Substituting $\nu=1/2,a=4\pi $ in Corollary \rm{\ref{cor2r}} and using \eqref{property}, we get the result.
 \end{proof}
 
    \begin{proof}[Theorem \rm{\ref{even1_based}} and its equivalence with Theorem \ref{thmeven1}][]
       The proof is similar to the proof of Theorem \rm{\ref{odd1_based}}. Here, we use \eqref{hur1} instead of \eqref{hur2} and $k\geq1$ to be an odd integer. 
       \par
       Next, we demonstrate that  Theorem \ref{thmeven1} and Proposition \ref{first paper}  are equivalent to Theorem \ref{even1_based}. \\
     \textbf{ Theorem \ref{thmeven1} $\Rightarrow$ Theorem \ref{even1_based} } It is sufficient to prove the theorem for $\theta=h/q$, where $q$ is prime and $0<h<q.$ 
	We begin our proof by considering the expression on the left-hand side of the identity in Theorem \ref{even1_based}. Employing \eqref{sin}, we have
				\begin{align}
					&\sum_{n=1}^{\infty} n^{\nu/2} K_{\nu}(a\sqrt{nx}) \sum_{d|n}d^k \cos \left(\frac{ 2\pi dh }{q}\right) =\sum_{n=1}^{\infty} n^{\nu/2} K_{\nu}(a\sqrt{nx}) \left(\sum_{\substack{d|n\\ q|d}} d^k  +\sum_{\substack{d|n\\ q\nmid d}}d^k \cos \left(\frac{ 2\pi dh }{q}\right)\right)\notag\\
					&=\sum_{m=1}^{\infty} (qm)^{\nu/2} K_{\nu}(a\sqrt{qmx}) \sum_{d|m}(qd)^k+ \sum_{n=1}^{\infty} n^{\nu/2} K_{\nu}(a\sqrt{nx}) \sum_{\substack{d|n\\ q\nmid d}}           \frac{d^k}{\phi(q)} \sum_{\chi \ even } \chi(d)  \chi(h)  \tau(\bar{\chi}) \notag\\
			  &=q^{\frac{\nu}{2}+k}\sum_{m=1}^{\infty} m^{\nu/2} K_{\nu}(a\sqrt{qmx}) \sum_{d|m}d^k-  \sum_{n=1}^{\infty} n^{\nu/2} K_{\nu}(a\sqrt{nx}) \sum_{\substack{d|n\\ q\nmid d}}           \frac{d^k}{\phi(q)}  \chi_0(d)           \notag\\&\ \ \ \ \ \ \ \ \ \ \ \ \ \ \ \ \ \ \ \ \ \ \ \ \ \ \ \ \ \ \ \ \ \ \ \ \ \ \ \ \ \ \ \ \ \ \ \ \ \ \ \  +\sum_{n=1}^{\infty} n^{\nu/2} K_{\nu}(a\sqrt{nx}) \sum_{\substack{d|n\\ q\nmid d}}           \frac{d^k}{\phi(q)} \sum_{\substack{\chi \neq \chi_0\\\chi even}}\chi(d)  \chi(h)  \tau(\bar{\chi}) \notag\\
				 &=q^{\frac{\nu}{2}+k}\sum_{m=1}^{\infty} m^{\nu/2} K_{\nu}(a\sqrt{qmx}) \sum_{d|m}d^k-  \sum_{n=1}^{\infty} n^{\nu/2} K_{\nu}(a\sqrt{nx})\frac{1}{ {\phi(q)}}\left(  \sum_{d|n}d^k- \sum_{\substack{d|n\\ q| d}}            d^k    \right)        \notag\\& \ \ \ \ \ \ \ \ \  \ \ \ \ \ \ \ \ \ \ \ \ \ \ \ \ \ \ \ \ \ \ \ \ \ \ \ \ \ \ \ \ \ \ \ \ \ \ \ \ \ \ \ \ \ \ \  +\frac{1}{\phi(q)} \sum_{\substack{\chi \neq \chi_0\\\chi even}}    \chi(h)  \tau(\bar{\chi}) \sum_{n=1}^{\infty} n^{\nu/2} K_{\nu}(a\sqrt{nx}) \sum_{ d|n }       d^k\chi(d) \notag\\ 
					&=\frac{q^{\frac{\nu}{2}+k+1}}{\phi(q)}\sum_{m=1}^{\infty} \sigma_k(m)m^{\nu/2} K_{\nu}(a\sqrt{qmx})  -  \frac{1}{ {\phi(q)}}\sum_{n=1}^{\infty}\sigma_k(n) n^{\nu/2} K_{\nu}(a\sqrt{nx})         \notag\\   &\ \ \ \ \ \ \ \ \ \  \ \ \ \ \ \ \ \ \ \ \ \ \ \ \ \ \ \ \ \ \ \ \ \ \ \ \ \ \ \ \ \ \ \ \ \ \ \ \ \ \ \ \ \ \ \ \   + \frac{1}{\phi(q)} \sum_{\substack{\chi \neq \chi_0\\\chi even}}    \chi(h)  \tau(\bar{\chi}) \sum_{n=1}^{\infty} \sigma_{k,\chi}(n)n^{\nu/2} K_{\nu}(a\sqrt{nx} ).  
					\label{even1.1}
				\end{align}
  
				Now, we first evaluate the first two sums on the right-hand side of \eqref{even1.1}. By Proposition \ref{first paper}, we have
				\begin{align}
					&\frac{q^{\frac{\nu}{2}+k+1}}{\phi(q)}\sum_{m=1}^{\infty} \sigma_k(m)m^{\nu/2} K_{\nu}(a\sqrt{qmx})  -  \frac{1}{ {\phi(q)}}\sum_{n=1}^{\infty}\sigma_k(n) n^{\nu/2} K_{\nu}(a\sqrt{nx})=-\frac{\Gamma(\nu) \zeta(-k)}{ 4\phi(q)} X^{\frac{\nu}{2}} \left( q^{k+1}-1    \right)     \notag\\ 
					&-\frac{\Gamma(1+\nu) }{4}X^{1+\frac{\nu}{2}} \delta_{k,1} +\frac{(-1)^{\frac{k+1}{2}}}{2\phi(q)} X^{\frac{\nu}{2}}\Gamma(\nu+k+1)\left(  2\pi X\right)^{k+1}\left(\sum_{n=1}^\infty \frac{ \sigma_{k }(n)}{\left(\frac{16\pi^2}{a^2q}\frac{n}{x}+1\right)^{\nu+k+1} }-\sum_{n=1}^\infty \frac{ \sigma_{k }(n)}{\left(\frac{16\pi^2}{a^2}\frac{n}{x}+1\right)^{\nu+k+1} }\right), \label{even1.2}
				\end{align}
		where $\delta_{k,1}$ is defined in \eqref{del_k}.		Now, we  examine the last sum on the right-hand side of \eqref{even1.1}. By  \eqref{thm3} of Theorem \eqref{thmeven1}, we have
				\begin{align}
					&\frac{1}{\phi(q)} \sum_{\substack{\chi \neq \chi_0\\\chi even}}    \chi(h)  \tau(\bar{\chi}) \sum_{n=1}^{\infty} \sigma_{k,\chi}(n)n^{\nu/2} K_{\nu}(a\sqrt{nx})  
					=\frac{(-1)^{\frac{k-1}{2}}k!q^{k+1}}{2(2\pi)^{k+1}\phi(q)} X^\frac{\nu}{2} \Gamma(\nu) \sum_{\substack{\chi \neq \chi_0\\\chi even}}    \chi(h)L(1+k,\Bar{\chi})\notag\\  &\ \ \ \ \ \ \ \ \  \ \ \ \ \ \ \ \ \ \ \ \ \ \ \ \ \ \ \ \ \ \ \ \ \ \ \ \ \ + \frac{(-1)^{\frac{k+1}{2}}}{2\phi(q)} X^\frac{\nu}{2}(2 \pi X)^{k+1}\Gamma(\nu+k+1)\sum_{\substack{\chi \neq \chi_0\\\chi even}}    \chi(h)\sum_{n=1}^\infty   \frac{ \Bar{\sigma}_{k,\Bar{\chi}}(n)}{\left(\frac{16\pi^2}{a^2q}\frac{n}{x}+1\right)^{\nu+k+1} }.\label{even1.3}
				\end{align}
				We consider
				\begin{align}
					\sum_{\substack{\chi \neq \chi_0\\ \chi\  even}}    \chi(h)L(1+k,\Bar{\chi})&= \sum_{\chi \ even}\chi(h)L(1+k,\Bar{\chi})-L(1+k, {\chi_0})\notag\\
					&=  \sum_{\chi\  even}\chi(h) \frac{1}{q^{k+1}}\sum_{r=1}^{q-1}\bar{\chi}(r) \zeta(k+1,r/q)- \left(\sum_{n=1}^\infty\frac{1}{n^{k+1}} -\sum_{n=1}^\infty\frac{1}{(nq)^{k+1}}\right)\notag\\
					&=   \frac{1}{q^{k+1}}\sum_{r=1}^{q-1}  \zeta(k+1,r/q)   \sum_{\chi\  even}\chi(h)\bar{\chi}(r)- \left( 1-\frac{1}{q^{k+1}}  \right)\zeta(k+1) \notag\\
					&=   \frac{\phi(q)}{2q^{k+1}} \left\{ \zeta(k+1,h/q)+\zeta(k+1,1-h/q) \right\} - \left( 1-\frac{1}{q^{k+1}}  \right)\zeta(k+1),   \label{even1.4}
				\end{align}where in the penultimate step, we used \eqref{Hurwitz} and \eqref{prop1}. Now, we examine the last expression in \eqref{even1.3}, and we obtain
				\begin{align}
					& \frac{1}{\phi(q)}\sum_{\substack{\chi \neq \chi_0\\\chi even}}    \chi(h)\sum_{n=1}^\infty   \frac{\Bar{\sigma}_{k,\Bar{\chi}}(n) }{\left(\frac{16\pi^2}{a^2q}\frac{n}{x}+1\right)^{\nu+k+1} }  
     = \frac{1}{\phi(q)} \sum_{n=1}^\infty   \frac{\sum_{d/n}d^k}{\left(\frac{16\pi^2}{a^2q}\frac{n}{x}+1\right)^{\nu+k+1} }    \sum_{\substack{\chi \neq \chi_0\\\chi even}}    \chi(h) \bar{\chi}(n/d)  \notag\\
					&= \frac{1}{\phi(q)} \sum_{n=1}^\infty   \frac{\sum_{d/n}d^k}{\left(\frac{16\pi^2}{a^2q}\frac{n}{x}+1\right)^{\nu+k+1} }    \left\{ \sum_{ \chi even} \chi(h) \bar{\chi}(n/d) -  \chi_0(n/d)\right\}   \notag\\
			 		&= \frac{1}{2} \sum_{d=1}^\infty d^k \sum_{\substack{r=1\\r\equiv \pm h(mod q) }}^\infty  \frac{1}{\left(\frac{16\pi^2}{a^2q}\frac{dr}{x}+1\right)^{\nu+k+1} } - \frac{1}{\phi(q)} \sum_{n=1}^\infty   \frac{ \left(\sigma_k(n)-\sigma_k(n/q)\right)}{\left(\frac{16\pi^2}{a^2q}\frac{n}{x}+1\right)^{\nu+k+1} }   \notag\\
					&=\frac{1}{2} \sum_{d=1}^\infty d^k \sum_{m=0 }^\infty \left\{   \frac{1}{\left(\frac{16\pi^2}{a^2q}\frac{d(mq+h)}{x}+1\right)^{\nu+k+1} }   +\frac{1}{\left(\frac{16\pi^2}{a^2q}\frac{d(mq+q-h)}{x}+1\right)^{\nu+k+1} }        \right\}  \notag\\ &\ \ \ \ \ \ \ \ \ \ \ \  \ \ \ \ \ \ \ \ \ \ \ \ \ \ \ \ \ \ \ \ \ \ \ \ \ \ \ \ \ - \frac{1}{\phi(q)} \sum_{n=1}^\infty   \frac{  \sigma_k(n) }{\left(\frac{16\pi^2}{a^2q}\frac{n}{x}+1\right)^{\nu+k+1} }  
					+\frac{1}{\phi(q)} \sum_{n=1}^\infty   \frac{  \sigma_k(n/q) }{\left(\frac{16\pi^2}{a^2q}\frac{n}{x}+1\right)^{\nu+k+1} }\notag\\
			&=\frac{1}{2} \sum_{d=1}^\infty d^k \sum_{m=0 }^\infty \left\{   \frac{1}{\left(\frac{16\pi^2}{a^2}\frac{d(m+h/q)}{x}+1\right)^{\nu+k+1} }   +\frac{1}{\left(\frac{16\pi^2}{a^2}\frac{d(m+1-h/q)}{x}+1\right)^{\nu+k+1} }        \right\}  \notag\\ &\ \ \ \ \ \ \ \ \ \  \ \ \ \ \ \ \ \ \ \ \ \ \ \ \ \ \ \ \ \ \ \ \ \ \ \ \ \ \  - \frac{1}{\phi(q)} \sum_{n=1}^\infty   \frac{  \sigma_k(n) }{\left(\frac{16\pi^2}{a^2q}\frac{n}{x}+1\right)^{\nu+k+1} }  
					+\frac{1}{\phi(q)} \sum_{r=1}^\infty   \frac{  \sigma_k(r) }{\left(\frac{16\pi^2}{a^2}\frac{r}{x}+1\right)^{\nu+k+1} }.
					\label{even1.5}
				\end{align}
				
				Substituting \eqref{even1.4}, \eqref{even1.5} in \eqref{even1.3}, we obtain
				\begin{align}
					&\frac{1}{\phi(q)} \sum_{\substack{\chi \neq \chi_0\\\chi even}}    \chi(h)  \tau(\bar{\chi}) \sum_{n=1}^{\infty} \sigma_{k,\chi}(n)n^{\nu/2} K_{\nu}(a\sqrt{nx}) \notag\\
					&=\frac{(-1)^{\frac{k-1}{2}}k! }{2(2\pi)^{k+1} } X^\frac{\nu}{2} \Gamma(\nu)\left\{  
					\frac{ 1}{2 } \left( \zeta(k+1,h/q)+\zeta(k+1,1-h/q) \right) -  \frac{ q^{k+1}- 1 }{\phi(q)}  \zeta(k+1) \right\}\notag\\
					& +\frac{(-1)^{\frac{k+1}{2}}}{4} X^\frac{\nu}{2}(2 \pi X)^{k+1}\Gamma(\nu+k+1)\sum_{d=1}^\infty d^k \sum_{m=0 }^\infty \left\{   \frac{1}{\left(\frac{16\pi^2}{a^2}\frac{d(m+h/q)}{x}+1\right)^{\nu+k+1} }  
					 +\frac{1}{\left(\frac{16\pi^2}{a^2}\frac{d(m+1-h/q)}{x}+1\right)^{\nu+k+1} }      \right\}  \notag\\
					&\  \ \ \ \  \ \ \ \ \ \ \ \ \ \ \ \ \ \ \ \ \ \ \ \ \ \ \ \ \ \ - \frac{(-1)^{\frac{k+1}{2}}}{2\phi(q)} X^\frac{\nu}{2}(2 \pi X)^{k+1}\Gamma(\nu+k+1) \sum_{n=1}^\infty   \frac{  \sigma_k(n) }{\left(\frac{16\pi^2}{a^2q}\frac{n}{x}+1\right)^{\nu+k+1} } \notag\\
					&   \ \ \ \ \ \  \ \ \ \ \ \ \ \ \ \ \ \ \ \ \ \ \ \ \ \ \ \ \ \ \ \   \ \ \ \ \ \ \ \ \  \ \ \ \ \ \ \ \ \ \ 
     + \frac{(-1)^{\frac{k+1}{2}}}{2\phi(q)} X^\frac{\nu}{2}(2 \pi X)^{k+1}\Gamma(\nu+k+1)\sum_{r=1}^\infty   \frac{  \sigma_k(r) }{\left(\frac{16\pi^2}{a^2}\frac{r}{x}+1\right)^{\nu+k+1} }.
					\label{even1.6}\end{align}
				Inserting   \eqref{even1.2} and \eqref{even1.6} into \eqref{even1.1}, we get the result.

    \textbf{ Theorem \ref{even1_based} $\Rightarrow$ Theorem \ref{thmeven1} }   Let $\theta=h/q$, and  $\chi$ be an even primitive non-principal character modulo $q$. Multiplying the identity \eqref{1290} in Theorem \ref{even1_based}
				 by $\bar{\chi}(h)/\tau(\bar{\chi})$, and then summing on $h$, $ 0<h<q$.   The remaining steps are similar to the proof of Theorem \ref{M1}.
    \end{proof}
		 \begin{proof}[Theorem \rm{\ref{even2_based}} and its equivalence with Theorem \rm{\ref{even2}}][] 
The proof is similar to the previous one. To avoid repetition, we skip the detail.
\end{proof}

    \begin{proof}[Theorem \rm{\ref{botheven_odd1_based}} and its equivalence with Theorem \rm{\ref{botheven_odd}}][]
      It is sufficient to prove the theorem for rationals $\theta=h_1/p$ and $\psi=h_2/q$ where $p$ and $q$ are primes, $0<h_1<p$ and $0<h_2<q$. 
Employing   Lemma \ref{eq:bessel} with $t=n$ and subsequently interchanging the summation and integration, we get for odd integer $k\geq 1$, 
\begin{align}\label{g11}
				&\sum_{n=1}^{\infty}  \sum_{d|n}d^k \sin \left(   \frac{ 2\pi dh_1}{p} \right) \sin \left(   \frac{ 2\pi nh_2}{dq} \right) n^{\nu/2} K_{\nu}(a\sqrt{nx}) \nonumber \\ 
     &= \frac{1}{2} \left(\frac{2}{a \sqrt{x}}\right)^{ \nu}  \frac{1}{2\pi i} \int_{(c)}  \Gamma(s) \Gamma(s + \nu) \left( \frac{4}{a^2 x } \right)^s  \sum_{n=1}^{\infty} \sum_{d|n}d^k \sin \left(   \frac{ 2\pi dh_1}{p} \right) \sin \left(   \frac{ 2\pi nh_2}{dq} \right) n^{-s} \mathrm{d}s \nonumber \\ 
      =& \frac{1}{2} \left(\frac{2}{a \sqrt{x}}\right)^{ \nu}  \frac{1}{2\pi i} \int_{(c)}  \Gamma(s) \Gamma(s + \nu) \left( \frac{4}{a^2 x } \right)^s  \sum_{d=1}^{\infty} {d^{k-s}} {\sin \left(   \frac{ 2\pi dh_1}{p} \right)} \sum_{m=1}^{\infty}  m^{-s}\sin \left(   \frac{ 2\pi mh_2}{q} \right)   \mathrm{d}s \nonumber \\ 
				=& \frac{p^k }{2}   X^{\nu/2} \frac{1}{2\pi i} \int_{(c)}  \Gamma(s) \Gamma(s + \nu)    \sum_{r_1=1}^p \zeta\left(s-k,\frac{r_1}{p}\right) \sin\left(\frac{2\pi r_1h_1}{p}\right)\sum_{r_2=1}^q \zeta\left(s,\frac{r_2}{q}\right) \sin\left(\frac{2\pi r_2h_2}{q}\right) (p^{-1}q^{-1}X)^s \mathrm{d}s,
			\end{align}
			where $c>k+1$ and $X = \frac{4}{a^2 x}$.  Next, we  invoke the functional relations     \eqref{hur2} and \eqref{reflection formula for gamma add} in \eqref{g11}, we obtain
 \begin{align}\label{g12}
				&\sum_{n=1}^{\infty}  \sum_{d|n}d^k \sin \left(   \frac{ 2\pi dh_1}{p} \right) \sin \left(   \frac{ 2\pi nh_2}{dq} \right) n^{\nu/2} K_{\nu}(a\sqrt{nx}) \nonumber \\ 
    &= \frac{ (-1)^{ \frac{k-1}{2} } }{8 (2\pi )^{k+1}}   X^{\nu/2} \frac{1}{2\pi i} \int_{(c)}   \Gamma(s + \nu) \Gamma(k+1-s) \left\{ \zeta\left(k+1-s,\frac{h_1}{p}\right)-\zeta\left(k+1-s, 1-{\frac{h_1}{p}}\right)\right\}
  \nonumber \\ 
 & \ \ \ \   \ \ \ \ \ \ \ \   \ \ \ \  \ \ \ \   \ \ \ \ \ \ \ \   \ \ \ \  \ \ \ \   \ \ \ \   \times \left\{ \zeta\left(1-s,\frac{h_2}{q}\right)-\zeta\left(1-s, 1-{\frac{h_2}{q}}\right)\right\}   ( 4\pi^2X)^s \mathrm{d}s.
 \end{align}
Now proceeding in the same way as in the proof of Theorem \ref{odd1_based}, we shift the line of integration to $\Re(s)=-d$ (with $0<d< \min \{1, \Re(\nu)\}$), then replace $1-s$ by $s$ to obtain
\begin{align}\label{g13}
				&\sum_{n=1}^{\infty}  \sum_{d|n}d^k \sin \left(   \frac{ 2\pi dh_1}{p} \right) \sin \left(   \frac{ 2\pi nh_2}{dq} \right) n^{\nu/2} K_{\nu}(a\sqrt{nx}) \nonumber \\ 
    &= \frac{ (-1)^{ \frac{k-1}{2} } \pi^2}{2 (2\pi )^{k+1}}      X^{1+\nu/2} \frac{1}{2\pi i} \int_{(1+d)}   \Gamma(\nu+1-s) \Gamma(k+s) \left\{ \zeta\left(k+s,\frac{h_1}{p}\right)-\zeta\left(k+s, 1-{\frac{h_1}{p}}\right)\right\}
  \nonumber \\ 
 & \ \ \ \   \ \ \ \ \ \ \ \   \ \ \ \  \ \ \ \   \ \ \ \ \ \ \ \   \ \ \ \  \ \ \ \   \ \ \ \   \times \left\{ \zeta\left(s,\frac{h_2}{q}\right)-\zeta\left(s, 1-{\frac{h_2}{q}}\right)\right\}   ( 4\pi^2X)^{-s} \mathrm{d}s \nonumber \\ 
 &= \frac{ (-1)^{ \frac{k-1}{2} } \pi^2}{2 (2\pi )^{k+1}}      X^{1+\nu/2}\sum_{m=0}^\infty \sum_{n=0}^\infty \left\{
 \right.\notag\\&\left.\ \ 
 (n+h_1/p)^{-k} \frac{1}{2\pi i}\int_{(1+d)}  \Gamma(\nu+1-s)    \Gamma(k+s)   \left( 4\pi^2 X(n+h_1/p)(m+h_2/q)\right)^{-s}\mathrm{d}s
 \right.\notag\\&\left.\ \ 
- (n+1-h_1/p)^{-k} \frac{1}{2\pi i}\int_{(1+d)}  \Gamma(\nu+1-s)    \Gamma(k+s)   \left( 4\pi^2 X(n+1-h_1/p)(m+ h_2/q)\right)^{-s} \mathrm{d}s   
 \right.\notag\\&\left.\ \ 
- (n+h_1/p)^{-k} \frac{1}{2\pi i}\int_{(1+d)}  \Gamma(\nu+1-s)    \Gamma(k+s)   \left( 4\pi^2 X(n+h_1/p)(m+1- h_2/q)\right)^{-s} \mathrm{d}s   
\right.\notag\\&\left.\ \ 
+ (n+1-h_1/p)^{-k} \frac{1}{2\pi i}\int_{(1+d)}  \Gamma(\nu+1-s)    \Gamma(k+s)   \left( 4\pi^2 X(n+1-h_1/p)(m+1- h_2/q)\right)^{-s} \mathrm{d}s     \right\}.
\nonumber \\ 
 \end{align}
We get the desired result by applying Lemma \ref{lem4} for $n =k$ and $a=1+\nu$ to each of these integrals in \eqref{g13}.

The equivalence of this theorem with  Theorem \ref{botheven_odd} can  be easily deduced. So, we leave the details for the reader.
 \end{proof}

    \begin{proof}[Theorem \rm{\ref{botheven_odd2_based}} and its equivalence with Theorem \rm{\ref{botheven_odd}}][] The proof is similar to Theorem \ref{botheven_odd1_based}, so we skip the proof. But we would show its equivalence with Theorem \ref{botheven_odd}. 
    
    \textbf{ Theorem \ref{botheven_odd} $\Rightarrow$ Theorem {\ref{botheven_odd2_based}} }
   It is sufficient to prove the theorem for rationals $\theta=h_1/p$ and $\psi=h_2/q$ where $p$ and $q$ are primes, $0<h_1<p$ and $0<h_2<q$. 
 We multiply both sides of identity \eqref{THM5} in Theorem \ref{botheven_odd} by $ \chi_1(h_1)  \tau(\bar{\chi_1}) /\phi(p)$  and   $ \chi_2(h_2)  \tau(\bar{\chi_2}) /\phi(q)$, then  sum on non-principal   primitive even character $\chi_1$ modulo $p$ and  $\chi_2$ modulo $q$. Using \eqref{sin}, the left-hand side of \eqref{THM5} becomes  \vspace{-3mm}\begin{align}\label{big1}
&\frac{1}{\phi(p)\phi(q)}\sum_{\substack{\chi_2\neq \chi_0\\\chi_2 \ even }}\chi_{2}(h_2)\tau(\Bar{\chi_{2}})\sum_{\substack{\chi_1\neq \chi_0\\\chi_1 \ even }}\chi_{1}(h_1)\tau(\Bar{\chi_{1}})   \sum_{n=1}^\infty \sigma_{k, {\chi_1}, {\chi_2}}(n) n^{\nu/2} K_{\nu}(a\sqrt{nx})        \notag\\
				 &=\frac{1}{\phi(p)\phi(q)}\sum_{n=1}^\infty n^{\nu/2} K_{\nu}(a\sqrt{nx})\sum_{d|n}d^k  \left\{ \sum_{\substack{\chi_2 \neq \chi_0\\\chi_2\ even}}\chi_{2}(h_2) { \chi_2}(n/d)\tau(\Bar{\chi_{2}}) \right\} \left \{ \sum_{\substack{\chi_1 \neq \chi_0\\\chi_1\ even}}\chi_{1}(h_1)  { \chi_1}(d)\tau(\Bar{\chi_{1}})  \right \}           \notag\\
 &= \sum_{n=1}^\infty n^{\nu/2} K_{\nu}(a\sqrt{nx})\sum_{d|n}d^k  \left\{ \frac{1}{\phi(q) }\sum_{ \chi_2\ even}\chi_{2}(h_2) { \chi_2}(n/d)\tau(\Bar{\chi_{2}})+  \frac{\chi_0(n/d)}{\phi(q)}  \right\}  \notag\\
  &\ \ \ \ \ \ \ \ \ \ \ \ \ \ \ \ \ \ \ \ \ \ \ \ \ \ \ \  \times\left \{\frac{1}{\phi(p) } \sum_{ \chi_1\ even}\chi_{1}(h_1)  { \chi_1}(d)\tau(\Bar{\chi_{1}}) + \frac{\chi_0(d)}{\phi(p)}   \right \}           \notag\\
&= \sum_{n=1}^\infty n^{\nu/2} K_{\nu}(a\sqrt{nx})\sum_{\substack{d|n\\ (p,d)=(q,n/d)=1}}d^k  \left\{  \cos \left( \frac{2\pi n h_2}{dq}\right) \cos \left( \frac{2\pi dh_1  }{p}\right)+\frac{1}{\phi(p) } \cos \left( \frac{2\pi nh_2  }{dq}\right) 
\right.\notag\\&\left.\ \ \ \ \ \ \ \ \ \ \ \ \ \ \   \ \ \ \ \ \ \ \ \ \  \ \ \ \ \ \ \  \ \ \ \ \ \ \ \ \ \ \ \ \ \ \ \ \ \ \ \ \ \ +\frac{1}{\phi(q) } \cos \left( \frac{2\pi dh_1  }{p}\right)   +\frac{1}{\phi(p) \phi(q) } \right\} \notag\\
&= \sum_{n=1}^\infty n^{\nu/2} K_{\nu}(a\sqrt{nx})\sum_{ d|n }d^k  \left\{  \cos \left( \frac{2\pi n h_2}{dq}\right) \cos \left( \frac{2\pi dh_1 }{p}\right)+\frac{1}{\phi(p) } \cos \left( \frac{2\pi nh_2  }{dq}\right)+\frac{1}{\phi(q) } \cos \left( \frac{2\pi dh_1  }{p}\right) 
\right.\notag\\&\left.\     \ \ \ \   +\frac{1}{\phi(p) \phi(q) } \right\} -\frac{p}{\phi(p) } \sum_{n=1}^\infty n^{\nu/2} K_{\nu}(a\sqrt{nx})\sum_{\substack{d|n\\  p|d}}d^k\left\{ \cos \left( \frac{2\pi nh_2  }{dq}\right)+\frac{1}{\phi(q) }   \right\} \notag\\
 & \ \ \ \   \ \ -\frac{q}{\phi(q) } \sum_{n=1}^\infty n^{\nu/2} K_{\nu}(a\sqrt{nx})\sum_{\substack{d|n\\  q| \frac{n}{d}}}d^k\left\{ \cos \left( \frac{2\pi dh_1  }{p}\right)+\frac{1}{\phi(p) }   \right\}   +\frac{pq}{\phi(p) \phi(q)} \sum_{n=1}^\infty n^{\nu/2} K_{\nu}(a\sqrt{nx})\sum_{\substack{d|n\\ p|d, \  q| \frac{n}{d}}}d^k \notag\\
 &= \sum_{n=1}^\infty n^{\nu/2} K_{\nu}(a\sqrt{nx})\sum_{ d|n }d^k     \cos \left( \frac{2\pi n h_2}{dq}\right) \cos \left( \frac{2\pi dh_1  }{p}\right) \notag\\
&+\left\{ \frac{1}{\phi(p) } \sum_{n=1}^\infty n^{\nu/2} K_{\nu}(a\sqrt{nx})\sum_{ d|n }d^k \cos \left( \frac{2\pi nh_2  }{dq}\right)- \frac{p^{\frac{\nu}{2}+k+1}}{\phi(p) } \sum_{m=1}^\infty m^{\nu/2} K_{\nu}(a\sqrt{mpx})\sum_{ d|m }d^k \cos \left( \frac{2\pi mh_2  }{dq}\right) \right\} \notag\\
&+\left\{\frac{1}{\phi(q) } \sum_{n=1}^\infty n^{\nu/2} K_{\nu}(a\sqrt{nx})\sum_{ d|n }d^k \cos \left( \frac{2\pi dh_1  }{p}\right) - \frac{q^{\frac{\nu}{2}+1}}{\phi(q) } \sum_{m=1}^\infty m^{\nu/2} K_{\nu}(a\sqrt{mqx})\sum_{ d|m }d^k \cos \left( \frac{2\pi dh_1  }{p}\right)   \right\} \notag\\
&+\left\{\frac{1}{\phi(p) \phi(q)} \sum_{n=1}^\infty n^{\nu/2} K_{\nu}(a\sqrt{nx}) \sum_{ d|n }d^k-\frac{p^{\frac{\nu}{2}+k+1}}{\phi(p) \phi(q)} \sum_{m=1}^\infty m^{\nu/2} K_{\nu}(a\sqrt{mpx})\sum_{ d|m }d^k  
\right.\notag\\&\left.\  \ \ \ \ - \frac{q^{\frac{\nu}{2}+1}}{\phi(p) \phi(q)} \sum_{m=1}^\infty m^{\nu/2} K_{\nu}(a\sqrt{mqx})\sum_{ d|m }d^k  + \frac{p^{\frac{\nu}{2}+k+1}q^{\frac{\nu}{2}+1}}{\phi(p) \phi(q) } \sum_{m=1}^\infty m^{\nu/2} K_{\nu}(a\sqrt{mpqx})\sum_{ d|m }d^k    \right\}.     
 \end{align}
Employing Theorem \ref{even2_based} with $\theta=h_2/q$, we evaluate the second and third terms on the right-hand side of \eqref{big1} as follows: 
\begin{align}\label{small1}
					& \frac{1}{\phi(p) } \sum_{n=1}^\infty n^{\nu/2} K_{\nu}(a\sqrt{nx})\sum_{ d/n }d^k \cos \left( \frac{2\pi nh_2  }{dq}\right)- \frac{p^{\frac{\nu}{2}+k+1}}{\phi(p) } \sum_{m=1}^\infty m^{\nu/2} K_{\nu}(a\sqrt{mpx})\sum_{ d|m }d^k \cos \left( \frac{2\pi mh_2  }{dq}\right) \notag\\
					=&\frac{(-1)^{\frac{k+1}{2}}(2 \pi X)^{k+1}}{4q^k \phi(p) } X^\frac{\nu}{2}  \Gamma(\nu+k+1)    \left\{ \sum_{r=1}^\infty  \sum_{\substack{d=1\\  d\equiv \pm h_2(q) } }^\infty  \frac{d^k}{\left(\frac{16\pi^2}{a^2q}\frac{rd}{x}+1\right)^{\nu+k+1} } - \sum_{r=1}^\infty  \sum_{\substack{d=1\\  d\equiv\pm h_2(q) } }^\infty  \frac{d^k}{\left(\frac{16\pi^2}{a^2pq}\frac{rd}{x}+1\right)^{\nu+k+1} }           \right\}
     \notag\\&+\frac{\Gamma(\nu) \zeta(-k)}{4\phi(p) } X^{\frac{\nu}{2}}\left( p^{k+1}-1 \right).
				\end{align}
     Using Theorem \ref{even1_based} with $\theta=h_1/p$, we evaluate the fourth and fifth terms on the right-hand side of \eqref{big1} as follows: 
     \begin{align}\label{small2}
					&\frac{1}{\phi(q) } \sum_{n=1}^\infty n^{\nu/2} K_{\nu}(a\sqrt{nx})\sum_{ d|n }d^k \cos \left( \frac{2\pi dh_1  }{p}\right) - \frac{q^{\frac{\nu}{2}+1}}{\phi(q) } \sum_{m=1}^\infty m^{\nu/2} K_{\nu}(a\sqrt{mqx})\sum_{ d|m }d^k \cos \left( \frac{2\pi dh_1  }{p}\right)\notag\\
					=&\frac{(-1)^{\frac{k+1}{2}}k! }{4(2\pi)^{k+1} } X^\frac{\nu}{2} \Gamma(\nu)\left\{  
					 \zeta(k+1,h_1/p)+\zeta(k+1,1-h_1/p )   \right\} \notag\\
					+& \frac{(-1)^{\frac{k+1}{2}}}{4\phi(q) } X^\frac{\nu}{2}(2 \pi X)^{k+1}\Gamma(\nu+k+1) 
     \left\{ \sum_{d=1}^\infty  \sum_{\substack{r=1\\  r\equiv\pm h_1(p) } }^\infty  \frac{d^k}{\left(\frac{16\pi^2}{a^2p}\frac{rd}{x}+1\right)^{\nu+k+1} } - \frac{1}{q^k}\sum_{d=1}^\infty  \sum_{\substack{r\equiv1\\  r=\pm h_1(p) } }^\infty  \frac{d^k}{\left(\frac{16\pi^2}{a^2pq}\frac{rd}{x}+1\right)^{\nu+k+1} }           \right\}.
				\end{align}
    Using  Proposition \ref{first paper}, we evaluate the last four terms on the right-hand side of \eqref{big1} as follows: 
     \begin{align}\label{small3}
					&\frac{1}{\phi(p) \phi(q)} \sum_{n=1}^\infty n^{\nu/2} K_{\nu}(a\sqrt{nx}) \sum_{ d|n }d^k-\frac{p^{\frac{\nu}{2}+k+1}}{\phi(p) \phi(q)} \sum_{m=1}^\infty m^{\nu/2} K_{\nu}(a\sqrt{mpx})\sum_{ d|m }d^k \notag\\  
  &- \frac{q^{\frac{\nu}{2}+1}}{\phi(p) \phi(q)} \sum_{m=1}^\infty m^{\nu/2} K_{\nu}(a\sqrt{mqx})\sum_{ d|m }d^k  + \frac{p^{\frac{\nu}{2}+k+1}q^{\frac{\nu}{2}+1}}{\phi(p) \phi(q) } \sum_{m=1}^\infty m^{\nu/2} K_{\nu}(a\sqrt{mpqx})\sum_{ d|m }d^k   \notag\\
					& \ \ \ \ =\frac{(-1)^{\frac{k+1}{2}}}{2\phi(p)\phi(q)} X^{\frac{\nu}{2}}\Gamma(\nu+k+1)
     \left(  2\pi X\right)^{k+1} \left\{ \sum_{n=1}^\infty \frac{ \sigma_{k }(n)}{\left(\frac{16\pi^2}{a^2}\frac{n}{x}+1\right)^{\nu+k+1} }-\sum_{n=1}^\infty \frac{ \sigma_{k }(n)}{\left(\frac{16\pi^2}{a^2p}\frac{n}{x}+1\right)^{\nu+k+1} }
			\right.\notag\\&\left.\  \ \ \ \ 
   -\frac{1}{q^k}\sum_{n=1}^\infty \frac{ \sigma_{k }(n)}{\left(\frac{16\pi^2}{a^2q}\frac{n}{x}+1\right)^{\nu+k+1} }+\frac{1}{q^k}\sum_{n=1}^\infty \frac{ \sigma_{k }(n)}{\left(\frac{16\pi^2}{a^2pq}\frac{n}{x}+1\right)^{\nu+k+1} }
   \right\} -\frac{\Gamma(\nu)\zeta(-k)}{4\phi(p)}X^\frac {\nu}{2}(p^{k+1}-1).  
   \end{align}
   
    Substitute \eqref{small1}, \eqref{small2} and \eqref{small3} into the right-hand side of \eqref{big1}, we deduce the left-hand side of \eqref{THM5} as folllows:
    \begin{align}\label{mainp1}
					&\frac{1}{\phi(p)\phi(q)}\sum_{\substack{\chi_2\neq \chi_0\\\chi_2 \ even }}\chi_{2}(h_2)\tau(\Bar{\chi_{2}})\sum_{\substack{\chi_1\neq \chi_0\\\chi_1 \ even }}\chi_{1}(h_1)\tau(\Bar{\chi_{1}})   \sum_{n=1}^\infty \sigma_{k, {\chi_1}, {\chi_2}}(n) n^{\nu/2} K_{\nu}(a\sqrt{nx})        \notag\\
				 &=\sum_{n=1}^\infty n^{\nu/2} K_{\nu}(a\sqrt{nx})\sum_{ d|n }d^k     \cos \left( \frac{2\pi n h_2}{dq}\right) \cos \left( \frac{2\pi dh_1  }{p}\right) \notag\\
				&+\frac{(-1)^{\frac{k+1}{2}}(2 \pi X)^{k+1}}{4q^k \phi(p) } X^\frac{\nu}{2}  \Gamma(\nu+k+1)    \left\{ \sum_{r=1}^\infty  \sum_{\substack{d=1\\  d\equiv\pm h_2(q) } }^\infty  \frac{d^k}{\left(\frac{16\pi^2}{a^2q}\frac{rd}{x}+1\right)^{\nu+k+1} } - \sum_{r=1}^\infty  \sum_{\substack{d=1\\  d\equiv \pm h_2(q) } }^\infty  \frac{d^k}{\left(\frac{16\pi^2}{a^2pq}\frac{rd}{x}+1\right)^{\nu+k+1} }           \right\}  \notag\\
    &+ \frac{(-1)^{\frac{k+1}{2}}k! }{4(2\pi)^{k+1} } X^\frac{\nu}{2} \Gamma(\nu)\left\{  
					 \zeta(k+1,h_1/p)+\zeta(k+1,1-h_1/p )   \right\} \notag\\
					&+ \frac{(-1)^{\frac{k+1}{2}}}{4\phi(q) } X^\frac{\nu}{2}(2 \pi X)^{k+1}\Gamma(\nu+k+1) 
     \left\{ \sum_{d=1}^\infty  \sum_{\substack{r=1\\  r\equiv \pm h_1(p) } }^\infty  \frac{d^k}{\left(\frac{16\pi^2}{a^2p}\frac{rd}{x}+1\right)^{\nu+k+1} } - \frac{1}{q^k}\sum_{d=1}^\infty  \sum_{\substack{r=1\\  r\equiv \pm h_1(p) } }^\infty  \frac{d^k}{\left(\frac{16\pi^2}{a^2pq}\frac{rd}{x}+1\right)^{\nu+k+1} }           \right\} \notag\\
    &+ \frac{(-1)^{\frac{k+1}{2}}}{2\phi(p)\phi(q)} X^{\frac{\nu}{2}}\Gamma(\nu+k+1)
     \left(  2\pi X\right)^{k+1} \left\{ \sum_{n=1}^\infty \frac{ \sigma_{k }(n)}{\left(\frac{16\pi^2}{a^2}\frac{n}{x}+1\right)^{\nu+k+1} }-\sum_{n=1}^\infty \frac{ \sigma_{k }(n)}{\left(\frac{16\pi^2}{a^2p}\frac{n}{x}+1\right)^{\nu+k+1} }
			\right.\notag\\&\left.\  \ \ \ \ 
   -\frac{1}{q^k}\sum_{n=1}^\infty \frac{ \sigma_{k }(n)}{\left(\frac{16\pi^2}{a^2q}\frac{n}{x}+1\right)^{\nu+k+1} }+\frac{1}{q^k}\sum_{n=1}^\infty \frac{ \sigma_{k }(n)}{\left(\frac{16\pi^2}{a^2pq}\frac{n}{x}+1\right)^{\nu+k+1} }
   \right\}. 
    \end{align}
 Next, we examine the right-hand side of \eqref{THM5} of Theorem \ref{botheven_odd}.  We find that
				\begin{align}\label{mainp2}
					&\frac{1}{\phi(p)\phi(q)}\sum_{\substack{\chi_2\neq \chi_0\\\chi_2 \ even }}\chi_{2}(h_2)\tau(\Bar{\chi_{2}})\sum_{\substack{\chi_1\neq \chi_0\\\chi_1 \ even }}\chi_{1}(h_1)\tau(\Bar{\chi_{1}}) \tau(\chi_1)\tau(\chi_2) \sum_{n=1}^\infty  \frac{\sigma_{k,\Bar{\chi_2},\Bar{\chi_1}}(n)}{\left(\frac{16\pi^2}{a^2pq}\frac{n}{x}+1\right)^{\nu+k+1} }          \notag\\
				 &=\frac{pq}{\phi(p)\phi(q)}\sum_{n=1}^\infty\sum_{d|n} \frac{d^k}{\left(\frac{16\pi^2}{a^2pq}\frac{n}{x}+1\right)^{\nu+k+1} }   \left\{\sum_{\chi_2\ even}\chi_{2}(h_2)\Bar{ \chi_2}(d)-\chi_0(d)\right\} \left\{\sum_{\chi_1\ even}\chi_{1}(h_1) \Bar{ \chi_1}(n/d)-\chi_0(n/d) \right\}           \notag\\
      &=\frac{pq}{\phi(p)\phi(q)}\sum_{d,r\geq 1}^\infty  \frac{d^k}{\left(\frac{16\pi^2}{a^2pq}\frac{dr}{x}+1\right)^{\nu+k+1} }   \left\{\sum_{\chi_2\ even}\chi_{2}(h_2)\Bar{ \chi_2}(d)-\chi_0(d)\right\} \left\{\sum_{\chi_1\ even}\chi_{1}(h_1) \Bar{ \chi_1}(r)-\chi_0(r) \right\}           \notag\\
     &= \frac{pq }{4}\sum_{\substack{d=1\\ d \equiv\pm h_2(q)}}^\infty  \sum_{\substack{r=1\\ r \equiv\pm h_1(p)}}^\infty \frac{d^k }{\left(\frac{16\pi^2}{a^2pq}\frac{dr}{x}+1\right)^{\nu+k+1} }   -\frac{pq }{2\phi(p) }\sum_{\substack{d=1\\ d \equiv\pm h_2(q)}}^\infty  \sum_{ \substack{r=1\\  p \nmid r}}^\infty \frac{d^k }{\left(\frac{16\pi^2}{a^2pq}\frac{dr}{x}+1\right)^{\nu+k+1} }   
					 \notag\\
					&   -  \frac{pq }{2\phi(q)}\sum_{\substack{d=1\\  q \nmid d}}^\infty  \sum_{\substack{r=1\\ r \equiv\pm h_1(p)}}^\infty \frac{d^k }{\left(\frac{16\pi^2}{a^2pq}\frac{dr}{x}+1\right)^{\nu+k+1} } 
					 +   \frac{pq }{\phi(p)\phi(q)}\sum_{\substack{d=1\\ q \nmid d}}^\infty \sum_{\substack{r=1\\ p \nmid r}}^\infty \frac{d^k }{\left(\frac{16\pi^2}{a^2pq}\frac{dr}{x}+1\right)^{\nu+k+1} }       \notag\\
					&= \frac{pq }{4}\sum_{\substack{d=1\\ d \equiv\pm h_2(q)}}^\infty  \sum_{\substack{r=1\\ r \equiv\pm h_1(p)}}^\infty \frac{d^k }{\left(\frac{16\pi^2}{a^2pq}\frac{dr}{x}+1\right)^{\nu+k+1} } \notag\\
					& +\left\{ -\frac{pq }{2\phi(p) }\sum_{\substack{d=1\\ d \equiv\pm h_2(q)}}^\infty  \sum_{ r=1 }^\infty \frac{d^k }{\left(\frac{16\pi^2}{a^2pq}\frac{dr}{x}+1\right)^{\nu+k+1} }   
					+  \frac{pq }{2\phi(p) }\sum_{\substack{d=1\\ d \equiv\pm h_2(q)}}^\infty  \sum_{ r=1 }^\infty \frac{d^k }{\left(\frac{16\pi^2}{a^2q}\frac{dr}{x}+1\right)^{\nu+k+1} }  \right\}    \notag\\   \notag\\
					& -  \frac{pq }{2\phi(q)}\sum_{ d=1 }^\infty  \sum_{\substack{r=1\\ r \equiv \pm h_1(p)}}^\infty \frac{d^k }{\left(\frac{16\pi^2}{a^2pq}\frac{dr}{x}+1\right)^{\nu+k+1} } 
					+ \frac{pq^{k+1} }{2\phi(q)}\sum_{ d=1 }^\infty  \sum_{\substack{r=1\\ r \equiv\pm h_1(p)}}^\infty \frac{d^k }{\left(\frac{16\pi^2}{a^2p}\frac{dr}{x}+1\right)^{\nu+k+1} }  
					\notag\\
     &+ \left\{   \frac{pq }{\phi(p)\phi(q)}\sum_{ d=1 }^\infty  \sum_{ r=1}^\infty \frac{d^k }{\left(\frac{16\pi^2}{a^2pq}\frac{dr}{x}+1\right)^{\nu+k+1} }      -  \frac{pq^{k+1} }{\phi(p)\phi(q)}\sum_{ d=1 }^\infty  \sum_{ r=1}^\infty \frac{d^k }{\left(\frac{16\pi^2}{a^2p}\frac{dr}{x}+1\right)^{\nu+k+1} }         
					\right.\notag\\&\left.\ 
					- \frac{pq }{\phi(p)\phi(q)}\sum_{ d=1 }^\infty  \sum_{ r=1}^\infty \frac{d^k }{\left(\frac{16\pi^2}{a^2q}\frac{dr}{x}+1\right)^{\nu+k+1} }                 
					+ \frac{pq^{k+1} }{\phi(p)\phi(q)}\sum_{ d=1 }^\infty  \sum_{ r=1}^\infty \frac{d^k }{\left(\frac{16\pi^2}{a^2}\frac{dr}{x}+1\right)^{\nu+k+1} }  \right\}.   
					 \end{align}   
			 Multiplying equation \eqref{mainp2} by $\frac{(-1)^{\frac{k+1}{2}}}{2p} X^{\frac{\nu}{2}+k+1}(\frac{2\pi}{q})^{k+1}\Gamma(\nu+k+1)$, we obtain the   right-hand side of \eqref{THM5}. We get the result by equating the resulting expression with \eqref{mainp1}. 
 
\textbf{ Theorem \ref{botheven_odd2_based}  $\Rightarrow$ Theorem \ref{botheven_odd} } Here, we consider $\chi_1$ and $\chi_2$ to be non-principal even primitive character modulo $p$ and $q$, respectively. Similarly, as above, one can prove this result.  
				\end{proof}


\begin{proof}[Theorem \rm{\ref{even-odd1_based}} and its equivalence with  Theorem \rm{\ref{even-odd}}][] The proof is similar to Theorem \ref{botheven_odd1_based}, so we skip the proof. But we would show its equivalence with Theorem \ref{even-odd}. 

\textbf{ Theorem \ref{even-odd}  $\Rightarrow$ Theorem \ref{even-odd1_based} }
Let us  multiply the identity \eqref{THM6} of Theorem \ref{even-odd} by $ \chi_1(h_1)  \tau(\bar{\chi_1}) /\phi(p)$  and $ \chi_2(h_2)  \tau(\bar{\chi_2}) /i\phi(q)$, then  take sum on non-principal even primitive characters $\chi_1$ modulo $p$ and odd primitive characters $\chi_2$ modulo $q$. So, from the left-hand side, we get
	  \begin{align} 
					&\frac{1}{i\phi(p)\phi(q)}\sum_{\chi_2\  odd}\chi_{2}(h_2)\tau(\Bar{\chi_{2}})\sum_{\substack{\chi_1\neq \chi_0\\\chi_1 \ even}}\chi_{1}(h_1) \tau(\Bar{\chi_{1}})  \sum_{n=1}^\infty \sigma_{k,\chi_1,\chi_2}(n)n^{\frac{\nu}{2}}K_\nu(a\sqrt{nx})\notag\\
					=&\frac{1}{i\phi(p)\phi(q)}\sum_{n=1}^\infty  n^{\frac{\nu}{2}}K_\nu(a\sqrt{nx})\sum_{d/n}d^k\sum_{\chi_2\  odd}\chi_{2}(h_2)\chi_{2}(n/d)\tau(\Bar{\chi_{2}})\sum_{\substack{\chi_1\neq \chi_0\\\chi_1 \ even}}\chi_{1}(h_1)\chi_{1}(d)\tau(\Bar{\chi_{1}})            \notag\\
					=&\sum_{n=1}^\infty  n^{\frac{\nu}{2}}K_\nu(a\sqrt{nx})\sum_{d|n}d^k \sin \left( \frac{2\pi nh_2}{dq}\right) \left\{ \frac{1}{\phi(p)}\sum_{ \chi_1 \ even}\chi_{1}(h_1)\chi_{1}(d)\tau(\Bar{\chi_{1}})+\frac{\chi_0(d)}{\phi(p)}  \right\} \notag\\
					=&\sum_{n=1}^\infty  n^{\frac{\nu}{2}}K_\nu(a\sqrt{nx})\sum_{\substack{d|n\\p\nmid d}}d^k \sin \left( \frac{2\pi nh_2}{dq}\right) \left\{ \cos \left( \frac{2\pi dh_1 }{p} \right)+\frac{1}{\phi(p)}  \right\} \notag\\
					=&\sum_{n=1}^\infty  n^{\frac{\nu}{2}}K_\nu(a\sqrt{nx})\sum_{d|n}d^k \sin \left( \frac{2\pi nh_2}{dq}\right) \left\{ \cos \left( \frac{2\pi dh_1 }{p} \right)+\frac{1}{\phi(p)}   \right\}
     \notag\\ &-\sum_{n=1}^\infty  n^{\frac{\nu}{2}}K_\nu(a\sqrt{nx})\sum_{\substack{d|n\\p|d}}d^k \sin \left( \frac{2\pi nh_2}{dq}\right) \left\{ 1+\frac{1}{\phi(p)}   \right\} \notag\\
					=&\sum_{n=1}^{\infty}  n^{\nu/2} K_{\nu}(a\sqrt{nx}) \sum_{d|n}d^k  \cos \left( \frac{2\pi dh_1 }{p} \right)\sin \left( \frac{2\pi nh_2}{dq}\right) + \frac{1}{\phi(p)} \sum_{n=1}^\infty   n^{\frac{\nu}{2}}K_\nu(a\sqrt{nx}) \sum_{d|n}d^k \sin \left( \frac{2\pi nh_2}{dq}\right) \notag\\
					&-\frac{p^{\frac{\nu}{2}+k+1}}{\phi(p)} \sum_{m=1}^\infty   m^{\frac{\nu}{2}}K_\nu(a\sqrt{pmx}) \sum_{d|m}d^k \sin \left( \frac{2\pi mh_2}{dq}\right),  \label{evenodd4.1}
				\end{align}
			where we used the  \eqref{cos}, \eqref{sin} in the  penultimate step. Applying Theorem \ref{odd2_based} for the last two terms of the right-hand side of \eqref{evenodd4.1}.  So, we rewrite the right-hand side of \eqref{evenodd4.1} as follows:
				\begin{align}
					&\sum_{n=1}^{\infty}  n^{\nu/2} K_{\nu}(a\sqrt{nx}) \sum_{d|n}d^k  \cos \left( \frac{2\pi dh_1 }{p} \right)\sin \left( \frac{2\pi nh_2}{dq}\right)   \notag \\& +\frac{(-1)^{\frac{k}{2}}}{4\phi(p)}  X^\frac{\nu}{2}(2 \pi X)^{k+1} \Gamma(\nu+k+1)\sum_{r=1}^{\infty}  \sum_{ m=0}^{\infty}\left(  \frac{  (m+h_2/q)^k  }{(1+\frac{16\pi^2r}{a^2x}(m+h_2/q))^{1+\nu+k} } - \frac{ (m+h_2/q)^k    }{(1+\frac{16\pi^2r}{a^2x}(m+1-h_2/q))^{1+\nu+k} } \right)\notag \\&-
					\frac{(-1)^{\frac{k}{2}}}{4\phi(p)}  X^\frac{\nu}{2}(2 \pi X)^{k+1}   \Gamma(\nu+k+1)\sum_{r=1}^{\infty}  \sum_{ m=0}^{\infty}\left(  \frac{  (m+h_2/q)^k  }{(1+\frac{16\pi^2r}{a^2px}(m+h_2/q))^{1+\nu+k} } 
					 -\frac{ (m+h_2/q)^k    }{(1+\frac{16\pi^2r}{a^2px}(m+1-h_2/q))^{1+\nu+k} } \right).\label{evenodd4.2}
				\end{align}
  Using \eqref{both}, \eqref{prop} and  \eqref{prop1} in the right-hand side of \eqref{THM6} of Theorem \ref{even-odd}, we have
			 \begin{align} 
					&-\frac{(-1)^{\frac{k}{2}} X^{\frac{\nu}{2}}}{2p\phi(p)\phi(q)}\left(\frac{2\pi X}{q}\right)^{k+1}   \Gamma(\nu+k+1) \sum_{n=1}^\infty  \frac{\sigma_{k,\Bar{\chi_2},\Bar{\chi_1}}(n)}{\left(\frac{16\pi^2}{a^2pq}\frac{n}{x}+1\right)^{\nu+k+1} } \sum_{\chi_2\  odd}\chi_{2}(h_2) \tau(\chi_2)\tau(\Bar{\chi_{2}})\sum_{\substack{\chi_1\neq \chi_0\\\chi_1 \ even}}\chi_{1}(h_1)  \tau(\chi_1) \tau(\Bar{\chi_{1}}) \notag \\ 
					&=\frac{(-1)^{\frac{k}{2}}pq }{2p\phi(p)\phi(q)}\left(\frac{2\pi X}{q}\right)^{k+1}   X^{\frac{\nu}{2}}\Gamma(\nu+k+1) \sum_{n=1}^\infty  \frac{ \sum_{d/n}d^k}{\left(\frac{16\pi^2}{a^2pq}\frac{n}{x}+1\right)^{\nu+k+1} } \sum_{\chi_2\  odd}\chi_{2}(h_2)\bar{\chi}_{2}(d) \notag\\
     & \ \ \ \ \ \ \ \   \ \ \ \ \ \ \ \   \ \ \ \ \ \ \ \   \ \ \ \ \ \ \ \   \ \ \ \ \ \ \ \   \ \ \ \ \ \ \ \   \ \ \ \ \ \ \ \ \times 
     \left\{ \sum_{ \chi_1 \ even}\chi_{1}(h_1) \bar{\chi}_{1}(n/d)- \bar{\chi}_0(n/d)\right\}  \notag \\ 
					&=\frac{(-1)^{\frac{k}{2}}q}{2\phi(p)\phi(q)}\left(\frac{2\pi X}{q}\right)^{k+1}   X^{\frac{\nu}{2}}\Gamma(\nu+k+1) \sum_{d=1}^\infty \sum_{r=1}^\infty \frac{ d^k}{\left(\frac{16\pi^2}{a^2pq}\frac{dr}{x}+1\right)^{\nu+k+1} } \sum_{\chi_2\  odd}\chi_{2}(h_2)\bar{\chi}_{2}(d) \notag\\
     & \ \ \ \ \ \ \ \   \ \ \ \ \ \ \ \   \ \ \ \ \ \ \ \   \ \ \ \ \ \ \ \   \ \ \ \ \ \ \ \   \ \ \ \ \ \ \ \   \ \ \ \ \ \ \ \ \times \left\{ \sum_{ \chi_1 \ even}\chi_{1}(h_1) \bar{\chi}_{1}(r)- \bar{\chi}_0(r)\right\}  \notag \\
					&=\frac{(-1)^{\frac{k}{2}}\left( 2\pi X\right)^{k+1}}{8 }    X^{\frac{\nu}{2}}\Gamma(\nu+k+1) \sum_{m,n\geq 0}^\infty  \left\{\frac{(n+h_2/q)^k}{\left(\frac{16\pi^2}{a^2}\frac{(n+h_2/q)(m+h_1/p)}{x}+1\right)^{\nu+k+1} } 
     \right.\notag\\&\left.\ \  \ \ \ \ \ \ \ \  \ \ \ \ \ \ \ \  \ \ \ \ \ \ \ \  \ \ \ \ \ \ \ \ -\frac{(n+1-h_2/q)^k}{\left(\frac{16\pi^2}{a^2}\frac{(n+1-h_2/q)(m+h_1/p)}{x}+1\right)^{\nu+k+1} } 
					\right.\notag\\&\left.\ \  \ \ \ \ \ \ \ \  -\frac{(n+1-h_2/q)^k}{\left(\frac{16\pi^2}{a^2}\frac{(n+1-h_2/q)(m+1-h_1/p)}{x}+1\right)^{\nu+k+1} }
					 +\frac{(n+h_2/q)^k}{\left(\frac{16\pi^2}{a^2}\frac{(n+h_2/q)(m+1-h_1/p)}{x}+1\right)^{\nu+k+1} }           \right\}    \notag \\ 
					&-\frac{(-1)^{\frac{k}{2}}q}{2\phi(p)\phi(q)}\left(\frac{2\pi X}{q}\right)^{k+1}   X^{\frac{\nu}{2}}\Gamma(\nu+k+1) \sum_{\substack{d,r\geq 1\\p\nmid r}}   \frac{ d^k}{\left(\frac{16\pi^2}{a^2pq}\frac{dr}{x}+1\right)^{\nu+k+1} } \sum_{\chi_2\  odd}\chi_{2}(h_2)\bar{\chi}_{2}(d).  
					\label{evenodd4.3}        \end{align}
				By \eqref{prop}, we examine the last term of right-hand side of \eqref{evenodd4.3}, we find
				\begin{align} 
					&\frac{1}{\phi(q)}\sum_{\substack{d,r\geq 1\\p\nmid r}}   \frac{ d^k}{\left(\frac{16\pi^2}{a^2pq}\frac{dr}{x}+1\right)^{\nu+k+1} } \sum_{\chi_2\  odd}\chi_{2}(h_2)\bar{\chi}_{2}(d) \notag \\ 
					= &\frac{1}{\phi(q)}\sum_{ d,r\geq 1 }   \frac{ d^k}{\left(\frac{16\pi^2}{a^2pq}\frac{dr}{x}+1\right)^{\nu+k+1} } \sum_{\chi_2\  odd}\chi_{2}(h_2)\bar{\chi}_{2}(d) -\frac{1}{\phi(q)}\sum_{ d,r\geq 1 }   \frac{ d^k}{\left(\frac{16\pi^2}{a^2q}\frac{dr}{x}+1\right)^{\nu+k+1} } \sum_{\chi_2\  odd}\chi_{2}(h_2)\bar{\chi}_{2}(d) \notag \\ 
					=&\frac{q^k}{2}\sum_{r=1}^\infty  \sum_{m=0}^\infty \left\{\frac{ (m+h_2/q)^k}{\left(\frac{16\pi^2}{a^2p}\frac{(m+h_2/q)r}{x}+1\right)^{\nu+k+1} } -  \frac{ (m+1-h_2/q)^k}{\left(\frac{16\pi^2}{a^2p}\frac{(m+1-h_2/q)r}{x}+1\right)^{\nu+k+1} }  \right\}\notag \\
					&-\frac{q^k}{2}\sum_{r=1}^\infty  \sum_{m=0}^\infty \left\{\frac{ (m+h_2/q)^k}{\left(\frac{16\pi^2}{a^2}\frac{(m+h_2/q)r}{x}+1\right)^{\nu+k+1} } -  \frac{ (m+1-h_2/q)^k}{\left(\frac{16\pi^2}{a^2}\frac{(m+1-h_2/q)r}{x}+1\right)^{\nu+k+1} }  \right\}.
					\label{evenodd4.4}\end{align}
				Equating \eqref{evenodd4.2}   \eqref{evenodd4.3}, and using \eqref{evenodd4.4} we get the result. 
				 
     \textbf{ Theorem \ref{even-odd1_based} $\Rightarrow$ Theorem  \ref{even-odd}}   Let $\theta=h_1/p$, $\psi=h_2/q$, and let $\chi_1$   be an even primitive character modulo $p$ and $\chi_2$ be an odd primitive character modulo $q$. Multiplying the identity in Theorem \ref{even-odd1_based}
				 by $\bar{\chi_1}(h_1)\bar{\chi_2}(h_2)/\tau(\bar{\chi_1})\tau(\bar{\chi_2})$, and then summing on $h_1$ and $h_2$, $ 0<h_1<p$,  $0<h_2<q$, one can prove that Theorem \ref{even-odd1_based}   implies Theorem \ref{even-odd}.
  \end{proof}

\begin{proof}[Theorem \rm{\ref{even-odd2_based}} and its equivalence with Theorem \rm{\ref{even-odd}} ][]
 The proof is  similar to the  previous one.
\end{proof}

\section{Proof of Cohen type identities }\label{proof of cohen identities...}

 \begin{proof}[Theorem \rm{\ref{oddcohen based}} and its equivalence with   Theorem \rm{\ref{oddcohen}}][]
This theorem can be proved using the methodology presented in \cite[p. 13, Section 3]{devika2023}. 
\\
Next, we will establish the equivalence of Theorem \ref{oddcohen based} and Theorem \ref{oddcohen}.\\
 \textbf{ Theorem \ref{oddcohen} $\Rightarrow$ Theorem  \ref{oddcohen based} }
 The double series on the right-hand side of the identity in Theorem \ref{oddcohen based} converges absolutely and uniformly on any compact interval for $\theta\in  (0,1)$. Therefore, it is sufficient to prove the theorem for $\theta=h/q$, where $q$ is prime and $0<h<q.$  Now we multiply the identity \eqref{THM7} of Theorem \eqref{oddcohen} by $ \chi(h)  \tau(\bar{\chi}) /i\phi(q)$, then take sum on odd primitive character $\chi$ modulo $q$. The left-hand side of the identity \eqref{THM7}  becomes \begin{align}
					 \frac{8 \pi x^{\nu/2}}{i\phi(q)} \sum_{\chi \ odd } \chi(h)  \tau(\bar{\chi})   \sum_{n=1}^\infty\sigma_{-\nu,\chi}(n)n^{\frac{\nu}{2}}K_\nu(4\pi\sqrt{nx}) = &\frac{8 \pi x^{\nu/2}}{i\phi(q)}\sum_{n=1}^\infty n^{\frac{\nu}{2}}K_\nu(4\pi\sqrt{nx}) \sum_{d|n}d^{-\nu} \sum_{\chi \ odd } \chi(d)  \chi(h)  \tau(\bar{\chi}) \notag\\
				 =&8 \pi x^{\nu/2} \sum_{n=1}^{\infty} n^{\nu/2} K_{\nu}(4\pi\sqrt{nx}) \sum_{d|n}d^{-\nu} \sin \left(\frac{2\pi d h}{q}\right),\label{ZX1}
				\end{align}
				where we have used the  identity \eqref{cos}.  Employing the functional equation \eqref{ll(s)}, the right-hand side of \eqref{THM7} transforms into
      \begin{align}
					& \frac{1}{i\phi(q)} \sum_{\chi \ odd } \chi(h)  \tau(\bar{\chi}) \ 8 \pi x^{\nu/2}  \sum_{n=1}^\infty\sigma_{-\nu,\chi}(n)n^{\frac{\nu}{2}}K_\nu(4\pi\sqrt{nx})\notag\\ 
					&=- \frac{\pi}{ q^{\nu-1}\sin(\frac{\pi\nu}{2})} \frac{1}{\phi(q)} \sum_{\chi \ odd } \chi(h)   L(1-\nu, \bar{\chi}) 
     +  \frac{1}{ x\ q^{\nu}\cos(\frac{\pi\nu}{2})} \frac{1}{\phi(q)} \sum_{\chi \ odd } \chi(h)   L(-\nu, \bar{\chi}) \notag\\ 
    &+ \frac{q^{1-\nu}  }{ \phi(q)}\sum_{\chi \ odd } \chi(h)   \left\{ \frac{2\zeta(\nu+1)L(1,{\bar{\chi}})(qx)^{\nu}} {\cos(\frac{\pi\nu}{2}) }  - \frac{2}{  \cos \left(\frac{\pi \nu}{2}\right)}\sum_{j=1}^{N} \zeta(2j)\ L(2j-\nu, \bar{\chi})(qx)^{2j-1} 
    \right.\notag\\&\left.\ \ 
     -\frac{2(qx)^{2N+1} }{ \cos\left(\frac{\pi \nu}{2}\right)   }\sum_{n=1}^{\infty}\bar{\sigma}_{-\nu, \bar{\chi}}(n)  \frac{ \left( n^{\nu+1-2N}-(qx)^{\nu+1-2N}\right)}{ n\ (n^2-(qx)^2)}
    \right\}.  \label{X1}
				\end{align}
    We observe using \eqref{Hurwitz} and \eqref{prop}
				\begin{align}
					\frac{1}{\phi(q)} \sum_{\chi \ odd } \chi(h) L(s, \bar{\chi})  = \frac{1}{2q^{s}}\left(\zeta(s, h/q) - \zeta(s, 1-h/q)\right).\label{X2}
				\end{align} 
     Lastly, the infinite series in \eqref{X1} becomes
				\begin{align}
					& \frac{1}{\phi(q)} \sum_{\chi \ odd } \chi(h) \sum_{n=1}^{\infty}\bar{\sigma}_{-\nu, \bar{\chi}}(n)  \frac{  n^{\nu+1-2N}-(qx)^{\nu+1-2N} }{ n\ (n^2-(qx)^2)}\notag\\
      &= \frac{1}{\phi(q)} \sum_{n=1}^{\infty} 
       \frac{  n^{\nu+1-2N}-(qx)^{\nu+1-2N} }{ n\ (n^2-(qx)^2)}  \sum_{d/n} d^{-\nu}  \sum_{\chi \ odd } \chi(h) \bar{\chi}\left(\frac{n}{d}\right)  \notag\\
       &=\frac{1}{2}\sum_{d=1}^{\infty}   d^{-\nu} \left\{\sum_{\substack{ m=1 \\ m\equiv  h(q)} }^\infty \frac{   (dm)^{\nu+1-2N}-(qx)^{\nu+1-2N} }{ dm\ (d^2m^2-(qx)^2)}-\sum_{\substack{ m=1 \\ m\equiv - h(q)} }^\infty \frac{   (dm)^{\nu+1-2N}-(qx)^{\nu+1-2N} }{ dm\ (d^2m^2-(qx)^2)}  \right\}\notag\\
        &=\frac{1}{2}\sum_{d=1}^{\infty}   d^{-\nu} \sum_{  r=0}^\infty
      \left\{    \frac{   (d(rq+h))^{\nu+1-2N}-(qx)^{\nu+1-2N} }{ d(rq+h)\ (d^2(rq+h)^2-(qx)^2)} 
      -    \frac{   (d(rq+q-h))^{\nu+1-2N}-(qx)^{\nu+1-2N} }{ d(rq+q-h)\ (d^2(rq+q-h)^2-(qx)^2)} 
      \right\}    \notag\\
      &=\frac{q^{\nu-2-2N}}{2}\sum_{d=1}^{\infty}   d^{-\nu-1} \sum_{  r=0}^\infty
      \left\{    \frac{   (d(r+h/q))^{\nu+1-2N}-x^{\nu+1-2N} }{ (r+h/q)\ (d^2(r+h/q)^2-x^2)} 
      -    \frac{   (d(r+1-h/q))^{\nu+1-2N}-x^{\nu+1-2N} }{ (r+1-h/q)\ (d^2(r+1-h/q)^2-x^2)} 
      \right\}.      \label{X3}
   \end{align}   
   Employing \eqref{X2}, \eqref{X3} in \eqref{X1}, we obtain
   \begin{align}\label{ZX2}
					& \frac{1}{i\phi(q)} \sum_{\chi \ odd } \chi(h)  \tau(\bar{\chi}) \ 8 \pi x^{\nu/2}  \sum_{n=1}^\infty\sigma_{-\nu,\chi}(n)n^{\frac{\nu}{2}}K_\nu(4\pi\sqrt{nx})  
					= \frac{1 }{\cos\left(\frac{\pi \nu}{2}\right)}     \zeta(\nu+1) \left(\zeta(1,h/q) - \zeta(1, 1-h/q)\right)x^\nu \notag\\
     & -\frac{\pi }{ 2 \sin\left(\frac{\pi \nu}{2}\right)}      \left(\zeta(1-\nu,\frac{h}{q}) - \zeta(1-\nu, 1-\frac{h}{q})\right)
     +\frac{1 }{ 2 x\cos\left(\frac{\pi \nu}{2}\right)}      \left(\zeta(-\nu,\frac{h}{q}) - \zeta(-\nu, 1-\frac{h}{q})\right) \notag\\
     & -\frac{1 }{  \cos\left(\frac{\pi \nu}{2}\right)}     \sum_{j=1}^N \zeta(2j) \left(\zeta(2j-\nu,\frac{h}{q}) - \zeta(2j-\nu, 1-\frac{h}{q})\right)x^{2j-1} 
       \notag\\
      &-\frac{ x^{2N+1}}{\cos(\frac{\pi\nu}{2}) }\sum_{d=1}^{\infty}   d^{-\nu-1} \sum_{  r=0}^\infty
      \left\{    \frac{   (d(r+h/q))^{\nu+1-2N}-x^{\nu+1-2N} }{ (r+h/q)\ (d^2(r+h/q)^2-x^2)} 
      -    \frac{   (d(r+1-h/q))^{\nu+1-2N}-x^{\nu+1-2N} }{ (r+1-h/q)\ (d^2(r+1-h/q)^2-x^2)} 
      \right\}.  
     \end{align} 
     Equating \eqref{ZX1} and \eqref{ZX2}, we get the result.\\
  \textbf{ Theorem \ref{oddcohen based}  $\Rightarrow$ Theorem \ref{oddcohen} }
      Analogous to proof of Theorem \ref{M1}, one can easily prove Theorem \ref{oddcohen}.
     \end{proof}

\begin{remark}
The other proofs of this section will be similar. We skip the proofs to avoid repetitions.  
\end{remark}
  
   \section{Proof of  Voronoi summation formulas}\label{proof of voronoi...}
 This section is devoted to the proof of  Voronoi-type summation formulas.
\begin{proof}[Theorem \rm{\ref{vor1.1}} and its equivalence with Theorem \ref{voro2}][] This theorem can be deduced using the methodology presented in \cite[p. 18, Section 4]{devika2023}. \\
Next, we establish the equivalence of Theorem \ref{voro2} and Theorem \ref{vor1.1}.\\
 \textbf{ Theorem \ref{voro2}  $\Rightarrow$ Theorem \ref{vor1.1} } 
 It is sufficient to prove the theorem for $\theta=h/q$, where $q$ is prime and $0<h<q$.
Now we multiply both sides of identity \eqref{vv1} in Theorem \ref{voro2} by $ \chi(h)  \tau(\bar{\chi}) /i\phi(q)$, then  sum on odd primitive character $\chi$ modulo $q$,  then
  the left-hand side of \eqref{voro2} becomes
 \begin{align}\label{vor01.1}
					 \frac{1}{i\phi(q)} \sum_{\chi \ odd } \chi(h)  \tau(\bar{\chi})   \sum_{\alpha<j<\beta}  \sigma_{-\nu,\chi}(j) f(j)&= \frac{1}{i\phi(q)}\sum_{\alpha<j<\beta}    \sum_{d| j}d^{-\nu}  \sum_{\chi \ odd } \chi(d)  \chi(h)  \tau(\bar{\chi})f(j) \notag\\
					&=\sum_{\alpha<j<\beta} \sum_{d|j}d^{-\nu}  \sin \left(\frac{2\pi d h}{q}\right)f(j) , 
				\end{align}
where we have used the  identity \eqref{cos}.  
  The right hand side of \eqref{voro2} becomes
\begin{align}\label{vor01.2}
&\frac{ 1 }{i\phi(q)} \sum_{\chi \ odd } \chi(h)  \tau(\bar{\chi})   \sum_{\alpha<j <\beta}    { {\sigma}_{-\nu, \chi }(j)} f(j)   =   \frac{ 1}{i\phi(q)}\sum_{\chi \ odd } \chi(h)  \tau(\bar{\chi})   L(1+\nu,\chi) \int_\alpha^\beta     {f(t) } \mathrm{d}t 
   \notag\\ 
   &-  \frac{2 \pi }{\phi(q)q^{\frac{\nu}{2}}}\sum_{\chi \ odd } \chi(h)        \sum_{n=1}^{\infty}\bar{\sigma}_{-\nu, \bar{\chi }}(n) \ n^{\nu/2} \int_\alpha^\beta     {f(t) } t^{-\frac{\nu}{2}}     \left\{  \left(   \frac{2}{\pi} K_{\nu}\left(4\pi  \sqrt{\frac{nt}{q}}\ \right) + Y_{\nu}\left(4\pi  \sqrt{\frac{nt}{q}}\ \right)   \right)  \sin \left(\frac{\pi \nu}{2}\right) 
   \right.\notag\\&\left.\ \ 
   \ \ \ \ \ \ \ \ \ \ \ \ \ \ \ \ \ \ \ \ \ \ \ \ \ \ \ \ \ \ \ \ \ \  \ \ \ \ \ \ \ \ \ \ \ \ \ \ \ \ \  \ \ \ \ \ \ \ \ \ \ \ \ \ \ \ \ \  \ \   \ \ \ \ \ \ \ \   -  J_{\nu}\left(4\pi  \sqrt{\frac{nt}{q}}\ \right)  \cos \left(\frac{\pi \nu}{2}\right)   \right\}  dt.
\end{align}
Utilising \eqref{ll(s)}, \eqref{Hurwitz}, \eqref{both} and \eqref{prop}, we get
\begin{align}\label{vor01.3}
 \frac{ 1}{i\phi(q)}\sum_{\chi \ odd } \chi(h)  \tau(\bar{\chi})   L(1+\nu,\chi)  =- {(2\pi)^{\nu}} \Gamma(-\nu) \sin\left(\frac{\pi \nu}{2}\right)\left\{   \zeta(-\nu,\frac{h}{q})- \zeta(-\nu,1-\frac{h}{q})\right\}.
    \end{align}
  We define for $y \in \mathbb{R},$
  \begin{align}\label{qq}
        Z_{\nu}\left(y\ \right)=    \left(   \frac{2}{\pi} K_{\nu}\left(y \right) + Y_{\nu}\left(y\right)   \right)  \sin \left(\frac{\pi \nu}{2}\right) 
       -  J_{\nu}\left(y \right)  \cos \left(\frac{\pi \nu}{2}\right).   
    \end{align}
Next, we consider
\begin{align}\label{vor01.4}
   &\frac{1}{\phi(q)}  \sum_{\chi \ odd } \chi(h)        \sum_{n=1}^{\infty}\bar{\sigma}_{-\nu, \bar{\chi }}(n) \ n^{\nu/2} Z_{\nu}\left(4\pi  \sqrt{\frac{nt}{q}}\ \right)\notag\\
   &= \frac{1}{\phi(q)} \sum_{n=1}^{\infty} 
      n^{\nu/2} Z_{\nu}\left(4\pi  \sqrt{\frac{nt}{q}}\ \right)  \sum_{d|n} d^{-\nu}  \sum_{\chi \ odd } \chi(h) \bar{\chi}\left(\frac{n}{d}\right)  \notag\\
     &=\frac{1}{2}\sum_{d=1}^{\infty}   d^{-\nu} \left\{\sum_{\substack{ r=1 \\ r\equiv  h(q)} }^\infty (dr)^{\nu/2} Z_{\nu}\left(4\pi  \sqrt{\frac{drt}{q}}\ \right)  -\sum_{\substack{ r=1 \\ r\equiv  -h(q)} }^\infty (dr)^{\nu/2} Z_{\nu}\left(4\pi  \sqrt{\frac{drt}{q}}\ \right)    \right\}\notag\\
       &=\frac{q^{\frac{\nu}{2}}}{2}\sum_{d=1}^{\infty}   d^{-\frac{\nu}{2} } \sum_{  m=0}^\infty
      \left\{  \left(m+\frac{h}{q}\right)^\frac{\nu}{2} Z_{\nu}\left(4\pi  \sqrt {d\left(m+\frac{h}{q}\right)t} \ \right)- \left(m+1-\frac{h}{q}\right)^\frac{\nu}{2} Z_{\nu}\left(4\pi  \sqrt {d\left(m+1-\frac{h}{q}\right)t} \ \right)   
      \right\}.  
\end{align}
Employing \eqref{vor01.3}, \eqref{qq} \eqref{vor01.4}, we deduce the expression for \eqref{vor01.2}. Then equating the resulting expression with \eqref{vor01.1}, we get the desired result.

\textbf{ Theorem \ref{vor1.1}  $\Rightarrow$ Theorem \ref{voro2} } Let $\theta=h/q$, and let $\chi$ be an odd primitive character modulo $q$. Multiplying the identity \eqref{vv} in Theorem \ref{vor1.1} by $\bar{\chi}(h)/\tau(\bar{\chi})$, and then summing on $h$, $ 0<h<q$,  the left-hand side of the identity \eqref{vv} becomes
     \vspace{-.1cm}
				\begin{align}\label{P1}
				    & \frac{1}{\tau(\bar{\chi})}\sum_{h=1}^{q-1}\bar{\chi}(h)\sum_{\alpha<j <\beta}      f(j) \sum_{d|j}d^{-\nu} \sin \left( 2\pi d   h/q \right)\notag\\
       &= \frac{1}{2 i \tau(\bar{\chi})} \sum_{\alpha<j <\beta}      f(j) \sum_{d|j}d^{-\nu}  \sum_{h=1}^{q-1}\bar{\chi}(h) \left( e^{2\pi i  d   h/q} - e^{-2\pi i  d   h/q}\right)\notag\\
     &= \frac{1}{2 i  } \sum_{\alpha<j <\beta}      f(j) \sum_{d|j}d^{-\nu}    \left(  \chi(d)-\chi(-d)\right)\notag\\
     &=i^{-1}\sum_{\alpha<j <\beta}    { {\sigma}_{-\nu, \chi }(j)} f(j) ,
				\end{align} 
where in the penultimate step, we used \eqref{gauss}. 
Employing \eqref{Hurwitz} and \eqref{both}, we can observe that
 \begin{align}\label{P2}
 \frac{1}{\tau(\bar{\chi})}\sum_{h=1}^{q-1}\bar{\chi}(h)\left(\zeta(-\nu,h/q) - \zeta(-\nu, 1- h/q)\right) =-2 q^{-\nu-1}\tau({\chi})L(-\nu,\Bar{\chi}).
\end{align}  
  Next, we consider
\begin{align}\label{P14}
&\frac{1}{\tau(\bar{\chi})}\sum_{h=1}^{q-1}\bar{\chi}(h)\sum_{d=1}^{\infty}d^{-\nu/2} \sum_{ m=0}^{\infty}\left\{  \left(m+h/q\right)^\frac{\nu}{2} Z_{\nu}\left(4\pi  \sqrt {d\left(m+h/q\right)t} \ \right) 
        - \left(m+1-h/q\right)^\frac{\nu}{2} Z_{\nu}\left(4\pi  \sqrt {d\left(m+1-h/q\right)t} \ \right)   
      \right\}  \notag\\
&=\frac{q^{-\nu/2}}{\tau(\bar{\chi})}\sum_{h=1}^{q-1}\bar{\chi}(h)\sum_{d=1}^{\infty}d^{-\nu/2} \sum_{\substack{r=1 \\ r \equiv h(q)}}^\infty r^{\nu/2} Z_\nu\left(4\pi  \sqrt{\frac{drt}{q}}\ \right)
-\frac{q^{-\nu/2}}{\tau(\bar{\chi})}\sum_{h=1}^{q-1}\bar{\chi}(h)\sum_{d=1}^{\infty}d^{-\nu/2}\sum_{\substack{r=1 \\ r \equiv -h(q)}}^\infty r^{\nu/2} Z_\nu\left(4\pi  \sqrt{\frac{drt}{q}}\ \right)\notag\\
&=\frac{2q^{-\nu/2}}{\tau(\bar{\chi})}  \sum_{d=1}^{\infty}  \sum_{ r=1  }^\infty  d^{-\nu/2} r^{\nu/2} \bar{\chi}(r)  Z_\nu\left(4\pi  \sqrt{\frac{drt}{q}}\ \right)   \notag\\
&=-\frac{2\tau({\chi}) }{q^{1+\nu/2}}   \sum_{ n=1  }^\infty    \Bar{\sigma}_{-\nu,\Bar{\chi} }(n) n^{\nu/2}  Z_\nu\left(4\pi  \sqrt{\frac{nt}{q}}\ \right)  .   
\end{align}
Substituting \eqref{P1}, \eqref{P14} and utilising \eqref{P2} in \eqref{vv}, one can complete the proof.
 \end{proof}
\begin{proof}[Theorem \rm{\ref{vor1.3}} and its equivalence with Theorem \rm{\ref{vore2}}][] The proof is similar to the previous one. Next, we demonstrate the equivalence of Theorem \ref{vore2} and Theorem \ref{vor1.3}.\\
 \textbf{ Theorem \ref{vore2} $\Rightarrow$ Theorem \ref{vor1.3} } It is sufficient to prove the theorem for rationals $\theta=h/q$, $0<h<q$. 
First, we  consider the following
\begin{align}\label{fan1}
 &  \sum_{\alpha<j <\beta}  \sum_{d|j}  {d}^{-\nu} \cos\left(\frac{2 \pi d h }{q}\right)f(j)   =\sum_{\alpha<j <\beta}   f(j) \left(\sum_{\substack{d|j\\ q|d}} d^{-\nu}  +\sum_{\substack{d|j\\ q\nmid d}}d^{-\nu} \cos \left(\frac{ 2\pi dh }{q}\right)\right)\notag\\
 &=q^{-\nu}\sum_{\frac{\alpha}{q}<m <\frac{\beta}{q}}  f(qm)\sum_{d|m}  d^{-\nu}
+\sum_{\alpha<j <\beta}  f(j)\sum_{\substack{d|j\\ q\nmid d}} \frac{d^{-\nu}}{\phi(q)}\sum_{\chi \ even } \chi(d)  \chi(h)  \tau(\bar{\chi}) \notag\\  
&=q^{-\nu}\sum_{\frac{\alpha}{q}<m <\frac{\beta}{q}} f(qm)\sum_{d|m}  d^{-\nu} 
-\sum_{\alpha<j <\beta}   f(j)  \sum_{\substack{d|j\\ q\nmid d}} \frac{d^{-\nu}}{\phi(q)}\chi_0(d) 
+ \sum_{\alpha<j <\beta}  f(j)  \sum_{\substack{d|j\\ q\nmid d}} \frac{d^{-\nu}}{\phi(q)} \sum_{\substack{\chi \neq \chi_0\\\chi\ even}}\chi(d)  \chi(h)  \tau(\bar{\chi})  \notag\\ 
&=q^{-\nu}\sum_{\frac{\alpha}{q}<m <\frac{\beta}{q}} f(qm)\sum_{d|m}  d^{-\nu}
-\sum_{\alpha<j <\beta}   f(j)   \frac{1}{\phi(q)} \left( \sum_{d|j}d^{-\nu} -\sum_{\substack{d|j\\ q| d}}d^{-\nu} \right) \notag\\
&\ \ \ \ +\frac{1}{\phi(q) }\sum_{\substack{\chi \neq \chi_0\\\chi even}}   \chi(h)  \tau(\bar{\chi})\sum_{\alpha<j <\beta}  \sigma_{-\nu,\chi} (j)  f(j)   \notag\\
&=\frac{q^{1-\nu}}{\phi(q)} \sum_{\frac{\alpha}{q}<m <\frac{\beta}{q}} \sigma_{-\nu } (m)f(qm) 
-\frac{1}{\phi(q)} \sum_{\alpha<j <\beta} \sigma_{-\nu } (j)  f(j)         +\frac{1}{\phi(q) }\sum_{\substack{\chi \neq \chi_0\\\chi\  even}}   \chi(h)  \tau(\bar{\chi})\sum_{\alpha<j <\beta} \sigma_{-\nu,\chi} (j)  f(j).   
  \end{align}

We will first estimate the sum of the right-hand side of \eqref{fan1}. 
Applying Proposition \ref{vorlemma1} with $f(x)=f(qx)$, then simplifying, we obtain
\begin{align}\label{fan2}    
          & \frac{q^{1-\nu}}{\phi(q)} \sum_{\frac{\alpha}{q}<m <\frac{\beta}{q}}   {\sigma}_{-\nu}(m)f(qm)   = \frac{1}{\phi(q)}  \int_{\alpha} ^{  \beta   }f(t) \left\{   \zeta(1-\nu) \ t^{-\nu}   +  q^{-\nu}  \zeta(\nu+1) \         \right\}  dt  \notag\\
+& \frac{2\pi q^{-\frac{\nu}{2}}}{\phi(q)} \sum_{n=1}^{\infty} \sigma_{-\nu}(n) n^{\nu/2}    \int_{\alpha} ^{  \beta   }f(t)  (t)^{-\frac{\nu}{2}}   \left\{  \left(   \frac{2}{\pi} K_{\nu}(4\pi   \sqrt{\frac{nt}{q}}) - Y_{\nu}(4\pi   \sqrt{\frac{nt}{q}})\right) \cos \left(\frac{\pi \nu}{2}\right) 
  -  J_{\nu}(4\pi   \sqrt{\frac{nt}{q}})  \sin \left(\frac{\pi \nu}{2}\right)  \right\}  dt.
           \end{align}
Similarly, one can get the second sum  of the right-hand side of \eqref{fan1}. Now, we focus on the last sum of  the right-hand side of \eqref{fan1}. By  \eqref{fann} of Theorem \ref{vore2}, the last sum of the right-hand side of \eqref{fan1} becomes

\begin{align}\label{fan3}
&\frac{1}{\phi(q) }\sum_{\substack{\chi \neq \chi_0\\\chi\  even}}   \chi(h)  \tau(\bar{\chi})\sum_{\alpha<j <\beta} \sigma_{-\nu,\chi} (j)  f(j)\nonumber\\
&  =    \frac{1}{\phi(q) }\sum_{\substack{\chi \neq \chi_0\\\chi\  even}}   \chi(h)  \tau(\bar{\chi})  L(1+\nu,\chi) \int_\alpha^\beta     {f(t) } \mathrm{d}t+ 2 \pi  \frac{q^{-\frac{\nu}{2}} }{\phi(q) }\sum_{\substack{\chi \neq \chi_0\\\chi\  even}}   \chi(h)     \notag\\  
 &   \times   \sum_{n=1}^{\infty}\bar{\sigma}_{-\nu, \bar{\chi }}(n) \ n^{\nu/2} \int_\alpha^\beta \frac{f(t) }{t^{\frac{\nu}{2}}}    \left\{  \left(   \frac{2}{\pi} K_{\nu}\left(4\pi  \sqrt{\frac{nt}{q}}\ \right) - Y_{\nu}\left(4\pi  \sqrt{\frac{nt}{q}}\ \right)   \right) \cos \left(\frac{\pi \nu}{2}\right) 
    -  J_{\nu}\left(4\pi  \sqrt{\frac{nt}{q}}\ \right)  \sin \left(\frac{\pi \nu}{2}\right)  \right\}  dt.
\end{align}
We define for $y \in \mathbb{R},$
  \begin{align}\label{qw}
        W_{\nu}\left(y\ \right)=    \left(   \frac{2}{\pi} K_{\nu}\left(y \right) - Y_{\nu}\left(y\right)   \right)  \cos \left(\frac{\pi \nu}{2}\right) 
       -  J_{\nu}\left(y \right)  \sin \left(\frac{\pi \nu}{2}\right).   
    \end{align}
 Now, we first observe 
\begin{align}\label{fan4}
 & \frac{1}{\phi(q) }\sum_{\substack{\chi \neq \chi_0\\\chi\  even}}   \chi(h)   \sum_{n=1}^{\infty}\bar{\sigma}_{-\nu, \bar{\chi }}(n) \ n^{\nu/2}  W_{\nu}\left(4\pi  \sqrt{\frac{nt}{q}}\ \right)  \notag\\  
 &=  \frac{1}{\phi(q) } \sum_{n=1}^{\infty} \sum_{d/n}d^{-\nu} n^{\nu/2}W_{\nu}\left(4\pi  \sqrt{\frac{nt}{q}}\ \right)   \left(\sum_{  \chi \  even }   \chi(h) \Bar{\chi}(n/d)-\chi_0(n/d)\right) \notag\\ 
 &= \frac{1}{\phi(q) } \sum_{d,r=1}^{\infty}d^{-\nu/2} r^{\nu/2}W_{\nu}\left(4\pi  \sqrt{\frac{drt}{q}}\ \right)\sum_{  \chi \  even }   \chi(h) \Bar{\chi}(r)- \frac{1}{\phi(q) } \sum_{n=1}^{\infty}  n^{\nu/2}W_{\nu}\left(4\pi  \sqrt{\frac{nt}{q}}\ \right)  \left(\sigma_{-\nu}(n)- \sigma_{-\nu}( n/q)\right) \notag\\ 
 &= \frac{1}{2}\sum_{d=1}^\infty  d^{-\nu/2} \sum_{\substack{ r=1 \\ r\equiv \pm  h(q)} }^\infty   r^{\nu/2}  W_{\nu}\left(4\pi  \sqrt{\frac{drt}{q}}\ \right) -\frac{1}{\phi(q) } \sum_{n=1}^{\infty}\sigma_{-\nu}(n)  n^{\nu/2}W_{\nu}\left(4\pi  \sqrt{\frac{nt}{q}}\ \right)   \notag\\ 
 &  \ \ + \frac{q^{\nu/2}}{\phi(q) } \sum_{n=1}^{\infty} \sigma_{-\nu}(n) n^{\nu/2}W_{\nu}\left(4\pi  \sqrt{ {nt} }\ \right)\notag\\ 
 &= \frac{q^{\nu/2}}{2}\sum_{d=1}^\infty  d^{-\nu/2}  \sum_{m=0}^\infty \left\{ \left(m+\frac{h}{q}\right)^{\nu/2}W_{\nu}\left(4\pi  \sqrt{d\left(m+\frac{h}{q}\right)t} \right)   + \left(m+1-\frac{h}{q}\right)^{\nu/2}W_{\nu}\left(4\pi  \sqrt{d\left(m+1-\frac{h}{q}\right)t} \right)\right\}  \notag\\ 
 &  \ \ -\frac{1}{\phi(q) } \sum_{n=1}^{\infty}\sigma_{-\nu}(n)  n^{\nu/2}W_{\nu}\left(4\pi  \sqrt{\frac{nt}{q}}\ \right)   + \frac{q^{\nu/2}}{\phi(q) } \sum_{n=1}^{\infty} \sigma_{-\nu}(n) n^{\nu/2}W_{\nu}\left(4\pi  \sqrt{ {nt} }\ \right),
\end{align}
where in the penultimate step, we used \eqref{prop1}
By utilising \eqref{1st_use}, \eqref{ll(s)}, \eqref{Hurwitz}, \eqref{both} and \eqref{prop1}, we arrive at
\begin{align}\label{fan5} 
 \frac{ 1}{\phi(q)}\sum_{\substack{\chi \neq \chi_0\\\chi\  even}} \chi(h)  \tau(\bar{\chi})   L(1+\nu,\chi)  = {(2\pi)^{\nu}} \Gamma(-\nu) \cos \left(\frac{\pi \nu}{2}\right)\left\{   \zeta(-\nu,\frac{h}{q})+ \zeta(-\nu,1-\frac{h}{q})\right\}+(1-\frac{1}{q^\nu})\frac{\zeta(1+\nu)}{\phi(q)}.\end{align}
 Inserting \eqref{fan4}, \eqref{fan5} into  \eqref{fan3},  we obtain an expression for the last sum of the right-hand side of \eqref{fan1}.
Then substituting this resulting expression together with \eqref{fan2} in \eqref{fan1}, we get the desired result.\\
  \textbf{Theorem \ref{vor1.3}  $\Rightarrow$ Theorem \ref{vore2}}
  Let $\theta=h/q$, and  $\chi$ be an even primitive non-principal character modulo $q$. Multiplying the identity \eqref{vor1234} in Theorem \ref{vor1.3}
				 by $\bar{\chi}(h)/\tau(\bar{\chi})$, and then summing on $h$, $ 0<h<q$, one can show that Theorem \ref{vor1.3} implies Theorem \ref{vore2}.
      \end{proof}

      \begin{remark}
The other proofs of this section will be almost similar to the proof of the previous theorems. We skip the proofs to avoid repetitions.  
\end{remark}
\section{Acknowledgement}
Research of the second author was supported by the University Grants Commission, Department of Higher Education, Government of India, under NTA Ref. no.  191620205105.  They are grateful to the referee for several helpful suggestions and comments.

		 	\bibliographystyle{IEEEtran}
			\bibliography{bibliography.bib}
		\end{document}